\documentclass[12pt,twoside]{article}
\usepackage[english]{babel}

\usepackage{ifthen}                
\newboolean{boldeqnitalic}
\setboolean{boldeqnitalic}{true}

\usepackage{caption}
\usepackage[intlimits]{amsmath}   
\usepackage{defjs2}                
\usepackage{color}
\usepackage{ifthen}
\usepackage{fancyhdr}
\usepackage{rotating}
\usepackage{multirow}
\usepackage{booktabs}              
\usepackage{amssymb}
\usepackage{amsbsy}
\usepackage{bm}                    
\usepackage{theorem}
\usepackage{upgreek}  
\usepackage{pifont}   
\usepackage{mathrsfs}
\usepackage{datetime}
\usepackage{tikz}
\usetikzlibrary{matrix}
\usepackage[all]{xy}
\usepackage{rotating}
\usepackage{subcaption}

\definecolor{dunkelgrau}{rgb}{0.8,0.8,0.8}
\definecolor{hellgrau}{rgb}{0.90,0.90,0.90} 
\newcommand\defi{\mathrel{\overset{\makebox[0pt]{\mbox{\normalfont\scriptsize\sffamily def}}}{=}}}
\usepackage[numbers]{natbib}  
\fancypagestyle{plain}{%
  \fancyhead{}%
  \fancyfoot[c]{\sffamily\thepage}%
}
\makeatletter                           
\def\cleardoublepage{\clearpage\if@twoside \ifodd\c@page\else
  \hbox{}
  \vspace*{\fill}
  \thispagestyle{empty}
  \newpage
  \if@twocolumn\hbox{}\newpage\fi\fi\fi}
\makeatother
\begin{document}
\unitlength1.0cm
\frenchspacing

\thispagestyle{empty}
 
\vspace{-3mm}
\begin{center}
  {\bf \large The Heterogeneous Multiscale Finite Element Method}\\[2mm]
  {\bf \large for the Homogenization of Linear Elastic Solids} \\[2mm]
  {\bf \large and a Comparison with the FE$^2$ Method} 
\end{center}

\vspace{4mm}
\ce{Bernhard Eidel, Andreas Fischer}

\vspace{4mm}

\ce{\small Heisenberg-Group, Institute of Mechanics, Department Mechanical Engineering}
\ce{\small University Siegen, 57068 Siegen, Paul-Bonatz-Str. 9-11, Germany}
\ce{\small e-mail: bernhard.eidel@uni-siegen.de, phone: +49 271 740 2224, fax: +49 271 740 2436}
\vspace{2mm}

\bigskip

%
%
%
%
%

\bigskip


\begin{center}
{\bf \large Abstract}

\bigskip

{\footnotesize
\begin{minipage}{14.5cm}
\noindent
The Heterogeneous Multiscale Finite Element Method (FE-HMM) is a two-scale FEM based on asymptotic homogenization for solving multiscale partial differential equations.
It was introduced in [W. E and B. Engquist, \emph{Commun. Math. Sci.}, 1 (2003), 87--132]. The objective of the present work is an FE-HMM formulation for the homogenization of linear 
elastic solids in a geometrical linear frame, and doing so, of a vector-valued field problem. 
A key ingredient of FE-HMM is that macrostiffness is estimated by stiffness sampling on heterogeneous microdomains in terms of a modified quadrature formula, 
which implies an equivalence of energy densities of the microscale with the macroscale. Beyond this coincidence with the Hill-Mandel condition, 
which is the cornerstone of the FE$^2$ method, we elaborate a conceptual comparison with the latter method. After developing an algorithmic framework we (i) assess
the existing a priori convergence estimates for the micro- and macro-errors in various norms, (ii) verify optimal strategies in uniform micro-macro mesh refinements 
based on the estimates, (iii) analyze superconvergence properties of FE-HMM, and (iv) compare FE-HMM with FE$^2$ by numerical results. 
\end{minipage}
}
\end{center}

{\bf Keywords:}
Heterogeneous multiscale method; Finite element method; Homogenization; Macro-to-micro modeling; Estimates; Superconvergence \hfill 

\section{Introduction}
\label{sec:intro}

\bigskip

Almost all matter is heterogeneous in its structural composition and so are the mechanical 
properties of solids heterogeneous. As a consequence, deformation and failure of solids 
and structures at large can be critically influenced by the heterogeneity at small scales, 
typically referred to as microscales. Homogenization theories provide a framework for the 
macroscopic modeling of microheterogeneous matter via averaging. Spatial homogenization problems
are typically described by partial differential equations (PDEs). For their solution on arbitrary domains
and for complex boundary conditions numerical methods must be used. Which among many numerical 
homogenization methods is best suited for particular problems is still a controversial issue 
and therefore a most active field of research. Little discussion is about the existing benefits 
and future promises of numerical homogenization; first, homogenization can reduce the computational 
costs or can give at all access to problems, which are prohibitive in full microresolution. 
Second, averaging out small-scale fluctuations of a homogenized solid can provide a better understanding 
of processes at large. Third, and complementary to the 
second aspect, the explicit sampling of properties and analysis of processes on representative 
microdomains allows for predictions on critical phenomena like the initiation of failure of 
structures at large. Hence, in two-scale numerical schemes one can achieve the one, the description
of processes at macroscopic scales without abandoning the other, the added value of small-scale 
processes, and finally, analyze their interactions.

The Heterogeneous Multiscale Finite Element Method, FE-HMM, is a numerical homogenization method for problems described
by PDEs with fast oscillations in the components of the stiffness tensor, the conductivity tensor, or alike.
Mathematical homogenization based on asymptotic expansion provides a sound theoretical basis for FE-HMM. It is 
guiding to a considerable extent the particular formulation of the numerical method in order to preserve 
theoretical characteristics. 

The physical problems treated so far in FE-HMM for spatial \footnote{The most general Heterogeneous Multiscale Method (HMM) introduced in \cite{E-Engquist-2003} also covers problems living on multiple time scales, see \cite{LiE2005}, \cite{E-Engquist-Li-Ren-VandenEijnden2007}.} two-scale problems are diffusion on rough surfaces, \cite{AssyrSchwab2005}, stationary and instationary heat conduction, \cite{AssyrNonnenmacher2009}, Darcy-flow, \cite{Assyr2009}, and transport problems, \cite{AbdulleAttinger2006}. In conclusion, existing FE-HMM implementations and simulations so far have been restricted to scalar-valued field problems.

The progress of FE-HMM and the underlying Heterogeneous Multiscale Method (HMM) since the seminal papers of E, Engquist and coworkers, \cite{E-Engquist-2003}, \cite{E-Engquist-Huang-2003}, is described in the overview papers \cite{E-Engquist-Li-Ren-VandenEijnden2007} and \cite{Assyr-etal2012}. The focus in \cite{Assyr2009} is the fully discrete analysis of multiscale PDEs covering elliptic, parabolic and advection diffusion problems.  
For elliptic problems, a comparison of the FE-HMM with the multiscale finite element method (MsFEM) as introduced in \cite{HouWu1997}, \cite{Hou-Wu-Cai1999} is presented in 

Abdulle presents in  \cite{Assyr2006} the mathematical analysis of FE-HMM for the elliptic case of linear elasticity. The focus is on the fully discrete convergence analysis
taking into account the discretization errors at both micro and macro levels. Simulation results are not considered in that reference. Abdulle \cite{Assyr2006} points out, that no convergence estimates of the overall numerical process taking into account the microscopic and the macroscopic discretization parameters 
have been given for micro-macro methods, even not for the linear case. Many of the cited micro-macro methods can be classified as FE$^2$ methods, which are equally two-level finite element methods for numerical homogenization. The FE$^2$ method has its origins in solid mechanics, \cite{GuedesKikuchi1990}, \cite{MichelMoulinecSuquet1999}, \cite{Miehe-etal-1999a}, \cite{Miehe-etal-1999b}, \cite{Fish-etal-damage1999}, \cite{FeyelChaboche2000}, \cite{Peric-etal-2010}, and has found considerable interest in academia and industry; as a versatile method FE$^2$ has been used in non-linear problems of elasticity and inelasticity. For recent, comprehensive overviews of the FE$^2$ method we refer to \cite{GeersKouznetsovaBrekelmans2010b}, \cite{Schröder2014} and \cite{Saeb-Steinmann-Javili-2016}. In order to account for size-dependency observed in materials science, Kouznetsova et al. \cite{Kouznetsova-etal2001}, \cite{Kouznetsova-etal2002} have introduced a second-order homogenization into FE$^2$. Beyond solid mechanics, FE$^2$ has continuously reached out for various physical branches including multifield problems with a coupling of mechanics with thermo-/electro-/magneto-fields, \cite{OezdemirBrekelmansGeers2008}, \cite{Schröder2009}, \cite{JaviliChatzigeorgiouSteinmann2013}, \cite{KeipSteinmannSchroeder-2014} to name but a few.

The theoretical cornerstone of the FE$^2$ method is the macro homogeneity condition or Hill-Mandel condition, \cite{Hill1963}, \cite{Hill1972}. This condition postulates the equality of macroscopic stress power with corresponding stress power on the microscale. Similarly, FE-HMM is a methodology based on a two-scale energy equivalence since its very first, most general conception in the HMM framework \cite{E-Engquist-2003}. Remarkably, the two methods have not yet been compared in more detail, and to the best of our knowledge, convergence properties for FE$^2$ covering both the microscale and the macroscale error have neither been derived nor have been measured in numerical tests. This is remarkable, since already one decade has past since the work of Abdulle \cite{Assyr2006} in which the author points out this gap.

In view of the addressed issues and open gaps, the present work has the following main aims: 

(i) We derive an FE-HMM formulation for linear elasticity in a geometrical linear frame and address aspects of its implementation. We point to the method's roots in asymptotic homogenization, which are summarized in the appendix for ready reference. Next, we continue with the FE-HMM characteristics of nonstandard numerical quadrature on the macroscale, which then leads to a detailed algorithmic framework. Therein, the novel aspect is the FE-HMM implementation for a vector-valued field problem putting the method into the context of solid mechanics. The setting of linear elasticity considerably eases a thorough comparison with related methods like FE$^2$, see point (ii) below.
In its extension to the vector-valued case of field problems for elliptic PDEs the present work stands on the shoulders of the FE-HMM framework for scalar-valued cases proposed by \cite{AssyrSchwab2005} and \cite{AssyrNonnenmacher2009}, \cite{Nonnenmacher-Disse-2011}. Moreover, the present work considerably takes profit from Abdulle's analysis in \cite{Assyr2006} and aims to complement that reference.
\\[1mm]
(ii) We elaborate a conceptual, algorithmic and numerical comparison of FE-HMM with the FE$^2$ method. The comparison will comprise conceptual aspects like the link of the underlying mathematical theory of asymptotic homogenization in FE-HMM with the physical postulate of energy density equivalence in FE$^2$. It covers the macro-micro kinematical coupling concept and, vice-versa, the micro-to-macro data transfer in terms of stiffness and stress.  
\\[1mm]
(iii) We assess the existing a priori error estimates in various norms ($L^2$-, $H^1$-, energy-norm) and for the coefficients of the homogenized elasticity tensor. We anticipate that the error analysis including a priori estimates are --for our understanding-- the strongest result and most valuable contribution of FE-HMM to the field of numerical homogenization, and the weakest spot of FE$^2$ to go without them. 
While the existence of error estimates in FE-HMM is remarkable, the content of the estimates is even more exciting, since they seem at a first glance to be at odds with familiar estimates from standard finite element methods. Here, we will analyze the inherent property of FE-HMM to exhibit superconvergence properties for the micro-FEM part in that the error in the $L^2$-norm and in the $H^1$-norm exhibit the same convergence order. Superconvergence in standard finite element methods, --if present at all-- refers to the non-standard property of stress and strain to converge pointwise in the same order as displacements as the primary variables. Superconvergence in standard FEM according to Barlow \cite{Barlow1976} requires several premises in that it is restricted to particular element shapes and, additionally, is restricted to particular element sites. For the microscale FEM part of FE-HMM in contrast, superconvergence does not require any premises. An additional benefit of error estimates in FE-HMM is that they allow for optimal uniform macro-/micro mesh-refinement strategies -- how to refine the micromesh for uniform macro mesh refinement, if the full convergence order shall be achieved but for minimal computational costs?

\section{The Heterogeneous Multiscale Finite Element Method}

\subsection{Model Problem of Linear Elasticity} 
\label{subsec:ModelProblemLinearElasticity}

We consider a body $\mathcal{B}$, a bounded subset of $\mathbb{R}^d$, $d=2,3$, with boundary 
$\partial \mathcal{B} = \partial \mathcal{B}_D \cup \partial \mathcal{B}_N$ where the Dirichlet boundary $\partial\mathcal{B}_D$ and the Neumann boundary $\partial\mathcal{B}_N$ are disjoint sets.
The closure of the body $\mathcal{B}$ is denoted by $\overline{\mathcal{B}}$.
The body, which exhibits an inhomogeneous microstructure, is subject to body forces $\bm f$ and surface tractions $\bar{\bm t}$ and in static equilibrium.
 
\subsubsection{The microproblem}
\label{subsec:StrongFormLinearElasticity}

The displacement $\bm u^{\epsilon} = (u_1^{\epsilon}, \ldots, u_d^{\epsilon})$ of the body is given by the solution of the system 
\begin{equation}
\label{eq:StrongFormMicro-2}
\renewcommand{\arraystretch}{1.6}
\left.\begin{array}{rcl}
 - \dfrac{\partial}{\partial x_j} \left(\mathbb{A}^{\epsilon}_{ijlm} \dfrac{\partial u^{\epsilon}_l}{\partial x_m} \right) 
   &=& f_i   \qquad \mbox{in} \quad \mathcal{B}  \\
   u_i^{\epsilon} 
   &=& \bar{u}_i \qquad \mbox{on} \quad \partial \mathcal{B}_{D} \\
   \left(\mathbb{A}^{\epsilon}_{ijlm} \dfrac{\partial u^{\epsilon}_l}{\partial x_m}\right) \, n_j 
   &=& \bar{t}_i   \qquad \mbox{on} \quad \partial \mathcal{B}_{N} \\
\end{array}\right\} \; 
\end{equation}
 The constitutive law is assumed to be linear elastic where $\mathbb{A}^{\epsilon}_{ijlm}$ is the fourth order elasticity tensor.  Superscript $\epsilon$ throughout indicates the dependency of suchlike marked quantities on the heterogeneity of the elastic material. Note that the body forces $\bm f$ and the traction vectors $\bar{\bm t}$ are assumed not to depend on $\epsilon$. In \eqref{eq:StrongFormMicro-2}$_{3}$, $\bm n=(n_1, \ldots, n_d)$ is the unit outward normal to $\partial \mathcal{B}$.  

For $\mathbb{A}^{\epsilon}_{ijlm}$ the following symmetries hold $\mathbb{A}^{\epsilon}_{ijlm}=\mathbb{A}^{\epsilon}_{jilm}=\mathbb{A}^{\epsilon}_{lmij}$ for any $i,j,l,m = 1, \ldots, d$.  
For the deformation kinematics geometrical linearity is assumed to hold with the linearized strain tensor $\bm \varepsilon$
\begin{eqnarray}
 \varepsilon_{ij}(\bm u^{\epsilon}) &=& \dfrac{1}{2} \left(\dfrac{\partial u_i^{\epsilon}}{\partial x_j} + \dfrac{\partial u_j^{\epsilon}}{\partial x_i} \right) \,.
 \label{eq:kinematical-relation}
\end{eqnarray}
It can be written in compact format by means of the linear differential operator $\bm L$ 
\begin{equation}
\label{eq:LinearDifferentialOperator-L}
  \bm \varepsilon (\bm u^{\epsilon}) = \bm L \, \bm u^{\epsilon} \, .
\end{equation}
  
For a finite element formulation the strong form \eqref{eq:StrongFormMicro-2} is transformed into a variational 
or weak form. Multiplying the strong form by a test function $\bm v \in \mathcal{V}$, using the Green formula yields the following variational formulation:

Find $\bm u^{\epsilon}$ such that 
\begin{equation}
\label{eq:weak-form-for-u-epsilon}
 \underbrace{\int_{\mathcal{B}} \mathbb{A}^{\epsilon} (\bm x) \, \bm \varepsilon(\bm u^{\epsilon}) : \bm \varepsilon(\bm v) \, dV}_{\displaystyle \defi B_{\epsilon} (\bm u^{\epsilon}, \bm v)} 
           = \underbrace{\int_{\mathcal{B}} \bm f \cdot \bm v \, dV \, + \,  \int_{\partial \mathcal{B}_N} \bar{\bm t} \cdot \bm v \,dA}_{\displaystyle \defi F(\bm v)} 
\end{equation}
which has to hold for all $\bm v \in \mathcal{V}$, where $\mathcal{V}$ is the space of admissible displacements, i.e. virtual displacements
that fulfill homogeneous Dirichlet boundary conditions
\begin{equation}
\label{eq:HilbertSpaceV}
  \mathcal{V}=\{\bm v; \bm v \in H^1(\mathcal{B})^d, \bm v|_{\partial \mathcal{B}_D} = \bm 0 \} \, .
\end{equation}
The existence and uniqueness of the solution of problem \eqref{eq:weak-form-for-u-epsilon} can be shown 
by use of the first Korn inequality and the Lax-Milgram theorem, \cite{Hughes-BOOK-2000}, \cite{Braess-BOOK-1997}. 

The direct numerical solution of \eqref{eq:weak-form-for-u-epsilon} by a standard finite element formulation is prohibitive for small $\epsilon$, since a proper account of the microheterogeneity of characteristic length $\epsilon$ requires an even finer finite element resolution with typical element size $h$, hence $h\ll \epsilon$. This is the main reason for numerical homogenization based on sampling in small regions of confined size instead of an accurate account of microstructure's heterogeneity everywhere.

The theoretical basis of FE-HMM is mathematical homogenization by asymptotic expansion, \cite{Bensoussan-Lions-Papanicolau-BOOK-1976}, \cite{Sanchez-Palencia-BOOK-1980}, \cite{Allaire1992}, \cite{Cioranescu-Donato-BOOK-1999}.
To put things into perspective and for ready reference we provide the main results of asymptotic homogenization for linearized elasticity in the Appendix \ref{subsec:Homogenization}. These contents are used at various places in the present paper to explain, how FE-HMM is constructed along the lines of asymptotic homogenization. 

\subsubsection{The macroproblem}
\label{subsec:Variational-FE-HMM-macro}

The strong form of the macroscopic/homogenized boundary value problem (BVP) reads 
\begin{equation}
\label{eq:Homogenized-Strong-Form}
\renewcommand{\arraystretch}{1.6}
\left.\begin{array}{rcl}
 - \dfrac{\partial}{\partial x_j} \left(\mathbb{A}^{0}_{ijlm} \dfrac{\partial u^{0}_l}{\partial x_m} \right) &=& \langle f_i \rangle  \qquad \mbox{in} \quad \mathcal{B}  \\
  u_i^{0} &=& \langle \bar{u}_i\rangle_{\Gamma} \qquad \mbox{on} \quad \partial \mathcal{B}_{D}    \\
  \left(\mathbb{A}^{0}_{ijlm} \dfrac{\partial u^{0}_l}{\partial x_m}\right) \, n_j &=& \langle \bar{t}_i\rangle_{\Gamma}   \qquad \mbox{on} \quad \partial \mathcal{B}_{N} \\
\end{array}\right\} \; 
\end{equation}
for a derivation see Sec.~\ref{subsec:Homogenization}. The macroscopic displacement is denoted by $u_{i}^0$ and $\mathbb{A}^{0}$ is the homogenized elasticity tensor. The term in brackets in \eqref{eq:Homogenized-Strong-Form}$_{1,3}$ is the macroscopic stress obtained by a volume average over the microdomain, hence it equals $\langle{\bm \sigma}^0\rangle$, see Sec.~\ref{subsec:Homogenization}. 

The values for the Dirichlet as well as Neumann boundary conditions in \eqref{eq:Homogenized-Strong-Form}$_{2,3}$ are obtained by surface averages according to \eqref{eq:SurfaceAveraging} of corresponding boundary conditions in \eqref{eq:StrongFormMicro-2}$_{2,3}$. Similarly,
$\langle f_i \rangle$ is the volume average of body forces in \eqref{eq:StrongFormMicro-2}$_{1}$.

The solution of the homogenized problem is obtained from the variational form 
\begin{equation}
 B_0 (\bm u^0, \bm v) =  \int_{\mathcal{B}} \mathbb{A}^0 \, \bm \varepsilon(\bm u^0): \bm \varepsilon(\bm v) \, dV   
                      =  \int_{\mathcal{B}} \bm f \cdot \bm v \, dV \, + \,  \int_{\partial \mathcal{B}_N} \bar{\bm t} \cdot \bm v \,  dA
                         \qquad \forall \, \bm v \in \mathcal V \, ,
 \label{eq:VariationalFormHomogenizedProblem}                         
\end{equation}
which follows from multiplying the strong form \eqref{eq:Homogenized-Strong-Form} by test functions $\bm v$ along with the application of Green's formula. For notational convenience we skip in \eqref{eq:VariationalFormHomogenizedProblem} and in the rest of the paper the averaging symbols $\langle \bullet \rangle$, $\langle \bullet \rangle_{\Gamma}$ for $\bm f$, $\bar{\bm u}$ and $\bar{\bm t}$ but keep in mind that these quantities follow from volume and surface averages, respectively. 
\\
There is an alternative route to derive \eqref{eq:VariationalFormHomogenizedProblem}; the asymptotic stress expansion \eqref{eq:AsympExp-sigma} is plugged into he variational form \eqref{eq:weak-form-for-u-epsilon}. The resultant weak form expressions are arranged according to their orders $\mathcal{O}\left(\epsilon^n\right), n=1, -1, -2$. Each of them is multiplied by $\epsilon^n$ with corresponding $n$ along with the limit of $\epsilon \rightarrow 0^+$. Doing so the above macroscale variational form can be identified. Hence, it is the same process as carried out for the identification of the macroscopic balance of linear momentum in Sec.~\ref{subsec:Homogenization}. 
 
Next we consider the piecewise linear continuous FEM in macro- and microspace, respectively. 
The domain $\mathcal{B}$ is a convex polygonal domain in order to avoid regularity issues.

We define a macro finite element space as 
\begin{equation}
 \mathcal{S}^p_{\partial \mathcal{B}_D}(\mathcal{B}, {\mathcal T}_H) = \left\{ \bm u^H \in H^1(\mathcal{B})^d; \bm u^H|_{\partial \mathcal{B}_D} = \bar{\bm u}; \bm u^H|_{K} \in {\mathcal{P}}^{p}(K)^d, \, \forall \, K \in {\cal T}_{H} \right\} \, ,
 \label{eq:MacroFESpace}
\end{equation}
where ${\mathcal P}^{p}$ is the space of (in the present work: linear, $p=1$) polynomials on the element $K$, 
${\mathcal T}_H$ the (quasi-uniform) triangulation of $\mathcal{B} \, \subset \, \mathbb{R}^d$. Index/superscript $H$ denotes the characteristic element size, 
with $H \gg \epsilon$ for efficiency. The space $\mathcal{S}^{p}_{\partial \mathcal{B}_D}$ is a subspace of $\mathcal{V}$ defined in \eqref{eq:HilbertSpaceV}.
 
For the solution of \eqref{eq:StrongFormMicro-2} in the macrodomain
we use the two-scale FEM framework of the FE-HMM as originally proposed in \cite{E-Engquist-2003}
and analyzed for elliptic PDEs in \cite{E-Ming-Zhang-2005}, and, with the focus on linear elasticity, in 
\cite{Assyr2006}.

The macrosolution of the FE-HMM is given by the following variational form:

Find $\bm u^H \in \mathcal{S}_{\mathcal{B}_D}(\mathcal{B}, \mathcal{T}_H)$ such that
\begin{equation}
\label{eq:VariationalFormulationHMM}
 B_H (\bm u^H, \bm v^H) = \int_{\mathcal{B}} \bm f \cdot \bm v^H  \, dV \, + \,  \int_{\partial \mathcal{B}_N} \bar{\bm t} \cdot \bm v^H \, dA
                          \qquad \forall \bm v^H \in \mathcal S_{\partial \mathcal{B}_D} (\mathcal{B}, \mathcal{T}_H) \, ,
\end{equation}
which reads as a standard finite element formulation.  

\subsection{The modified macro bilinear form of FE-HMM}
\label{subsec:Modified-Macro-Bilinear-Form}

If the homogenized constitutive 
tensor $\mathbb{A}^{0}(\bm x)$ is explicitly known, the bilinear form $ B_H (\bm u^H, \bm v^H)$ can be calculated
using standard numerical quadrature according to \eqref{eq:ModifiedBilinearForm-1}, where $\bm x_{K_l}$ and $\omega_{K_l}$ are the quadrature points and 
quadrature weights, respectively
\begin{eqnarray}
 B_H (\bm u^H, \bm v^H) &=& \sum_{K\in \mathcal T_H} \sum_{l=1}^{N_{qp}} \omega_{K_l}\cdot
                            {\color{black}\left[\mathbb{A}^{0}(\bm x_{K_l})\, \bm \varepsilon(\bm u^H(\bm x_{K_l})) : \bm \varepsilon(\bm v^H(\bm x_{K_l})) \right]}
                            \label{eq:ModifiedBilinearForm-1} \\
                        &\approx& \sum_{K\in \mathcal T_H} \sum_{l=1}^{N_{qp}} \omega_{K_l}\cdot \left[ \dfrac{1}{{\color{black}|K_{\delta}|}} 
                            {\color{black} \int_{K_{\delta}} \mathbb{A}^{\epsilon} (\bm x) \, \bm \varepsilon(\bm u^h_{K_{\delta}}) : \bm \varepsilon(\bm v^h_{K_{\delta}}) \, dV} \right] \, .
                            \label{eq:ModifiedBilinearForm-2}                                
\end{eqnarray}

Since $\mathbb{A}^{0}(\bm x)$ is typically not known for heterogeneous materials, the ansatz of FE-HMM is to approximate 
the virtual work expression at point $\bm x_{K_l}$ in the semidiscrete form \eqref{eq:ModifiedBilinearForm-1} by another bilinear form using the known microheterogeneous elasticity tensor $\mathbb{A}^{\epsilon}$, see \eqref{eq:ModifiedBilinearForm-2}. According to this approximation, the solution $\bm u_{K_{\delta}}^h$ is obtained on microsampling domains $K_{\delta}=\bm x_{K_l} + \delta \, [-1/2, +1/2]^d$, $\delta \geq \epsilon$, which are each centered at the quadrature points $\bm x_{K_l}$ of $K$, $l=1, \ldots, N_{qp}$. For a visualization see Fig.~\ref{fig:MicMac-Problem-Meshing-etc}. These microsampling domains with volume $|K_{\delta}|$ provide the additive contribution to the stiffness matrix of the macro finite element.  
\begin{Figure}[htbp]
   \begin{minipage}{16.0cm}  
   \includegraphics[width=15.0cm]{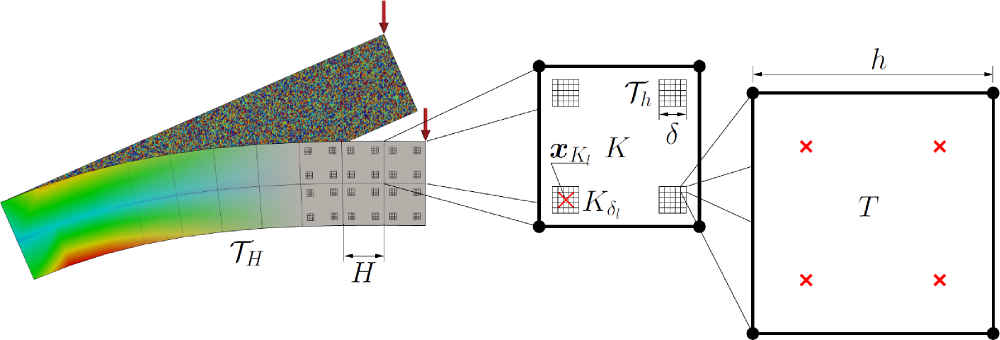}  
\\[2mm]
{\hspace*{3.5cm} (a) \hspace*{4.8cm} (b) \hspace*{3.2cm} (c)}
   \end{minipage}
\caption{Two-scale finite element method: (a) Macroscopic BVP with macrodiscretization, (b) one macro finite element with microdomains centered at the quadrature points, (c) micro finite element with standard quadrature points. 
\label{fig:MicMac-Problem-Meshing-etc}}
\end{Figure}

{\bf Remark 1} \quad
From a physical point of view, the approximation of \eqref{eq:ModifiedBilinearForm-1} by \eqref{eq:ModifiedBilinearForm-2} indicates that the FE-HMM is a numerical homogenization scheme that is based on the equality of the macroenergy density with the microenergy density.   

From a mathematical point of view, the replacement of the pointwise integrand value by another bilinear form indicates that FE-HMM is a modified numerical quadrature formula. Since it is this setting that enables the fully discrete a priori error analysis, it is a key FE-HMM ingredient. The inserted product of test functions leads to a priori estimates of the convergence order of the microerror, which seem to be at odds with familiar results for FEM; for details see Sec.~\ref{subsec:a-priori-error-estimates} along with its references to original work. 

{\bf On Notation} \quad The microdomains attached to the macro quadrature points $l$ in the macro element $K$ are distinguishable by index $l$. In order to avoid too heavy notation we will use $K_{\delta}$ instead of $K_{\delta_l}$ in the 
rest of the paper. 
 
\subsection{Variational formulation of the microproblem}
\label{subsec:Variational-FE-HMM-micro}

It can be shown that the FE-HMM microproblem resembles the discrete version of the cell problem of asymptotic expansion, if it is formulated for each microdomain $K_{\delta}$ in $K$ with $l=1, \ldots, N_{qp}$, $K \in \mathcal{T}_H$ like this:
\\[2mm]
Find $\bm u^h_{K_{\delta}}$ such that the conditions for macro-micro coupling and for the micro bilinear form \eqref{eq:micro-problem-vers1} are fulfilled:
\begin{equation}
\label{eq:micro-problem-vers1}
\renewcommand{\arraystretch}{1.6}
\left.\begin{array}{rcl}
 \left(\bm u^h_{K_{\delta}} -  \bm u^{H}_{lin, K_{\delta}} \right) &\in& \mathcal{S}^q (K_{\delta}, \mathcal{T}_h)  \\[2mm] 
 B_{K_{\delta}}(\bm u^h_{K_{\delta}}, \bm w^h_{K_{\delta}})  &:=& 
  \displaystyle{\int_{K_{\delta}}} \mathbb{A}^{\epsilon} (\bm x) \, \bm \varepsilon(\bm u^h_{K_{\delta}}) : \bm \varepsilon(\bm w^h_{K_{\delta}}) \, dV = 0  \\ 
 & &  \forall  \, \bm w^h_{K_{\delta}} \in \mathcal{S}^q (K_{\delta}, \mathcal{T}_h)   \, ,  
\end{array} \quad \right\} \; 
\end{equation} 
 
where the micro finite element space $\mathcal{S}^q_{}(K_{\delta}, \mathcal{T}_h)$ is defined by 
\begin{equation}
   \mathcal{S}^q (K_{\delta}, \mathcal{T}_h) = \{ \bm w^h \in \mathcal{W}(K_{\delta}); \bm w^h \in (\mathcal{P}^q(T))^d, \, T\in \,\mathcal{T}_h \}\, .
   \label{eq:Periodic-micro-FE-space}
\end{equation}
In \eqref{eq:Periodic-micro-FE-space} $\mathcal{T}_h$ is a quasi-uniform discretization of the sampling domain $K_{\delta}$ with mesh size $h \ll \varepsilon$ resolving the finescale, $\mathcal{W}(K_{\delta})$ denotes the coupling condition or boundary conditions used for computing the microfunctions, and $\mathcal{P}^q$ is the space of polynomials on the element $T$. In the present work we restrict to linear shape functions, hence $q=1$ and consider periodic boundary conditions (PBCs), see \eqref{eq:asymptotic_homo_PBCs_def}.
 
The linearization of $\bm u^H$ in \eqref{eq:micro-problem-vers1}$_1$ is carried out at the quadrature point $\bm x_{K_l}$  
\begin{equation}
 \bm u^{H}_{lin, K_{\delta}} = \bm u^H (\bm x_{K_l}) + (\bm x - \bm x_{K_l}) \cdot \nabla \bm u^H(\bm x_{K_l}) \, .
 \label{eq:linearization_uH}
\end{equation}

{\color{black}}
{\bf Remark 2} \quad Note that \eqref{eq:micro-problem-vers1} resembles for PBC the discrete version of the cell problem of asymptotic homogenization, \eqref{eq:Weak-form-cell-problem}, if $\bm u^h_{K_{\delta}} - \bm u^{H}_{lin, K_{\delta}} = \bm \chi^h_{K_{\delta}}$ and if $\bm u^{H}$ follows the linearization of \eqref{eq:linearization_uH}.
Then, replacing $\bm u^h_{K_{\delta}}$, \eqref{eq:micro-problem-vers1}$_2$ reads as
\begin{equation}
  \label{eq:Weak-form-cell-problem-FEHMM}
  \int_{K_{\delta}} \mathbb{A}^{\epsilon} (\bm x) \, \bm \varepsilon(\bm \chi^h_{K_{\delta}} + \bm u^H_{lin, K_{\delta}} ) : \bm \varepsilon(\bm w^h_{K_{\delta}}) \, dV
  = 0 \qquad \forall  \, \bm w^h_{K_{\delta}} \in \mathcal{S}^q (K_{\delta}, \mathcal{T}_h) \, .
\end{equation}
The linearized $\bm u^{H}_{K_{\delta}}$ is the FE-HMM counterpart of the quantity $\bm I^{lm}$ in asymptotic homogenization; they both induce a {homogeneous} deformation in the microdomain/unit cell and \eqref{eq:Weak-form-cell-problem-FEHMM} coincides with \eqref{eq:Weak-form-cell-problem}. Furthermore, we hint already here at the coincidence of FE-HMM
with the FE$^2$ method with respect to the above superposition of a homogeneous deformation with periodic fluctuations, which is frequently referred to as strain-driven, first order computational homogenization.


\subsection{Bases for spaces}       
\label{subsec:bases-for-spaces}

For the solution of \eqref{eq:micro-problem-vers1} a basis $\{N_I^H\}_{I=1}^{M_{mac}}$ for the macro finite element space $\mathcal{S}^1_0(\mathcal{B}, \mathcal{T}_H)$ is employed in order to represent the macrosolution $\bm u^H$ of \eqref{eq:VariationalFormulationHMM}. Similarly, a basis $\{N_i^h\}_{i=1}^{M_{mic}}$ of the micro finite element space $\mathcal{S}^1_0(K_{\delta}, \mathcal{T}_h)$, \eqref{eq:Periodic-micro-FE-space}, is introduced in order to represent the solution $\bm u^h$ of a microproblem. In \eqref{eq:Macro-displacement-vector-in-the-FE-basis}, $M_{mac}$ denotes the number of nodes of the macrodomain, and $M_{mic}$ in \eqref{eq:micro-displacement-vector-in-the-fe-basis} denotes the number of nodes of each microdomain. 
Hence, the macro- and the microsolution follow the representation
\begin{eqnarray}
\bm u^H &=& \sum_{I=1}^{M_{mac}} N_I^H \, \bm d_I^H\, , \qquad \bm d_I^H= \{d_{I,x_1}^H, d_{I,x_2}^H, d_{I,x_3}^H\}^T \, ,
\label{eq:Macro-displacement-vector-in-the-FE-basis}\\
\bm u^h &=& \sum_{i=1}^{M_{mic}} N_i^h \, \bm d_i^h\, , \qquad \bm d_i^h= \{d_{i,x_1}^h, d_{i,x_2}^h, d_{i,x_3}^h\}^T \, ,
\label{eq:micro-displacement-vector-in-the-fe-basis}
\end{eqnarray}
where $\bm d_I^H$ is the displacement vector of macronode $I$, and $\bm d_i^h$ is the displacement vector for micronode $i$.

The node-based shape functions $\{N_{I/i}^{H/h}\}$ exhibit the standard properties (i) partition of unity: $\sum_{j \in \mathcal{N}} N_j(\bm X_i) \equiv 1 \quad \forall \, \bm X_i \in \mathcal{B}/K_l$ and (ii) compact support $N_j(\bm X_{j^{\prime}}) = \delta_{jj^{\prime}} \quad \forall \, j,j^{\prime} \in \mathcal{N}$, where $\mathcal{N}$ is the set of finite element nodes. 
 

\subsection{Macrostiffness calculation} 
\label{subsec:bottom-up-macrostiffness-calculation}

According to \eqref{eq:ModifiedBilinearForm-2}, FE-HMM can be seen as a bottom-up multiscale method in that the required information for estimating stiffness on the macrolevel is obtained from microsampling domains $K_{\delta}, l=1,\ldots, N_{qp}$.
  
The macro bilinear form $B^e_H(\bm u^H, \bm v^H)$ is the virtual work of internal forces in a macro finite element $e$. The corresponding bilinear form in terms of the shape functions $B^e_H(\bm N_I^H, \bm N_J^H)$ extracts the stiffness matrix contribution $\bm k^{e,mac}_{IJ}$ for macronodes $I, J$, a $d \times d$ matrix. Consequently, we replace in \eqref{eq:ModifiedBilinearForm-1} the displacements $\bm u^H, \bm v^H$ by macro shape function matrices $\bm N_I^H, \bm N_J^H$ and obtain\footnote{Note that matrix $\bm N^H_{I}$ is a $d \times d$ diagonal matrix with elements $N^H_I$; hence for $d=3$ it holds $\bm N^H_I=\mbox{diag}(N^H_I, N^H_I, N^H_I)$. Similarly, matrix $\bm u^{h(I)}$ in \eqref{eq:ModifiedBilinearForm-2-for-varphiH} is a $d \times d$ diagonal matrix diag$(u_x, u_y, u_z)$.} 
\begin{equation}
\bm k^{e,mac}_{IJ} = B_H^{e} (\bm N_I^H, \bm N_J^H) = \sum_{l=1}^{N_{qp}} \dfrac{\omega_{K_l}}{|K_{\delta}|} 
                                        \int_{K_{\delta}} (\bm L \bm u^{h(I)}_{K_{\delta}})^T \mathbb{A}^{\epsilon} (\bm x) \, \bm L \bm u^{h(J)}_{K_{\delta}} \, dV \, .
                                        \label{eq:ModifiedBilinearForm-2-for-varphiH}    
\end{equation}
In \eqref{eq:ModifiedBilinearForm-2-for-varphiH} $\bm u^{h(I)}_{K_{\delta}}$ is the counterpart of $\bm u_{K_{\delta}}^h$ in \eqref{eq:micro-problem-vers1}. It is the dimensionless solution of the microproblem on $K_{\delta}$, which is driven by the shape function $N_I^H$ at macronode $I$. In the following, we add $x_i, i=1, \ldots, d$ to account for the vector-valued field problem of dimension $d$. Consequently, $\bm u^{h(I,x_i)}_{K_{\delta}}$ is the microsolution driven by a macroelement unit-displacement state $\bm u^{H(I,x_i)}_{lin, K_{\delta}}$ at node $I$ in $x_i$-direction.

For stiffness calculation, problem \eqref{eq:micro-problem-vers1} is to be reformulated as follows: Find $\bm u^{h(I,x_i)}_{K_{\delta}}$ on each $K_{\delta}$ such that the conditions for the kinematical coupling and for the micro bilinear form are fulfilled 
\begin{equation}
\label{eq:micro-problem-vers2}
\renewcommand{\arraystretch}{1.6}
\left.\begin{array}{rcl}
 \left(\bm u^{h(I,x_i)}_{K_{\delta}} - {\bm u}^{H(I,x_i)}_{lin, K_{\delta}} \right) &\in& \mathcal{S}^q (K_{\delta}, \mathcal{T}_h)  \\[2mm]
 B_{K_{\delta}}(\bm u^{h(I,x_i)}_{K_{\delta}}, \bm w^{h(I,x_i)}_{K_{\delta}})  &:=&  
  \displaystyle{\int_{K_{\delta}}} (\bm L \bm u^{h(I,x_i)}_{K_{\delta}})^T \mathbb{A}^{\epsilon} (\bm x) \, \bm L \bm w^{h(I,x_i)}_{K_{\delta}} dV = 0  \\ 
 & &  \forall  \, \bm w^{h(I,x_i)}_{K_{\delta}} \in \mathcal{S}^q_{per} (K_{\delta}, \mathcal{T}_h)\,     
\end{array} \quad \right\} \; .
\end{equation} 

In order to realize the coupling of the macrodisplacement field ${\bm u}^{H(I,x_i)}_{lin, K_{\delta}}$ with the yet unknown microdisplacement field 
$\bm u^{h(I,x_i)}_{K_{\delta}}$,
the two fields are expanded into the same basis $\{N_i^h\}_{i=1}^{M_{mic}}$ of $\mathcal{S}^1(K_{\delta}, \mathcal{T}_h)$,
\begin{eqnarray}
   {\bm u}^{H(I,x_i)}_{lin, K_{\delta}} &=& \sum_{m=1}^{M_{mic}} \, N^h_{m, K_{\delta}} {\bm d}^{H(I,x_i)}_{m} \, , 
   \quad {\bm d}^{H(I,x_i)}_{m} = \left({d}^{H(I,x_i)}_{m,x_1}, {d}^{H(I,x_i)}_{m,x_2}, {d}^{H(I,x_i)}_{m,x_3} \right)^T \label{eq:varphi-H-by-psi} \, ,
    \\
   \bm u^{h(I,x_i)}_{K_{\delta}} &=& \sum_{m=1}^{M_{mic}} N^h_{m, K_{\delta}} \, \bm d^{h(I,x_i)}_{m} \, ,
   \quad \bm d^{h(I,x_i)}_{m} = \left( d^{h(I,x_i)}_{m,x_1}, d^{h(I,x_i)}_{m,x_2}, d^{h(I,x_i)}_{m,x_3}\right)^T \label{eq:varphi-h-by-psi} \, .
\end{eqnarray}

The solution of the microproblems for the minimizers $\bm d^{h(I,x_i)}$ is presented in Sec.~\ref{subsec:Solution-of-microproblems}. We continue with the calculation of the macroelement stiffness matrix in \eqref{eq:ModifiedBilinearForm-2-for-varphiH} and obtain
\begin{eqnarray}          
\bm k^{e,mac}_{IJ} 
&=& B^e_H\left[\bm N_I^H, \bm N_J^H\right]  \nonumber\\
&=& \label{eq:k-mac-element-2} 
           \sum_{l=1}^{N_{qp}} \dfrac{\omega_{K_l}}{|K_{\delta}|} 
           \left(\bm d^{h(I,x_i)}\right)^T 
           \, \sum_{T \in \mathcal{T}_h} \,
           \Big(
           \sum_{m=1}^{n_{node}} \sum_{n=1}^{n_{node}} 
                 \, \int_{T} \bm B_m^{e,T} \mathbb{A}^{\epsilon} \, \bm B^e_n \, dV 
           \bm d^{h(J,x_i)}_n  \Big) \\    
          &=& \label{eq:k-mac-element-5}
           \sum_{l=1}^{N_{qp}} \dfrac{\omega_{K_l}}{|K_{\delta}|} 
           \, \left( \bm d^{h(I)} \right)^T \, \bm K^{mic}_{K_{\delta}} \,  \bm d^{h(J)}  \, , \nonumber 
\end{eqnarray}          
where $\bm d^{h(I)}=\left(\,\bm d^{h(I,x_1)} | \bm d^{h(I,x_2)} | \bm d^{h(I,x_3)} \, \right)$ for $d=3$. The assembly of $\bm k^{e,mac}_{IJ}$ yielding $\bm k^{e,mac}$ implies the arrangement of $\bm d^{h(I)}$ in different columns for $I=1,\ldots, N_{node}$, which yields the transformation matrix $\bm T_{K_{\delta}}$. Doing so we obtain 
\begin{eqnarray}
 \bm k^{e,mac}_{K} &=& \sum_{l=1}^{N_{qp}} \dfrac{\omega_{K_l}}{|K_{\delta}|} 
               \, \, \bm T^{T}_{K_{\delta}}\, \bm K^{mic}_{K_{\delta}} \, \bm T_{K_{\delta}} \label{eq:k-mac-element-6}   \\
 \mbox{with} \quad \bm T_{K_{\delta}} &=& \bigg[ \Big[ \big[ \bm d^{h(I,x_i)} \big]_{i=1,\ldots,d} \Big]_{I=1, \ldots, N_{node}} \bigg] \label{eq:k-mac-element-7} \,.        
\end{eqnarray}

Equations \eqref{eq:k-mac-element-6} and \eqref{eq:k-mac-element-7} give insight into important 
characteristics of $\bm T_{K_l}$. First, the matrix dimensions 
\begin{eqnarray}
 \bm k^{e,mac}_K      &\in& \mathbb{R}^{(N_{node} \cdot d)\times (N_{node} \cdot d)}  \nonumber\\
 \bm T_{K_{\delta}}         &\in& \mathbb{R}^{(M_{mic} \cdot d) \times (N_{node} \cdot d)}            \\
 \bm K^{mic}_{K_{\delta}}  &\in& \mathbb{R}^{(M_{mic}\cdot d) \times (M_{mic}\cdot d)}     \nonumber 
\end{eqnarray}
reveal that  $\bm T_{K_{\delta}}$ ''compresses''\footnote{Following the terminology of HMM, $\bm T_{K_{\delta}}$ can be called a stiffness compression operator, see \cite{E-Engquist-2003}, Sec.2.1.} the total micro stiffness matrix $\bm K^{mic}_{K_{\delta}}$ to the dimensions of the macroelement stiffness matrix $\bm k^{e,mac}_K$. For that reason $\bm T_{K_{\delta}}$ is a micro-macro stiffness transfer operator along with a model reduction/coarse-graining. Second, the transformation matrix $\bm T_{K_{\delta}}$ is built up by the column vectors $\bm d^{h(I,x_i)}$, $x_i\,|\,i=1, \ldots, d$, $I=1, \ldots, N_{node}$. 
 
\begin{Figure}[htbp]
   \begin{minipage}{16.0cm}  
    \begin{center}
           \includegraphics[height=6.2cm, angle=0]{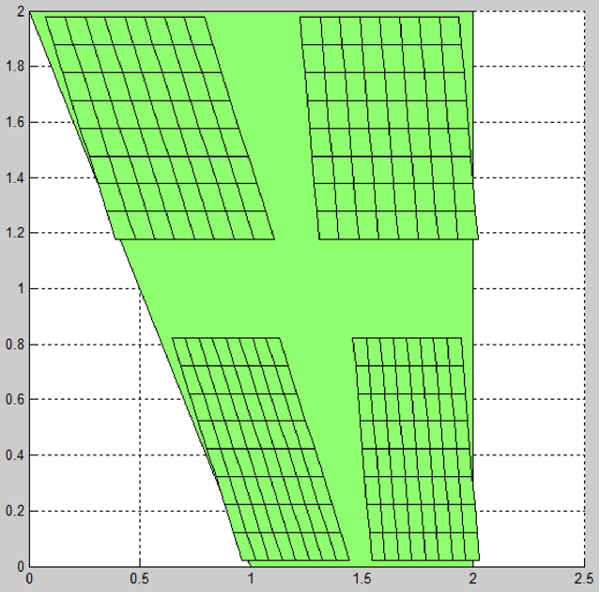}  \hspace*{6mm}
          \includegraphics[height=6.2cm, angle=0]{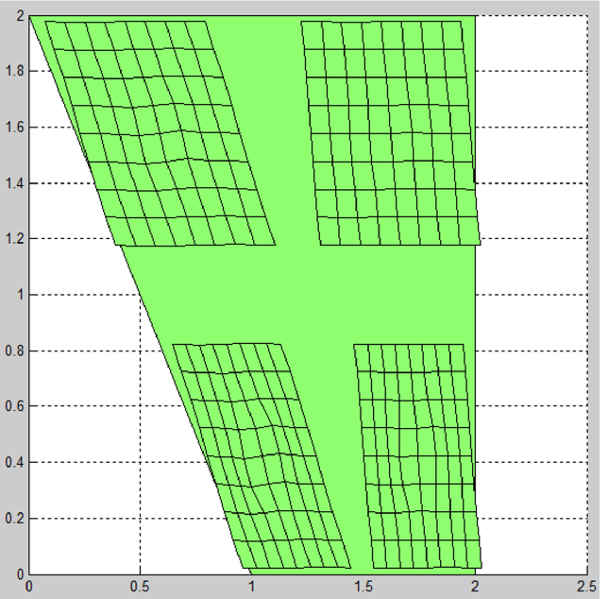}  
    \end{center}
   \end{minipage} 
\caption{Macro-to-micro and periodic kinematical coupling: (left) exemplary unit displacement state at a macronode $I$ in an $x_i$-direction mediated by a macro shape function $N_I$ and the corresponding linearized displacements ${u}^{H(I, x_i)}_{lin, K_{\delta}}$ on the microdomains which is a homogeneous deformation state; (right) the microdomains after energy minimization subject to the additional constraint of PBC result in the microsolutions ${u}^{h(I, x_i)}_{lin, K_{\delta}}$.
 \label{fig:Determination-of-varphi-H}}
\end{Figure}

The macroelement shape functions each represent a unit displacement state for each macro node  $I$ in each direction of space $x_i\, |\, i=1, \ldots, d$.
They drive the microproblem in terms of the corresponding nodal values ${\bm d}^{H(I,x_i)}_{m}, m=1, \ldots, M_{mic}$ via \eqref{eq:varphi-H-by-psi} in each microdomain to evaluate the macroelement stiffness $\bm k^{e,mac}_{IJ}$. Each unit displacement state in $x_i$-direction induces in ${\bm d}^{H(I,x_i)}$ nonzero components only in $x_i$, for example
${\bm d}^{H(I,x_i)}|_{i=2} = \left[ 0, {d}^{H(I,x_2)}_{1,x_2}, 0, \hdots, 0, {d}^{H(I,x_2)}_{M_{mic},x_2}, 0\right]^T$.
  
Figure \ref{fig:Determination-of-varphi-H} (left) visualizes a unit displacement state, the corresponding uniform deformation of microdomains and (right) the microsolution  obtained for the additional constraint of PBC\footnote{The problem deals with a stiff inclusion in a softer matrix, analyzed in detail in Sec.\ref{subsec:matrix-inclusion-problem}.}.  
  
\newcommand\myeq{\mathrel{\overset{\makebox[0pt]{\mbox{\normalfont\sffamily !}}}{=}}}

\subsection{Solution of the microproblems} 
\label{subsec:Solution-of-microproblems}

While the element stiffness matrix of a macro finite element is computed by the modified quadrature formula expanding the collocation 
points to sampling domains, the element stiffness matrix $\bm k^{e,mic}_{mn}$ of a micro finite element is computed by standard 
Gauss-Legendre quadrature in the $\xi$-, $\eta$-, $\zeta$-parameter space on the unit cube (of side length 2) 
\begin{eqnarray}
   \bm k^{e,mic}_{mn} &\approx& \sum_{k=1}^{n_{gp}} \omega_{k} \, \bm B^{e\,T}_m(\xi_k,\eta_k, \zeta_k) \, \mathbb{A}^{\epsilon} \, \bm B^e_n(\xi_k,\eta_k, \zeta_k) \, \, \mbox{det} \, \bm J(\xi_k, \eta_k, \zeta_k) \, ,
                   \label{eq:k-e-mic-4}  
\end{eqnarray}
where $\bm B^{e}_m = \bm L\, \bm N_m$ and where $\bm J(\xi_k, \eta_k, \zeta_k)$ is the Jacobian for an isoparametric finite element formulation. 
 
The assembly of the element microstiffness matrices on the sampling domain $K_{\delta}$ yields the total microstiffness matrix  
\begin{equation}
     \bm K^{mic}_{K_{\delta}} = \Assem \, \bm k_{K_{\delta}}^{e,mic} \, 
     \label{eq:Micro-stiffness-matrix-total}     
\end{equation}
required for the solution of a microproblem.
 
The problem \eqref{eq:micro-problem-vers2} can be rewritten as a minimization problem
\begin{equation}
  \int_{K_{\delta}} (\bm L \bm u^{h(I,x_i)}_{K_{\delta}})^T \mathbb{A}^{\epsilon} \, \bm L \bm u^{h(I,x_i)}_{K_{\delta}} \, dV 
                \quad \rightarrow \quad \mbox{min} \label{eq:micro-p-min} 
\end{equation}
over all functions $\bm u_{K_{\delta}}^{h(I,x_i)} \in \mathcal{S}^1(K_{\delta}, \mathcal{T}_h)$ satisfying $\bm u_{K_{\delta}}^{h(I,x_i)} - {\bm u}_{K_{\delta}}^{H(I,x_i)} \in \mathcal{S}^1_{per}(K_{\delta}, \mathcal{T}_h)$.
   
For the solution of the microproblem \eqref{eq:micro-p-min}, we choose the method of Lagrange multipliers and thus arrive at a saddlepoint problem, see for example \cite{Bertsekas-BOOK-1996}. It is convenient for our purposes to express the strain energy by the micro domain stiffness matrix and the microdisplacement vector. 

Doing so, the problem reads: Minimize for each node on a macro finite element $I=1, \ldots, N_{node}$ and for each direction of space $x_i\,|\,i=1, \ldots, d$ the functional 
\begin{equation}
 \mathcal{L} (\bm d^{h(I,x_i)}, \bm \lambda^{(I,x_i)}) = \dfrac{1}{2} \left(\bm d^{h(I,x_i)}\right)^T \bm{K}_{K_{\delta}}^{mic} \, \bm d^{h(I,x_i)}
 + \bm \lambda^{(I,x_i)\,T} \, \bm G \, \left(\bm d^{h(I,x_i)} - \overline{\bm d}^{H(I,x_i)}\right) \, ,
 \label{eq:Lagrange-functional}
\end{equation}
where $\bm G$ contains the kinematical coupling conditions and reads for $d=3$ 
\begin{equation}
\bm G = 
\left[ \begin{array}{ccccccc}
 b_1 & 0   & 0   & \, \ldots \ldots \, & b_{M_{mic}} & 0 & 0   \\
 0   & b_1 & 0   & \, \ldots \ldots \, & 0 & b_{M_{mic}} & 0   \\
 0   & 0   & b_1 & \, \ldots \ldots \, & 0 & 0 & b_{M_{mic}}   \\[2mm] 
 \multicolumn{7}{c}{\overline{\bm G}}                          \\
\end{array} \right] \, .
\label{eq:D-matrix-3d}
\end{equation}
In $\bm G$, the first $d$ rows contain the normalization condition for the fluctuations, the macro-micro coupling condition, to which the first $d$ Lagrange multipliers are associated; for details see \eqref{eq:normalization-condition-discrete-version}--\eqref{eq:Choice-for-1st-condition-of-coupling}. The submatrix $\overline{\bm G}$ contains the $L$ non-redundant  periodic coupling conditions for adjacent nodes at opposite edges/faces in 2d/3d on the microdomain boundary with the associated $L \cdot d$ Lagrange multipliers. The vector of Lagrange multipliers $\bm \lambda^{(I,x_i)} \in \mathbb R^{(1+L)\cdot d}$ reads for $d=3$ 
\begin{equation}
\bm \lambda^{(I,x_i)}=\{\lambda^{(I,x_i)}_{0,x_1}, \lambda^{(I,x_i)}_{0,x_2}, \lambda^{(I,x_i)}_{0,x_3}, \lambda^{(I,x_i)}_{1,x_1}, \lambda^{(I,x_i)}_{1,x_2}, \lambda^{(I,x_i)}_{1,x_3}, 
\ldots ,\lambda^{(I,x_i)}_{L,x_1}, \lambda^{(I,x_i)}_{L,x_2}, \lambda^{(I,x_i)}_{L,x_3}\}^T \, . 
\end{equation}
The first variations of $\mathcal{L}$ with respect to its variables give the stationarity condition 
\begin{equation}
\left[ \begin{array}{cc}
 \bm K^{mic}_{K_{\delta}} & \bm G^T \\
 \bm G  & \bm 0
\end{array} \right] 
\left[ \begin{array}{c}
\bm d^{h(I,x_i)} \\
\bm \lambda^{(I,x_i)}\\
\end{array}\right]
=
\left[ \begin{array}{c}
\bm 0   \\
\bm G \, {\bm d}^{H(I,x_i)}    \\
\end{array}\right]  \, \, \mbox{for} \, \, I=1, \ldots, N_{node}, \, \, i=1, \ldots, d \, .
\label{Solve4alpha-lambda}
\end{equation}
This set of linear equations is to be solved for $\bm d^{h(I,x_i)}$ and $\bm \lambda^{(I,x_i)}$. Since the coefficient matrix in \eqref{Solve4alpha-lambda} is constant for all macro unit displacement states $(I,x_i)$, the solution vectors as well as the right-hand sides for $I=1, \ldots, N_{node}$ and $x_i, i=1, \ldots, d$ are augmented to full matrices, hence, $\bm d^{h(I,x_i)} \rightarrow \bm T$, $\bm \lambda^{(I,x_i)} \rightarrow \bm \Lambda$, ${\bm d}^{H(I,x_i)} \rightarrow \bm d^H$.   
 
The coefficients $b_i$, $i=1,\ldots, M_{mic}$ in \eqref{eq:D-matrix-3d} follow from the kinematical macro-micro coupling condition, \eqref{eq:micro-problem-vers2}$_1$,  
which can be written as 
\begin{equation}
  \int_{K_{\delta}} \left(\bm u^{h(I,x_i)}_{K_{\delta}} -  {\bm u}^{H(I,x_i)}_{lin, K_{\delta}} \right) \, dV = \bm c \, .
  \label{eq:normalization-condition-discrete-version} 
\end{equation}
Since \eqref{eq:normalization-condition-discrete-version} is the discrete counterpart of the normalization condition \eqref{eq:Normalization-for-chi} for the periodic fluctuations, the particular choice of the constant is inconsequential for the microsolution. Here we choose $\bm c = \bm 0$.

With \eqref{eq:varphi-h-by-psi}, it holds  
\begin{eqnarray}
   \int_{K_{\delta}} \bm u^{h(I,x_i)} \, d V &=& \sum_{T \in \mathcal{T}_h} \sum_{m=1}^{M_{mic}} \bm d_m^{h(I,x_i)}
                       \underbrace{\int_{T} \bm N^h_m \, dV}_{=:b_m} 
                       =  \sum_{T \in \mathcal{T}_h} \bm d^{h(I,x_i)} \cdot \bm b \\
\int_{K_{\delta}} {\bm u}^{H(I,x_i)}_{lin} \, dV &=& \sum_{T \in \mathcal{T}_h} \sum_{m=1}^{M_{mic}} \bm d_m^{H(I,x_i)}
                             \int_{T} \bm N^h_m \, dV 
                       =  \sum_{T \in \mathcal{T}_h} {\bm d}^{H(I,x_i)}  \cdot \bm b \\                       
\mbox{where} \quad  b_m &=& \int_{T} \bm N^h_m \, dV =\dfrac{1}{n_{node}} \, |T| 
                    \label{eq:Choice-for-1st-condition-of-coupling}\\
                        &=& \dfrac{1}{4} \, |T| \quad \mbox{for} \, d=2, \mbox{quadrilaterals and}\, q=1 \, , \label{eq:Choice-for-1st-condition-of-coupling-Result2d} \nonumber \\
                        &=& \dfrac{1}{8} \, |T| \quad \mbox{for} \, d=3, \mbox{hexahedra and}\, q=1 \, . \label{eq:Choice-for-1st-condition-of-coupling-Result3d} \nonumber
\end{eqnarray}

Using the method of Lagrange multipliers, the PBC are not strongly imposed on the FE space but rather weakly enforced. 
In \cite{AssyrSchwab2005} it is found that despite of this weak coupling of scales, the error estimates of \cite{E-Ming-Zhang-2005} 
still hold, where in the latter reference, the PBC are strongly imposed on the FE space.

\subsection{The homogenized elasticity tensor and its efficient computation}
\label{subsec:postproc_elasticity_tensor}

The direct stiffness sampling on microdomains circumvents the necessity to calculate the homogenized elasticity tensor $\mathbb{A}^{0}$ and therefore the necessity of an explicit macroscopic constitutive model. Notwithstanding, calculating $\mathbb{A}^{0}$ in terms of its numerical approximation $\mathbb{A}^{0,h}$ offers a second route to numerical homogenzation. In its derivation we follow \cite{Assyr2006} and propose an efficient way of its computation. 

A reference simplex $\hat{K}$ with its local vectorial basis $\hat{N}^H_{0,j} = (1-\xi_1-...-\xi_d) \otimes e_j , \hat{N}^H_{i,j} = \xi_i \otimes e_j , i,j=1,...,d$ is considered. With $\hat{u}_{i,j}^h$ being the solution to \eqref{subsec:Variational-FE-HMM-micro} such that 
$\hat{u}_{i,j}^h - \hat{N}_{i,j}^H \in \mathcal{S}_{\text{per}}(K_{\delta}, \mathcal{T}_h)$ we can compute the numerical homogenized elasticity tensor as follows
\begin{equation}
\begin{split}
\dfrac{1}{\vert K_{\delta} \vert} \displaystyle
\int_{K_{\delta}}^{} \mathbb{A}^{\epsilon}(x) \bm \varepsilon(\hat{u}_{i,j}^h) \colon \bm \varepsilon(\hat{u}_{l,m}^h) \: d \bm \xi \ & = \ \dfrac{1}{\vert Y \vert} \displaystyle \int_{Y}^{} \mathbb{A}^{0,h} \bm \varepsilon(\hat{N}_{i,j}^H) \colon \bm \varepsilon(\hat{N}_{l,m}^H) \: d \bm y \\
& = \ \mathbb{A}^{0,h} \bm \varepsilon(\hat{N}_{i,j}^H) \colon \bm \varepsilon(\hat{N}_{l,m}^H)\,.
\end{split}
\label{eq_hom_all}
\end{equation}

With 
\begin{equation}
\bm \varepsilon(\hat{N}_{i,j}^H) = 1/2 (e_i \otimes e_j + e_j \otimes e_i), \qquad i,j = 1,..,d \, , 
\label{eq:Macro-unit-strain-state-representation}
\end{equation}
\eqref{eq_hom_all} yields the coefficients of the homogenized elasticity tensor 
\begin{equation}
\mathbb{A}_{ijlm}^{0,h} \ = \ \dfrac{1}{\vert K_{\delta} \vert} \displaystyle
\int_{K_{\delta}}^{} \mathbb{A}^{\epsilon}(x) \bm \varepsilon(\hat{u}_{i,j}^h) \colon \bm \varepsilon(\hat{u}_{l,m}^h) \: d \bm x, \qquad i,j,l,m = 1,...,d  \, . 
 \label{eq:A0-approximation-by-h}
\end{equation}
The error introduced by approximating $\mathbb{A}_{ijlm}^{0}$ by $\mathbb{A}_{ijlm}^{0,h}$ follows
$|\mathbb{A}_{ijlm}^{0,h} - \mathbb{A}_{ijlm}^{0}| \leq C (h/\epsilon)^{2q}$, see \cite{Assyr2009}, Sec. 3.2.2. 

The homogenized elasticity tensor is calculated by inserting selected macro shape functions of a reference element in the variational formulation of the microproblem, which leads to a macroscopic unit strain state as shown in \eqref{eq:Macro-unit-strain-state-representation}. By doing so the homogenized elasticity tensor can be calculated coefficientwise by evaluating the corresponding strains on the microdomain.
 
The homogenized elasticity tensor $\mathbb{A}^{0,h}$ can be calculated alongside the FE-HMM calculation most efficiently from the already available results for the calculation of the macroelement stiffness matrix. For that purpose the microdisplacements are used, which follow from the unit displacement states of a macroelement.

\begin{Figure}[htbp]
 \begin{minipage}{16.0cm}  
\hspace*{2mm}
\includegraphics[height=2.7cm]{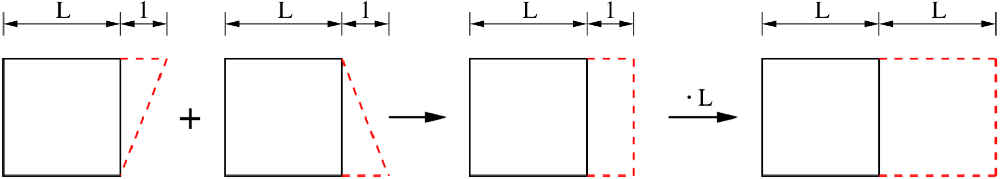}
 \end{minipage}
\caption{Superposition of two unit displacement states and the making of a unit strain state.}
\label{fig:Einheitsverschiebung_1}
\end{Figure}

While the notation above is limited to simplex reference elements, the application of a unit strain state can be easily expanded to other elements. The unit displacement states have to be combined and multiplied by a corresponding factor, such that they yield the corresponding unit strain states on the macrodomain. \\
Figure \ref{fig:Einheitsverschiebung_1} exemplarily visualizes this procedure for a quadrilateral. Two unit displacement states of the nodes on the right macroelement side are superimposed. The resultant constant strain state is multiplied with the element length and results in a unitary strain state.
 
\textbf{Remark 3} \quad The described computation of the homogenized elasticity tensor is only valid for linear elasticity. For nonlinear elasticity, the microstrains have to be calculated from the microdisplacements, which result from the macro unit displacement state.

\subsection{The macrosolution}
  
The global macroquantities, the stiffness matrix $\bm K^{mac}$, the force vector $\bm F$ and the displacement vector $\bm D$ follow from a standard assembly of the corresponding macroelement matrices/vectors 
\begin{equation}
      \label{eq:Assembly-operation}
      \bm K^{mac} = \Assem \, \bm k^{e,mac} \, , \qquad \bm F = \Assem \, \bm f^{e} \, , \qquad  \bm D = \Assem \, \bm d^{e,H}\,. 
\end{equation}

Prescribed displacements $\bm D= \bm D_d$ on Dirichlet boundaries $\partial \mathcal{B}_D$ are eliminated by static condensation 
\begin{equation}
\left[ \begin{array}{cc}
 \bm K_{11} &  \bm K_{12} \\
 \bm K_{21} &  \bm K_{22}
\end{array} \right] 
\left[ \begin{array}{c}
\bm D \\
\bm D_d \\
\end{array}\right]
=
\left[ \begin{array}{c}
\bm F_I   \\
\bm F_d   \\
\end{array}\right] \, .
\label{static-condensation}
\end{equation}
The first set of equations in \eqref{static-condensation} gives the desired macro solution $\bm K_{11}\, \bm D = \bm F_I - \bm K_{12} \bm D_d$.

\subsection{Postprocessing}
\label{subsec:Postprocessing}

After the solution of the macroproblem, the microproblem can be finally solved.
Now, the microproblems in a macro finite element are driven by the true, linearized 
macrodisplacements, $\bm u^H_{lin, K}$, which results in the microsolution vector $\bm u^h$.
The postprocessing is carried out in the following steps
\begin{eqnarray}
   \bm u^h &=& \sum_{i}^{M_{mic}} N_i^h \bm d_i^h \, ,  
     \label{eq:uh-final}
     \\[1mm]
   \bm u^{Hh} (\bm x) &=&  \bm u^H(\bm x_{K_l}) + \bm u^h \qquad \mbox{for} \,\, \bm x \in K \in \mathcal{T}_H \, ,
   \label{eq:MicroscaleDisplacements}
    \\[3mm]
\varepsilon_{ij}(\bm u^H) &=& \dfrac{1}{2} \left(\dfrac{\partial u_i^H}{\partial x_j} + \dfrac{\partial u_j^H}{\partial x_i} \right) \, ,
     \label{eq:Macro-strain-H}
     \\
     \varepsilon_{ij}(\bm u^{Hh}) &=& \dfrac{1}{2} \left(\dfrac{\partial u_i^{Hh}}{\partial x_j} + \dfrac{\partial u_j^{Hh}}{\partial x_i} \right) \, , 
     \label{eq:Micro-strain-Hh}
     \\[3mm]
\bm \sigma(\bm u^{Hh}) &=& \mathbb{A}^{\epsilon} \, \bm \varepsilon(\bm u^{Hh}) = \mathbb{A}^{\epsilon} \, \bm \varepsilon(\bm u^{h}) \, ,
                               \label{eq:FineScaleStressTensor}   
                               \\
\bm \sigma(\bm u^{H}) &=& \dfrac{1}{\vert K_{\epsilon} \vert} \displaystyle \int_{K_{\epsilon}}^{} \bm \sigma(\bm u^{h}) \, dV  
                     \label{eq:HomogenizedStressTensor_1}  
                     \\
                 &=& \mathbb{A}^{0,h} \, \bm \varepsilon(\bm u^{H}) \, . 
                 \label{eq:HomogenizedStressTensor_2}     
\end{eqnarray}
The true resultant microdisplacements are obtained by the superposition of the macrosolution with the microsolution at the centre of the microdomain $\bm x_{K_l}$  
according to \eqref{eq:MicroscaleDisplacements}. From the macrosolution $\bm u^H$ we obtain the macrostrains $\bm \varepsilon(\bm u^H)$, \eqref{eq:Macro-strain-H}, and from the resultant microsolution $\bm u^{Hh}$ the microstrains $\bm \varepsilon(\bm u^{Hh})$, \eqref{eq:Micro-strain-Hh}.

Next, micro stresses $\bm \sigma(\bm u^{Hh})$, \eqref{eq:FineScaleStressTensor}, and the macro stresses $\bm \sigma(\bm u^{H})$, \eqref{eq:HomogenizedStressTensor_1}, are calculated. 

The second equality in \eqref{eq:FineScaleStressTensor} is due to \eqref{eq:MicroscaleDisplacements}. The equality of macrostress according to \eqref{eq:HomogenizedStressTensor_1} with \eqref{eq:HomogenizedStressTensor_2} follows from the fact that macrostress is defined as the volume average of microstress in asymptotic homogenization, and that the homogenized elasticity tensor is identified in this averaging process, see Sec.~\ref{subsec:Homogenization}.

\subsection{Algorithm}

For convenience, we summarize the FE-HMM algorithm for the geometrically linear, linear elastic case in Tab.~\ref{tab:Algorithm}.
Step 1 can be understood as a stiffness preprocessor using already the microsolver. Step 2 is referred to as the macrosolver, 
step 3 uses again the microsolver for the postprocessing.  

\begin{Table}[htbp]
\begin{tabular}{ p{15.5cm}}
\hline\\[-5.05mm]
\hline\\[-4.0mm]
\end{tabular}

\begin{enumerate}
 \item {\bf Macrostiffness calculation by stiffness sampling on microdomains} \\[2mm]
      For all macro finite elements $K \in \mathcal{T}_H$:
      \begin{enumerate}
        \item For each microdomain $K_{\delta}$, $l=1,\ldots, N_{qp}$\,: \\[1mm]
              Compute $\bm K^{mic}_{K_{\delta}}$ \eqref{eq:Micro-stiffness-matrix-total} from $\bm k^{e,mic}_{ij}$ \eqref{eq:k-e-mic-4}. 
        \item For all macronodes $I$ of element $K$, $I=1, \ldots, N_{node}$\,: \\
              \hspace*{6mm} For all directions $x_i, i=1, \ldots, d$\,: \\
              \hspace*{12mm} On each microdomain $K_{\delta}$, $l=1,\ldots, N_{qp}$\,: \\[1mm]
              \hspace*{18mm} -- Calculate ${\bm u}^{H(I,x_i)}_{lin, K_{\delta}}$ and ${\bm d}^{H(I,x_i)}$. \\[1mm]
              \hspace*{18mm} -- Solve for $\bm d^{h(I,x_i)}$, \eqref{Solve4alpha-lambda}, and assemble it in $\bm T_{K_l}$, \eqref{eq:k-mac-element-7}. \\ 
        \item Compute $\bm k^{e,mac}$ from $\bm T_{K_{\delta}}$, $l=1,\ldots, N_{qp}$, and $\bm K^{mic}_{K_{\delta}}$, \eqref{eq:k-mac-element-6}. \\
              (Option:  Calculate $\mathbb{A}^{0,h}$, \eqref{eq:A0-approximation-by-h}).   
       \end{enumerate}
 \item {\bf Macroproblem}     
      \begin{enumerate}
        \item By assembly obtain $\bm K^{mac}$, $\bm F$, and $\bm D$, \eqref{eq:Assembly-operation}.
        \item Solve for $\bm D$, \eqref{static-condensation}, obtain $\bm u^H$.
      \end{enumerate}
 \item {\bf Postprocessing on microdomains and on the macrodomain}   \\[2mm]
        {For all macro finite elements $K \in \mathcal{T}_H$}\, :  
            \begin{enumerate}
                \item For each microdomain $K_{\delta}, l=1,\ldots, N_{qp}$\, : 
                   \begin{enumerate}
                     \item From $\bm u^H$ calculate ${\bm d}^H_{lin}$. 
                     \item Solve for $\bm d^h$, \eqref{Solve4alpha-lambda}.
                     \item Calculate $\bm u^{Hh}$, \eqref{eq:MicroscaleDisplacements}, and $\bm \varepsilon(\bm u^{Hh})$, \eqref{eq:Micro-strain-Hh}.
                     \item Compute $\bm \sigma(\bm u^{h})$, \eqref{eq:FineScaleStressTensor}. 
                   \end{enumerate}
            \item Calculate $\bm \sigma(\bm u^{H})$ by averaging microstresses, \eqref{eq:HomogenizedStressTensor_1}. \\
                  (or, optionally, via $\mathbb{A}^{0,h}$, \eqref{eq:HomogenizedStressTensor_2}). 
            \end{enumerate}
\end{enumerate} 

\begin{tabular}{ p{15.5cm}}
\hline\\[-5.05mm]
\hline\\[-2.0mm]
\end{tabular}
 \caption{Algorithm for FE-HMM in linear elasticity. \label{tab:Algorithm}}
\end{Table}
 
\section{Convergence and a priori error estimates}
\label{subsec:a-priori-error-estimates}
 
For the fully discrete case \cite{Assyr2006} provides convergence results.  
In the periodic case, the macrosolution $\bm u^H$ approximates the homogenized solution $\bm u^0$, see Theorem 4.1 in \cite{Assyr2006}.    
Since $\bm u^0$ does not capture the small scale oscillations of $\bm u^{\epsilon}$, its numerical approximation $\bm u^H$ equally can not.  
But similar to the first order correction according to \eqref{eq:definition_u_1-by-corrector} and the corresponding enriched approximation \eqref{eq:estimate_convergence_u_eps-against-u_01-in-H1-norm}, the reconstructed microsolution $\bm u^{Hh}$ according to \eqref{eq:MicroscaleDisplacements} 
approximates the true fine scale solution $\bm u^{\epsilon}$ in the $H^1$-norm, see Theorem 4.2. in \cite{Assyr2006}.
 
\subsection{A priori estimates} 
\label{subsec:A-priori-estimates}

In view of reliable solution data for engineering decisions, a key property of numerical methods 
is the existence of error estimates with guaranteed bounds.  

FE-HMM exhibits a priori as well as a posteriori estimates for various types of partial differential 
equations. For elliptic PDEs we refer to \cite{E-Ming-Zhang-2005} and \cite{Ohlberger2005}, for the 
(elliptic) case of linear elasticity in a geometrical linear setting to \cite{Assyr2006}.
A posteriori error indicators along with adaptive mesh refinement have been 
presented in \cite{AssyrNonnenmacher2011} and \cite{AssyrBai2013}.

We restrict to a priori estimates for the elliptic problem of linear elastic solids in a geometrical linear frame provided in 
\cite{Assyr2006}. In Sec.~\ref{sec:NumericalExamples} we assess these estimates for uniform discretizations.

The total error of the FE-HMM-method can be decomposed into three parts
\begin{equation}
   || \bm u^0 - \bm u^H || \, \leq \, \underbrace{|| \bm u^0 - \bm u^{0,H} ||}_{\displaystyle e_{mac}} 
                             \, + \, \underbrace{|| \bm u^{0,H} - \widetilde{\bm u}^H ||}_{\displaystyle e_{mod}}
                             \, + \, \underbrace{|| \widetilde{\bm u}^H - \bm u^H ||}_{\displaystyle e_{mic}} \, ,
   \label{eq:Error-decomposition-mac-mod-mic}
\end{equation}
where $e_{mac}$, $e_{mod}$, $e_{mic}$ are the macro error, the modeling error, and the micro error.

Here, $\bm u^0$ is the solution of the homogenized equation, \eqref{eq:Homogenized-Strong-Form}, 
$\bm u^H$ is the FE-HMM solution, $\bm u^{0,H}$ is the standard FEM solution of problem \eqref{eq:VariationalFormHomogenizedProblem}, 
therefore obtained through exact $\mathbb{A}^{0}$; and $\widetilde{\bm u}^H$ is the FE-HMM solution obtained through exact microfunctions 
(in $W(K_{\delta})$).

For sufficiently regular problems, the following a priori estimates hold in the $L^2$-norm, the $H^1$-norm and the energy-norm: 
\begin{eqnarray}
   || \bm u^0 - \bm u^H ||_{L^2(\mathcal{B})} &\leq& C\left( H^{p+1} + \left(\dfrac{h}{\epsilon}\right)^{2q} \right) + e_{mod} \, .
    \label{eq:Total-Error-estimate-L2} \\
   || \bm u^0 - \bm u^H ||_{H^1(\mathcal{B})} &\leq& C\left( H^p + \left(\dfrac{h}{\epsilon}\right)^{2q} \right) + e_{mod} \, ,
   \label{eq:Total-Error-estimate-H1} \\
   || \bm u^0 - \bm u^H ||_{A(\mathcal{B})} &\leq& C\left( H^p + \left(\dfrac{h}{\epsilon}\right)^{2q} \right) + e_{mod} \, ,
   \label{eq:Total-Error-estimate-Energy}  
\end{eqnarray}

{\bf Remark 4} 

(i) Abdulle \cite{Assyr2009} underlines that for the fully discrete analysis of FE-HMM the setting of the modified quadrature formula according to \eqref{eq:ModifiedBilinearForm-2} is essential. The reason is that the error analysis crucially relies on conditions (ellipticity condition and approximation condition, first studied by Ciarlet and Raviart, see \cite{Ciarlet-BOOK-1978}) which ensure that a FEM with numerical quadrature converges with the same rate to the exact solution than the same FEM with exact integration. 
 
(ii) Equations \eqref{eq:Total-Error-estimate-L2} and \eqref{eq:Total-Error-estimate-H1} provide relations of the macro error and the micro error. They enable strategies to achieve the optimal convergence order for minimal computational costs in micro-macro uniform discretizations. Details will be discussed in Sec.~\ref{sec:NumericalExamples} dealing with numerical examples.

\subsection{Superconvergence}
\label{subsec:superconvergence}

The convergence order of $2q$ of the FE-HMM micro error according to \eqref{eq:Total-Error-estimate-L2} is already in the $L^2$-norm at odds with familiar results of standard FEM. The fact however, that the micro error in the $H^1$-norm scales in the same order as in the $L^2$-norm is even more unusual, since it indicates the same order of a derivative as the primary quantity. This phenomenon is called \emph{superconvergence}. We contrast this type of standard FE-HMM superconvergence with the (non-standard) superconvergence in single-scale FEM. For the latter Barlow \cite{Barlow1976} has shown that superconvergence exists only at particular sites of finite elements of rectangular shape and, moreover, that the number and locus of these points depend on the polynomial order of the shape functions. Zienkiewicz and Zhu presented in \cite{SPR}, \cite{SPR2} a procedure for the transfer of the superconvergence property from superconvergent, inner element points to element nodes referred to as ''superconvergent patch recovery'' (SPR). Based on these superconvergent nodal values the same authors constructed a refinement indicator for adaptive remeshing. 
 
Here, we analyze superconvergence for the numerical homogenization by FE-HMM. 

\subsubsection{Macro FEM and micro FEM} The macro error of FE-HMM in the $L^2$-norm (here: of displacements) scales in the order $\mathcal{O}(H^{p+1})$.
The corresponding macro error in the $H^1$-norm and in the energy norm (hence of strain and stress) scale in the order $\mathcal{O}(H^{p})$, one order below the displacements. 
\\
Following the analysis of \cite{Barlow1976}, for $p=1$, superconvergence of strain and stress --hence order $p+1=2$-- can be expected in the element center of macro finite elements, if they exhibit rectangular shapes, which implies that the error in the $H^1$-norm and in the energy norm exhibit the same convergence order as the error in the $L^2$-norm.
\\[2mm]
Since the micro errors of FE-HMM in the $H^1$-norm and in the energy norm exhibit the same convergence order as the error in the $L^2$-norm, $\mathcal{O}((h/\varepsilon)^{2q})$, superconvergence is standard for the micropart of FE-HMM, which therefore is not restricted to rectangular element shapes and not to Barlow's nominally superconvergent points. 

\subsubsection{The superconvergent patch recovery for the macro FEM}

\begin{figure}[htbp]
\centering
\includegraphics[width=5.5cm, angle=0]{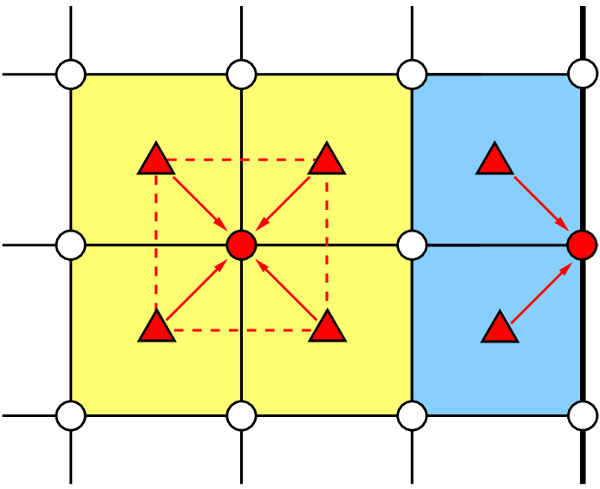} 
\caption{Recovery of nodal stresses from stresses of surrounding superconvergent points (marked by a $\triangle$) here, for $p=1$ in the element center of rectangular quadrilaterals for both a node in the bulk and a node at the surface.}
\label{fig:SPR_Elements}
\end{figure}

In order to verify superconvergence of the macro FEM in the FE-HMM context, the rationale of the SPR according to \cite{SPR} is applied; strain and stress are calculated at superconvergent element sites, for $p=1$ in the center of an element. Next, these values are transferred by a least-square procedure to the finite element node in the direct neighborhood. Elements having such a node in common are referred to as the patch in the superconvergent recovery procedure. For a  visualization see Fig.~\ref{fig:SPR_Elements}. 

Stresses on the patch are prescribed component-wise by
\begin{equation}
\sigma_p^\star = \ \bm P \,\bm a  
\label{eq_ch7_48}
\end{equation}
with 
\begin{equation}
\bm P \ = \ \left[ 1, \ x, \ y, \ xy \right] \quad \text{and} \quad \bm a \ = \ \left[ a_1, \ a_2, \ a_3, \ a_4 \right] \, .
\end{equation}
Vector $\bm P$ contains polynomial terms of the ansatz functions, here for $d=2$ and the case of linear quadrilaterals. For the determination of the 
unknown vector $\bm a$ the function
\begin{align}
F(\bm a) \ &= \ \sum_{i=1}^{n} \, (\sigma_h(x_i,\, y_i) - \sigma_p^\star (x_i, \, y_i))^2 \nonumber \\
&= \ \sum_{i=1}^{n} \, (\sigma_h (x_i, \, y_i) - \bm P (x_i, \, y_i)\mathbf{a} )^2
\end{align}
has to be minimized. Therein, $(x_i, \, y_i)$ are the coordinates of the superconvergent points, $n$ is the number of superconvergent points of the total patch and $\sigma_h (x_i, \, y_i)$ are the stresses in these superconvergent points. Minimization of $F(\bm a)$ implies that $\bm a$ fulfills the condition
\begin{equation}
\sum_{i=1}^{n} \, \bm P^T (x_i, \, y_i) \, \bm P (x_i, \, y_i) \, \bm a \ = \ \sum_{i=1}^{n} \, \bm P^T (x_i, \, y_i) \, \sigma_h (x_i, \, y_i) \, ,
\end{equation}
which can be solved for $\bm a$ 
\begin{equation}
\bm a \ = \ \bm A^{-1} \, \bm b
\end{equation}
with 
\begin{equation}
\bm A \ = \ \sum_{i=1}^{n} \, \bm P^T (x_i, \, y_i) \, \bm P(x_i, \, y_i) \quad \text{and} \quad \bm b \ = \ \sum_{i=1}^{n} \, \bm P^T (x_i, \, y_i) \, \sigma_h (x_i, \, y_i) \, .
\end{equation}
Stresses in the central node of the patch can be recovered by inserting its nodal coordinates $(x_N, \, y_N)$ into the $\bm P$-Vector in \eqref{eq_ch7_48}.
\begin{figure}[htbp]
\centering
\includegraphics[width=0.65\linewidth]{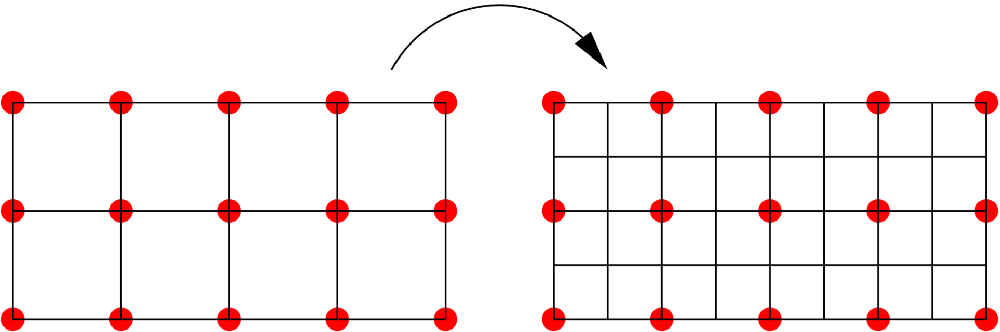}
\\[6mm]

\caption{Exemplary hierarchical refinement of a finite element mesh preserving the initial nodal points.}
\label{fig:Verfeinerung}
\end{figure}
Summarizing, the numerical assessment of the superconvergence property of stresses and strains can be calculated as described in Tab.~\ref{tab:AlgorithmSuperconvergence}

\begin{Table}[htbp]
\begin{tabular}{ p{15.5cm}}
\hline\\[-5.05mm]
\hline\\[-4.0mm]
\end{tabular}

\begin{enumerate}
	\item Computation of the solution of a BVP for various meshes, which must be hierarchical to ensure a proper comparison of stress and strain  at identical nodes, for an example see Fig.~\ref{fig:Verfeinerung}.
	\item Evaluation of stresses and strains in the superconvergent points and recovery of nodal values by means of the superconvergent patch recovery procedure.
	\item Selection of the solution obtained on the finest mesh as the ''overkill''-solution.
	\item Pointwise calculation of the convergence order of the nodal error in the energy-norm.
\end{enumerate}

\begin{tabular}{ p{15.5cm}}
\hline\\[-5.05mm]
\hline\\[-2.0mm]
\end{tabular}
 \caption{Algorithm for assessing superconvergence. \label{tab:AlgorithmSuperconvergence}}
\end{Table}
 

\section{Conceptual Comparison of FE-HMM with the FE$^2$ Method}
\label{sec:FE-HMM-vs-FEsquared}

The aim of the present section is a conceptual comparison of the FE-HMM with the FE$^2$ method, 
which reveals similarities, but also differences in some noticable aspects. 

\subsection{The FE$^2$-Method}
\label{subsec:FEsquared}

To put things into perspective we briefly reiterate the building blocks of the FE$^2$ method.
 
The basic idea of computational homogenization in the FE$^2$ method is to solve at each integration point of the macroproblem the associated microproblem. For the case of the first-order strain-driven computational homogenization the macroproblem provides input to the problem at the microscale in terms of macroscopic deformation and, vice-versa, the microproblem provides stiffness information in terms of the macroscopic tangent and stress. The latter is obtained from averaging microstresses, which implies that a constitutive law does not exist on the macroscale.

The theoretical cornerstone of the FE$^2$ method is the Hill-Mandel condition, which postulates  
that the stress power on the macrolevel must be equal to the average microscopic stress power in a so-called Representative Volume Element, RVE 
\begin{equation}
    \overline{\bm \sigma} : \dot{\overline{\bm \varepsilon}} = \dfrac{1}{V} \int_{\text{RVE}} \bm \sigma : \dot{\bm \varepsilon} \, dV \, 
    \qquad      \longleftrightarrow   \qquad 
    \dfrac{1}{V} \int_{\text{RVE}} \bm \sigma : \dot{\bm \varepsilon} \, dV - \overline{\bm \sigma} : \dot{\overline{\bm \varepsilon}} = 0 \, ,
    \label{eq:Hill-Mandel-condition1}
\end{equation}
 where macrostress $\overline{\bm \sigma}$ is calculated as the volumetric mean of microstresses $\bm \sigma$ according to 
\begin{equation}
    \overline{\bm \sigma} = \dfrac{1}{V} \int_{\text{RVE}} \bm \sigma \, dV \, .
   \label{eq:FEsquared-macrostress-definition}
\end{equation}

Boundary conditions for the $\text{RVE}$, which are consistent to the Hill-Mandel postulate, are the constraint condition $\dot{\bm \varepsilon}:= \dot{\overline{\bm \varepsilon}}$ on the whole $\text{RVE}$ (Voigt condition), the constraint condition of constant stress $\bm \sigma = \overline{\bm \sigma}$ on the 
whole $\text{RVE}$ (Reuss condition). Moreover, linear Dirichlet and linear Neumann conditions each fulfill the condition as well. Finally, PBC are consistent 
with the postulate. They are visualized in Fig.~\ref{fig:PBCs}. Macrostrain $\overline{\bm \varepsilon}$ yields a homogeneous deformation on the RVE, which is superimposed by periodic fluctutations $\widetilde{\bm w}$. 
\begin{Figure}[htbp]
   \begin{minipage}{16.0cm}  
    \begin{center}
         \includegraphics[width=11.0cm, angle=0]{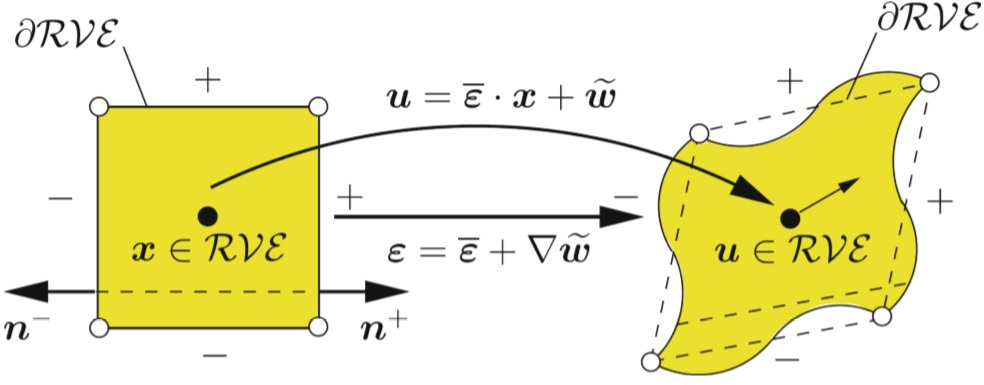}    
    \end{center}
   \end{minipage} 
   \\[8mm]
\caption{Periodic boundary conditions (PBC) on the boundary of the microsampling domain $\partial \text{RVE}$.
\label{fig:PBCs}}
\end{Figure}
For a proper definition of PBC, the boundary of the $\text{RVE}$ is divided into pairwise periodic parts $\partial \text{RVE}^+$ and $\partial \text{RVE}^-$ with a corresponding sign convention for position vectors $\bm x^+$ and $\bm x^-$ and normal vectors $\bm n^+ = - \bm n^-$. For pairwise periodicity of the displacements, i.e. $\widetilde{\bm w}^+ = \widetilde{\bm w}^- = \widetilde{\bm w}$ along with  
\begin{equation}
    \int_{\text{RVE}} \widetilde{\bm \sigma} : \dot{\widetilde{\bm \varepsilon}} \, dV 
    =
    \int_{\partial\text{RVE}} \widetilde{\bm t} : \dot{\widetilde{\bm w}} \, dA 
    =  
    \int_{\partial\text{RVE}} (\bm t - \overline{\bm \sigma} \bm n) \cdot (\dot{\bm u} - \dot{\overline{\bm \varepsilon}} \bm x) \, dA  = 0\, ,
   \label{eq:LinearBCs-Requirement}
    \end{equation}
yields
\begin{equation}
    \int_{\partial\text{RVE}} \widetilde{\bm t} \cdot \dot{\widetilde{\bm w}} \, dA 
    =  
    \int_{\partial \text{RVE}^+} \widetilde{\bm t}^+ \cdot \dot{\widetilde{\bm w}} \, dA 
    + 
    \int_{\partial \text{RVE}^-} \widetilde{\bm t}^- \cdot \dot{\widetilde{\bm w}} \, dA 
    =
    \int_{\partial\text{RVE}^-} (\widetilde{\bm t}^+ + \widetilde{\bm t}^-) \cdot \dot{\widetilde{\bm w}} \, dA  = 0\, ,
   \label{eq:Periodic-BCs}
\end{equation}
 
which implies the additional condition $\widetilde{\bm t}^+ = - \widetilde{\bm t}^-$. 
Thus, the PBC satisfying the macrohomogeneity condition can be given as
\begin{equation}
 \widetilde{\bm w}^+ = \widetilde{\bm w}^- \qquad \mbox{and} \qquad \widetilde{\bm t}^+ = - \widetilde{\bm t}^- \qquad \mbox{on} \quad \partial \text{RVE} \, .
\end{equation}
 
The homogenized elasticity tensor $\overline{\mathbb{A}}$ is calculated for each macro quadrature point.
If static condensation is applied, $\overline{\mathbb{A}}$ can be calculated according to 
\begin{equation}
 \dfrac{d\, \overline{\bm \sigma}}{d\, \overline{\bm \varepsilon}} =: \overline{\mathbb{A}} = \mathbb{A}^{\text{Voigt}} - \dfrac{1}{|K_l|} \overline{\bm L}^T (\bm K_{K_{l}}^{mic})^{-1} \overline{\bm L} \, ,
 \label{eq:FEsquare-homogenized-elasticity-tensor}
\end{equation}
with
\begin{equation}
 \mathbb{A}^{\text{Voigt}} = \langle \mathbb{A} \rangle = \dfrac{1}{{|K_l|}} \sum_{T \in \mathcal{T}_h} \int_T \mathbb{A}^{\varepsilon} \, dV \, , 
 \qquad \overline{\bm L} = \Assem \, \overline{\bm l_e}\, ,  \quad
 \overline{\bm l_e}  = \int_T \bm B^T \mathbb{A}^{\varepsilon} \, dV \, .
 \label{eq:FEsquare-A_voigt_L_Le-definitions}
\end{equation}
The corresponding system of linear equations which has to be solved to obtain $\overline{\mathbb{A}}$  is
\begin{equation}
  \bm K_{K_{l}}^{\text{mic}} \, \bm X = \overline{\bm L} \, .
\end{equation}
The number of right hand sides follows from the dimension of $\overline{\bm L}$, hence depends on ${\bm B}$ acccording to \eqref{eq:FEsquare-A_voigt_L_Le-definitions}$_{2,3}$. The number is 3 for $d=2$, and 6 for $d=3$.

\subsection{Comparison}
\label{subsec:Comparison-FEHMM-FE2}
\subsubsection{Homogeneity condition. Energy equivalence} 
While the FE$^2$-method is built on the Hill-Mandel condition \eqref{eq:Hill-Mandel-condition1}, FE-HMM is built on asymptotic expansion.
It is a straightforward excercise to show that for PBC, asymptotic homogenization theory fulfills the Hill-Mandel condition.
With the strain influence function already introduced in \eqref{eq:Elastic-strain-influence-function} a direct calculation of microstrain from macrostrain is obtained
\begin{equation}
   \varepsilon_{ij} = E_{ij}^{kl} (\bm y) \, \overline{\varepsilon}_{kl} (\bm x) = \left( \mathbb{I}_{ijkl} + \dfrac{\partial \chi_i^{kl}}{\partial y_j}\right) \overline{\varepsilon}_{kl} \, .
   \label{eq:ElasticStrainInfluenceFunction}
\end{equation}
With \eqref{eq:ElasticStrainInfluenceFunction} we obtain
\begin{eqnarray}
      \dfrac{1}{V} \int_{V} \sigma_{ij} \, \varepsilon_{ij} \,dV 
               &=& \dfrac{1}{V} \int_{V} \sigma_{ij} \left( \mathbb{I}_{ijkl} + \dfrac{\partial \chi_i^{kl}}{\partial y_j}\right) \,dV \overline{\varepsilon}_{kl} 
                   \label{eq:HillMandel-fromMathHomog-1} \nonumber\\
               &=& \underbrace{\dfrac{1}{V} \int_{V} \sigma_{ij} \,dV}_{\displaystyle =\, \overline{\sigma}_{ij}} \, \overline{\varepsilon}_{ij} 
                   + \dfrac{1}{V} \underbrace{\int_{V} \sigma_{ij} \dfrac{\partial \chi_i^{kl}}{\partial y_j} \,dV}_{\displaystyle =\,0} \, \overline{\varepsilon}_{kl} \label{eq:HillMandel-fromMathHomog-2}\\
               &=& \overline{\sigma}_{ij} \, \overline{\varepsilon}_{ij} \, . \label{eq:HillMandel-fromMathHomog-3} \nonumber
\end{eqnarray}
 
For the evaluation of the second integral in \eqref{eq:HillMandel-fromMathHomog-2} we have used integration by parts

\begin{equation}       
          \int_{V} \sigma_{ij} \dfrac{\partial \chi_i^{kl}}{\partial y_j} \,dV \, \overline{\varepsilon}_{kl} =
          \int_{\partial V} \chi_i^{kl} \sigma_{ij} n_j \, d\gamma \, \overline{\varepsilon}_{kl} - \int_{V}  \chi_i^{kl} \dfrac{\partial \sigma_{ij}}{\partial y_j} \,dV \, \overline{\varepsilon}_{kl} 
          \label{eq:HillMandel-fromMathHomog-Integration-by-parts} = 0 
\end{equation}      
with $\bm n$ being the unit normal to the boundary. 
Here we exploited periodicity, namely that $\chi_i^{kl}$ is equal on the opposite site of the unit cell, whereas surface tractions $t_i= \sigma_{ij} n_j$ are equal and of opposite sign on opposite sites of the unit cell. The second integral in \eqref{eq:HillMandel-fromMathHomog-Integration-by-parts} vanishes, since it holds ${\partial \sigma_{ij}}/{\partial y_j} = 0$ for equilibrium. 

Equivalence of energy densities in FE-HMM is enforced in the frame of a modified quadrature rule for the calculation of the macrostiffness matrix by stiffness sampling on microdomains, see~\eqref{eq:ModifiedBilinearForm-2}. Already in its very first, most general proposal, the Heterogeneous Multiscale Method (HMM) for variational problems was built on energy equivalence in terms of a fine-scale to coarse-scale stiffness transfer, see Sec. 2.2 in \cite{E-Engquist-2003}. For that reason we can call the FE-HMM energy equivalence ''E-macrohomogeneity condition'' for the initials of energy and the authors. 
 
\subsubsection{Kinematics of coupling conditions} 
Both FE-HMM and FE$^2$ attach the microdomains to the macrodomain at the quadrature points of macro finite elements as visualized in Fig.~\ref{fig:MicMac-Problem-Meshing-etc}. The explicit kinematical coupling follow for FE$^2$ from the energy equivalence condition \eqref{eq:Hill-Mandel-condition1}, and for FE-HMM from the cell problem of asymptotic homogenization. Periodic boundary conditions meet these requirements such that the macro-micro kinematical coupling effectively coincides for FE$^2$ and FE-HMM; in FE$^2$ strain at the macro quadrature point $\overline{\bm \varepsilon}$ induces a homogeneous deformation in the attached RVE with the corresponding displacement field $\bm u= \overline{\bm \varepsilon} \cdot \bm X$. In FE-HMM the homogeneous deformation in the microdomain is induced by the linearized macrodisplacement field $\bm u^H_{lin}$, which is superimposed by periodic fluctuations following from constrained energy minimization in the RVE and microdomain, respectively. In the above nomenclature, the periodic fluctuations are referred to as $\widetilde{\bm w}$ in FE$^2$, in FE-HMM as the difference $\bm u^h - \bm u^{H}_{lin}$, and in asymptotic homogenization as $\bm \chi^{lm}$. 

The enforcement of coupling conditions have been realized in the FE$^2$ method either by static condensation/elimination, \cite{MichelMoulinecSuquet1999}, by the method of Lagrange multipliers, \cite{Miehe-Koch-2002}, or by the Penalty method, \cite{TemizerWriggers2008}. In FE-HMM, the coupling was realized by Lagrange multipliers in \cite{AssyrSchwab2005}, \cite{AssyrNonnenmacher2009}, \cite{EidelFischer2016}. For a discussion of the role of different types of boundary conditions on the RVE we refer to \cite{Yue-E2007}, \cite{LarssonRunesson-etal-2011}.

\subsubsection{Micro-to-macro coupling/bottom-up data transfer} 
\label{subsubsec:Mic-mac-coupling}
In strain-driven two-level finite element frameworks, the microscale must pass over stress and stiffness information to the macroscale. 
For the particular case of linear elasticity, the deformation-dependent tangential material stiffness boils down to the constant elasticity tensor. 

Stiffness transfer in FE$^2$ is carried out by the effective fourth-order constitutive tensor $\overline{\mathbb{A}}$, although its computation may differ. It can be carried out by numerical perturbation as proposed in \cite{FeyelChaboche2000}, by Lagrange multiplier techniques, \cite{Miehe-Koch-2002}, or by static condensation, \cite{MichelMoulinecSuquet1999}. The effective constitutive tensor is passed over to the corresponding macro quadrature point, where it enters the macro element stiffness matrix. 

For FE-HMM in practice, the bottom-up stiffness transfer can be carried out by two different means, either by (i) a direct stiffness matrix transfer mediated by micro-to-macro transformation matrices or by (ii) the homogenized constitutive (here: elasticity) tensor, where the latter coincides with the FE$^2$ method. We anticipate that the two variants of stiffness transfer show quantative agreement, see Sec.~\ref{sec:NumericalExamples}. Notwithstanding, it is only the first approach of the stiffness matrix transfer which enables due to its underlying modified quadrature formula the derivation of unified a priori estimates and therefore is a key, indispensible FE-HMM ingredient, see Sec. \ref{subsubsec:Apriori-estimates}.

\subsubsection{Data-driven versus model-driven method} 
\label{subsubsec:Data-model-driven-methods}
Closely related to the bottom-up data transfer is the discrimination between data- and model-driven methods. Within the FE$^2$ method, a macroscopic constitutive model is not defined since not required. Instead, the macroscale is fed by homogenized microdata in terms of stress and stiffness, see Sec.\ref{subsubsec:Mic-mac-coupling}. In this respect FE$^2$ is a data-driven method, in view of the absent macro constitutive equations it can even be called an ''equation-free'' method. 

Similarly, the sole constitutive input for FE-HMM is $\mathbb{A}^{\epsilon}$, which lives on the microscale. It is passed over to the macroscale via the transformed micro stiffness matrix, and macrostresses are calculated by averaging microstresses. Of course, the macroscale solver is a finite element method but the physical model is not declared.

\subsubsection{A priori error estimates} 
\label{subsubsec:Apriori-estimates}
FE-HMM is endowed with a priori estimates for the macro error and the micro error in the $L^1$-, $H^2$- and the energy-norm as well as with estimates for the coefficients of the homogenized elasticity tensor, see Sec.~\ref{subsec:A-priori-estimates}. 
 
While the a priori estimates of the macro error in FE-HMM directly follow from the approximation condition for standard FEM, \cite{Ciarlet-BOOK-1978}, the a priori estimate for the micro error and its propagation to the macroscale crucially relies on the FE-HMM setting as a modified quadrature rule and on the method's roots in asymptotic homogenization, see \cite{Assyr2005}, \cite{Assyr2006}. FE$^2$ lacks these roots and, as a consequence, no unified a priori estimate covering both the macro and the micro error exists. Moreover, there are, to the best of our knowledge, even no numerical convergence tests available for FE$^2$.  

Unified a priori estimates allow in FE-HMM for micro-macro uniform refinement strategies, which are optimal in that the full convergence order can be achieved for a minimum of computational costs. This is a considerable benefit of FE-HMM compared with FE$^2$. In numerical tests related to scalar-valued field problems it was shown that the a priori estimates are quite sharp, \cite{Assyr2009}, \cite{AssyrNonnenmacher2009}. For the present case of linear elasticity, the a priori estimates and derived optimal micro-macro uniform refinement strategies will be assessed in Sec.~\ref{sec:NumericalExamples}.

\subsubsection{Numerical costs}
\label{subsubsec:Numerical-costs}
The method of Lagrange multipliers for stiffness calculation in FE-HMM requires to solve the microproblem for $N_{node} \times d$ right hand sides. Since the microstiffness matrix is independent of the unit displacement states, it is a linear set of equations with $N_{node} \times d$ different right hand sides thus reducing the computational complexity. In FE$^2$ using the static condensation method, the number of right hand sides in the microproblem is 3 for $d=2$, and 6 for $d=3$, hence cheaper than the FE-HMM approach using Lagrange multipliers. Of course, the solution of the microproblems in FE$^2$ by the method of Lagrange multipliers is equally more expensive than static condensation, see \cite{Miehe-Koch-2002}. Moreover, if static condensation is used in FE-HMM, the numerical costs are the same as in FE$^2$. If FE-HMM takes the route using the homogenized
elasticity tensor $\mathbb{A}^{0,h}$, the numerical effort is exactly the same as in FE$^2$.
 
In either case, linear problems, here for both geometrical linearity and material linearity, are solved not only by FE-HMM but also by FE$^2$ in one single step, no iterations are required, see~\cite{Assyr2006}, Sec.~1; for that reason the algorithm in Tab.~\ref{tab:Algorithm} similarly applies to the FE$^2$ method. 

Beyond parallelization several approaches can reduce the computational complexity of two-level FEM. Methods based on Fast-Fourier Transforms (FFT) as introduced in \cite{Moulinec-Suquet-1994}  use a direct, point-wise discretization of the Lippmann-Schwinger equation, with improvements in \cite{MichelMoulinecSuquet1999} and \cite{Schneider-Ospald-Kabel-2015}. In an effort to reduce the computation time by a modified modeling \cite{SchröderBalzaniBrands2010} introduce into FE$^2$ the rationale of a statistically similar representative volume element (SSRVE) in order to replace the true microstructure in its full geometrical complexity by a surrogate, which resembles the original one by geometrical features as analyzed by different geometrical similarity measures.

\subsubsection{Limitations}
\label{subsubsec:Limitations}
The classical FE$^2$ method has proven to be equally suitable for linear and non-linear problems. It has been applied to non-linear elasticity, to various inelastic constitutive models, and to multifield problems with a coupling of mechanics with thermo-/electro-/magneto-fields. The versatility of FE$^2$ is its big plus; the extension of FE-HMM to finite deformation non-linear elasticity is a non-trivial task, most notably to provide the mathematical analysis in terms of estimates for that case. 

In the frame of first-order computational homogenization methods FE$^2$ and FE-HMM exhibit inherent limitations. Inhomogeneous deformation modes such as bending cannot be properly represented on the microscale. As a consequence, the microscale is not appropriate to capture geometrical size effects in material behavior. This issue was overcome by second order computational homogenization as introduced into FE$^2$ by \cite{Kouznetsova-etal2001}, \cite{Kouznetsova-etal2002}. 
For a thorough discussion of the limitations and future challenges we refer to \cite{GeersKouznetsovaBrekelmans2010a}. Alternative to the second-order homogenization to capture the size effect, one can introduce surface/interface energies into homogenization which are more intuitive and can be physically
motivated as described in \cite{Duan-Wang-Huang-Karihaloo2005} and \cite{JaviliSteinmannMosler2017}.

\subsubsection{Relation to other methods}
\label{subsubsec:Relation-to-other-methods}
There is a remarkable connection of FE$^2$ and FE-HMM to the fully nonlocal, cluster-based Quasicontinuum Method (CQC), which was introduced by \cite{KnapOrtiz2001}. In CQC, nanosampling domains in the shape of spherical clusters are used to approximate the sum over all lattice sites of a crystal by a weighted sub-sum. For a modification to an energy-based approach of this concept endowing the method with a variational structure see \cite{EidelStukowski2009}, \cite{Eidel2009}. The nanosampling domains of CQC are strongly coupled to the displacement field of the finite elements. Additionally, the clusters are coupled by the nonlocal interactions of the atoms in the sampling domain with other atoms within the cutoff-radius. The concept of scale-separation is a necessary theoretical requirement in FE$^2$ and FE-HMM. Quite in contrast, CQC as an atomistic-continuum coupling method aims at a seamless transition between the scales. For CQC, the sole constitutive input originates from the fine scale (potential of atomic interactions), and consequently, no explicit coarse-scale constitutive model exists. It is remarkable that in this particular characteristics, CQC is close to FE$^2$ and FE-HMM, but is in contrast to other atomistic-continuum coupling methods, which equally require a coarse-scale/macro model, \cite{Eidel-etal-2010}.

In conclusion, innovative quadrature formula (FE-HMM) and summation rules (CQC), see \cite{AmelangVenturiniKochmann2015}, have been advancing multiscale modeling to a considerable extent. 
 
\section{Numerical Examples}   
\label{sec:NumericalExamples} 
This section assesses the numerical performance of the FE-HMM method. The results and corresponding errors
are measured in the following norms
\begin{eqnarray}
\mbox{$L^2$-norm:} \quad  || \bm u ||_{L^2(\Omega)} &:=& \sqrt{\int_{\Omega} \bm u : \bm u \, dV } \, ,
\label{eq:L2-norm} \\
\mbox{$H^1$-norm:} \quad  || \bm u ||_{H^1(\Omega)} &:=&  
             \sqrt{ \left( \sum_{i,j=1}^d \int_{\Omega} \left( \dfrac{\partial u_i}{\partial x_j} \right)^2 \, dV
                 +  \sum_{i=1}^d \int_{\Omega} \left(u_i\right)^2 \, dV \right) }\, ,
\label{eq:Hilbert-norm} \\
\mbox{energy-norm:} \quad  ||\,\bm u\,||_{A(\Omega)} &:=& \sqrt{\int_{\Omega} \mathbb{A} \, \bm \varepsilon(\bm u) : \bm \varepsilon(\bm u) \, dV} 
                                             = \sqrt{\bm d^T \, \bm K \, \bm d }\, ,
\label{eq:energy-norm} \\
\mbox{maximum-norm:} \quad  ||\,\bm u\,||_{\infty} &:=& \mbox{sup}_{\bm x \in \Omega} | \bm u| \, .
\label{eq:max-norm}  
\end{eqnarray}

Four benchmark problems subject to plane strain conditions will be analyzed, 
\begin{enumerate}
 \item[(i)] a matrix-inclusion problem, \\[-6mm]
 \item[(ii)] a microstructure with an analytical solution for homogenization,\\[-6mm]
 \item[(iii)] a non-uniformly periodic microstructure,\\[-6mm]
 \item[(iv)] a uniformly periodic tessellation inspired by M.C. Escher.
\end{enumerate}

Key aspects are 
\begin{enumerate}
 \item[(A)] the convergence order of the simulation results, where the errors are calculated by means of an accurate reference solution (''overkill''-solution) on very fine meshes. Doing so, the a priori estimates will be assessed.
 \item[(B)] Based on the unified error estimates covering the macro error and the micro error, tests will be carried out on how full convergence orders can be achieved for minimal costs in uniform micro-macro mesh refinements.
 \item[(C)] Superconvergence employing the superconvergent patch recovery will be tested.
 \item[(D)] Simulation results of FE-HMM are compared with those of FE$^2$, thus complementing the conceptual comparison of Sec.~\ref{subsec:Comparison-FEHMM-FE2}.
 \item[(E)] Finally, the reconstruction of the microsolution from the macrosolution shall demonstrate the capability of postprocessing small scale features, which can be used to investigate local phenomena like stress concentrations, which eventually may initiate inelastic deformations and failure mechanisms.
\end{enumerate}

To improve the readability but still properly document the present simulation results, we provide the numerical convergence data in tabular form as a supplement and restrict here to the more telling convergence diagrams.  
 
The elasticity tensor will be given in standard Voigt notation. Moreover, we replace the $x_i$, $i=1,\ldots,d$ coordinate system by a $x$-,$y$-,$z$-coordinate system.

\subsection{Matrix-inclusion problem for beam-bending} 
\label{subsec:matrix-inclusion-problem}

\begin{Figure}[htbp]
   \begin{minipage}{16.0cm}  
   \centering
\includegraphics[height=5.5cm, angle=0]{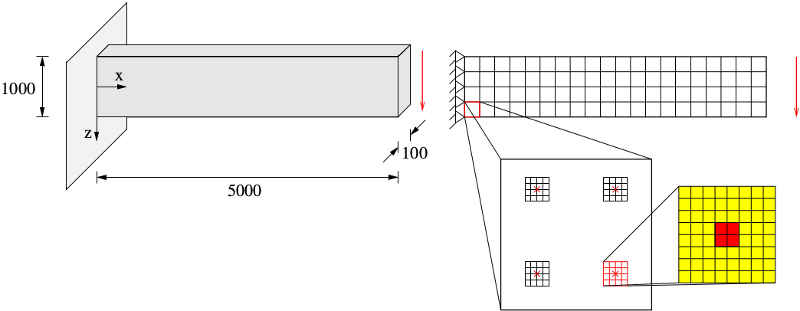}    
   \end{minipage}
\caption{Matrix-inclusion material for beam bending. Geometry, boundary conditions, loading and finite element discretizations.
\label{fig:MatrixInclusion2d-Geometry}}
\end{Figure}

In the first numerical example we consider the microstructure of a stiff inclusion in a soft matrix, hence
a problem with uniformly periodic coefficients, $\delta = \epsilon$, such that periodic boundary 
conditions are applied to the microdomains.  The Young's moduli of the inclusion $E_i=100000$ $\text{N/mm}^2$ 
and the matrix phase $E_m=40000$ $\text{N/mm}^2$ exhibit the contrast of $E_i/E_m=2.5$, for the Poisson's ratio it holds
$\nu=0.2$. The volume ratio of the inclusion phase is $V_i/V_{tot}=1/16$.  

The beam as displayed in Fig.~\ref{fig:MatrixInclusion2d-Geometry} exhibits length $l_1=5000~\text{mm}$, height $l_2=1000~\text{mm}$, 
and thickness $t=100~\text{mm}$. The external line load is $f=1~\text{N/mm}^2$. The side length of the square microcell is 
$\epsilon=5~\text{mm}$, hence $\epsilon/l_1=5/5000 \ll 1$ and $\epsilon/l_2=5/1000 \ll 1$. 

\subsubsection{Macro- and microconvergence} 
\begin{Figure}[htbp]
   \begin{minipage}{16cm}  
      \centering
     \includegraphics[width=6.0cm, angle=0]{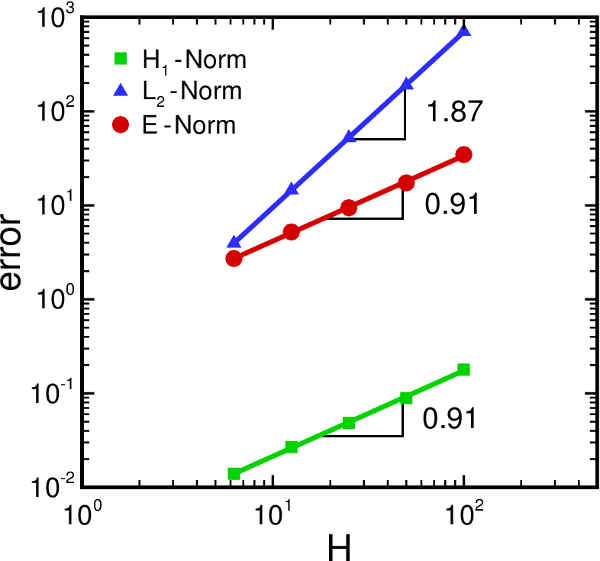} \hspace*{4mm}
     \includegraphics[width=6.0cm, angle=0]{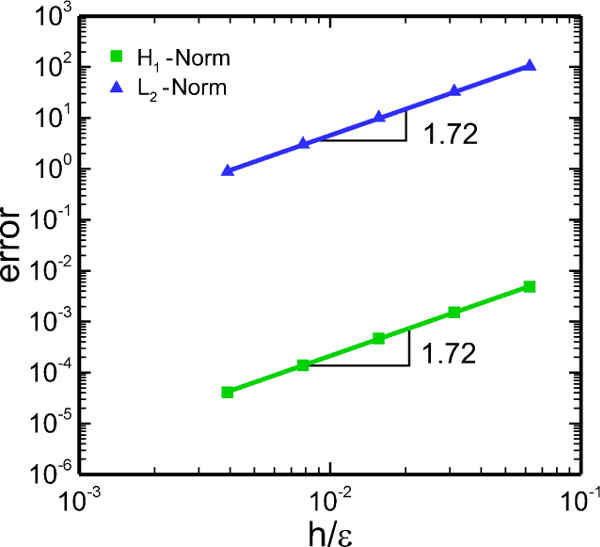}  
   \end{minipage}
\caption{Matrix-inclusion problem for beam bending. (Left:) Macro- and (right:) microconvergence.
\label{fig:MatrixInclusion2d-Macro-Micro-Convergence}}
\end{Figure}

The simulation results 
are obtained along with a reference solution with 3200 $\times$ 640 
macroelements and $32 \times 32$ microelements. The results turn out to be invariant with respect to a variation of $\epsilon$.
Figure \ref{fig:MatrixInclusion2d-Macro-Micro-Convergence} (left) shows the total error as a function of the macrodiscretization while keeping the microdiscretization fixed; the mean convergence orders in the $L^2$-norm of almost 2, and in the $H^1$- and in the energy-norm of almost 1 confirm the a priori estimates of Sec.~\ref{subsec:A-priori-estimates}.
 
Similarly, the convergence of the micro error  
in right diagram of Fig.~\ref{fig:MatrixInclusion2d-Macro-Micro-Convergence} is in reasonable agreement with the theoretical predictions. The deviation (1.72 versus 2) is due to the magnitude of the stiffness jump at the inclusion-matrix interface, which lowers the regularity of the microproblem. This assertion is underpinned by the local convergence distribution of the error in the energy-norm as shown in the right of Fig.~\ref{fig:MatrixInclusion2d-A0,h-Convergence}; at the matrix-inclusion interphase, most notably  at the corners, the convergence order is considerably reduced. The reference solution for microerror calculation is obtained using $200 \times 40$ macroelements and $2048 \times 2048$ microelements. 
 
\subsubsection{Homogenized elasticity tensor} 
\begin{Figure}[htbp]
   \begin{minipage}{16cm}  
      \centering
     \includegraphics[width=6.0cm, angle=0]{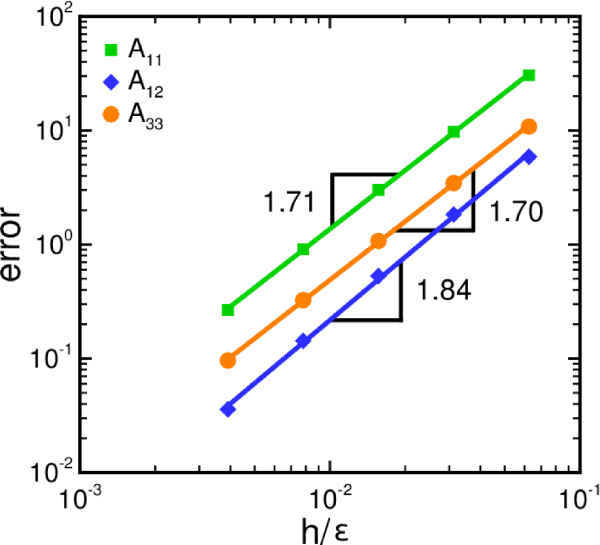} \hspace*{10mm}  
     \includegraphics[width=6.0cm, angle=0]{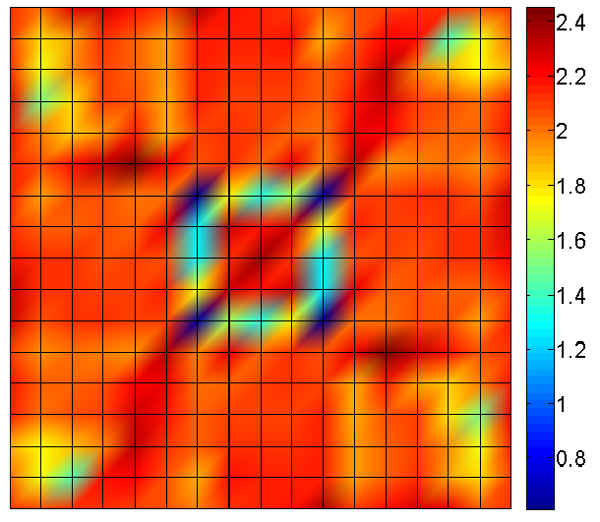}
   \end{minipage}
\caption{Matrix-inclusion problem for beam bending. (Left) convergence of $\mathbb{A}^{0,h}$ components, (right) local convergence order of the error in the energy-norm on the microdomain at [4980.3 mm, 105.3 mm].
\label{fig:MatrixInclusion2d-A0,h-Convergence}}
\end{Figure}

Next, we assess the convergence estimate for the components of the approximate homogenized elasticity tensor $\mathbb{A}_{}^{0,h}$ where the reference solution $\mathbb{A}_{}^{0}$ is obtained on a microdomain discretized by $2048 \times 2048$ elements. The convergence as a function of $h/\epsilon$ is shown in the left diagram of Fig.~\ref{fig:MatrixInclusion2d-A0,h-Convergence}. The mean slope varies between 1.70 for coefficient $\mathbb{A}_{33}$ and 1.84 for coefficient $\mathbb{A}_{12}$ and well agrees with the convergence of the microerrors in the $L^2$-norm and the $H^1$-norm. The deviation from the theoretical order 2 is caused by the aforementioned reduced regularity of the microproblem for its material stiffness contrast. 

\subsubsection{Comparison of FE-HMM with FE$^2$} 

In practice, the use of $\mathbb{A}^{0,h}$ in the FE-HMM simulation is an alternative route to the standard FE-HMM solution method. The difference between the solutions is quantitatively analyzed in Tab.~\ref{tab:MatrixInclusion2d-FEHMM-vs-FE2}; the normed percental deviation is throughout less than $5\times 10^{-10}$ and independent of the macrodiscretization.

Since the FE$^2$ method draws in the homogenization process on the fourth order constitutive tensor, here $\mathbb{A}^{0,h}$, it coincides with the FE-HMM byproduct of the homogenized elasticity tensor. The simulation results for FE-HMM along with the micro-macro stiffness matrix transfer show quantitative agreement with the FE$^2$, which underpin the identity of the methods by numerical means as well.
 
\begin{table}[htbp]
\center
\renewcommand{\arraystretch}{1.2}
\begin{tabular}{lrrrrr}
\hline \\[-5mm]         
& microelements                             &  \multicolumn{4}{c}{$64 \times 64$ }  \\
& macroelements                             &  $50 \times 10$  & $100 \times 20$  & $200 \times 40$  &  $400 \times 80$   \\[1mm]
\hline \\[-5mm]                       
FE-HMM   &  $|| \bm u ||_{\infty} $         &   11.8018   &  11.8498   &  11.8630  &  11.8667   \\
         &  $|| \bm u ||_A$                 &   1080.36   &  1082.50   &  1083.07  &  1083.22   \\[1mm]
\hline \\[-5mm] 
FE$^2$ $\mathrel{\widehat{=}}$  &  $|| \bm u ||_{\infty}$  &   11.8018  &  11.8498  &  11.8630  & 11.8667  \\
FE-HMM($\mathbb{A}^{0,h}$)   &  deviation               &    4.09 $\cdot 10^{-11}$   &  -6.72  $\cdot 10^{-13}$ &  4.68 $\cdot 10^{-13}$  &  3.58 $\cdot 10^{-11}$       \\ 
                             &  $|| \bm u ||_A$         &   1080.36  &   1082.50 &  1083.07  & 1083.22  \\
                             &  deviation               &  2.18 $\cdot 10^{-11}$  &   -3.38 $\cdot 10^{-11}$ & 2.35 $\cdot 10^{-10}$  &  1.96 $\cdot 10^{-11}$       \\[1mm]
\hline
\end{tabular}
\caption{{Matrix-inclusion problem for beam bending.} Comparison of the standard FE-HMM with the FE-HMM($\mathbb{A}^{0,h}$)/FE$^2$ solution. 
\label{tab:MatrixInclusion2d-FEHMM-vs-FE2}}
\end{table}

\subsubsection{Postprocessing on selected microdomains} 
\begin{Figure}[htbp]
   \begin{minipage}{16cm}  
   \centering
     \includegraphics[height=4.2cm, angle=0]{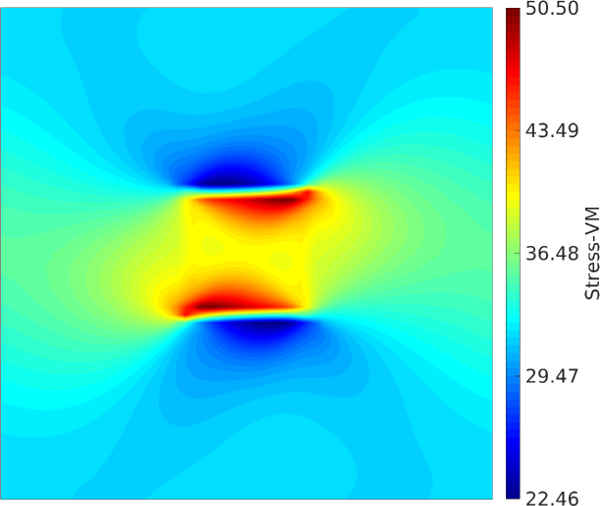} \hspace*{12mm}
     \includegraphics[height=4.2cm, angle=0]{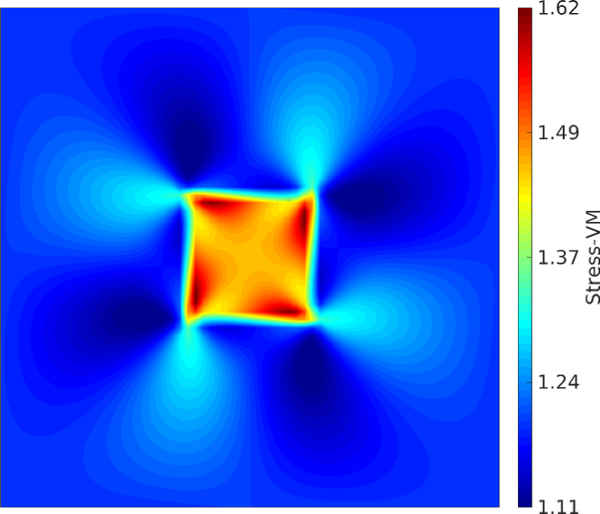}  
   \end{minipage}
\caption{Matrix-inclusion problem for beam bending. Von-Mises stress distribution (left) on the microdomain at $x=21.1,y=21.1$ (mm)
and (right) on the microdomain at $x=4978.9,y=978.9$ (mm).
\label{fig:MatrixInclusion2d-VonMisesStress-micro-Macro-micro}}
\end{Figure}
Figure \ref{fig:MatrixInclusion2d-VonMisesStress-micro-Macro-micro} displays the von-Mises stress distributions on two different
microdomains, in the lower left and in the upper right of the beam. 
The stiff inclusion attracts stress, while the stress in the soft matrix is considerably smaller. Right at the interface of the inclusion with the matrix, a stress jump can be observed; the fact that the stiffer inclusion carries over large stresses {shielding} the softer matrix in terms of the matrix' strong local reduction of stress is referred to as stress-shielding.
\begin{Figure}[htbp]
\begin{minipage}{16.5cm}  
	\centering
	{\includegraphics[height=4.0cm]{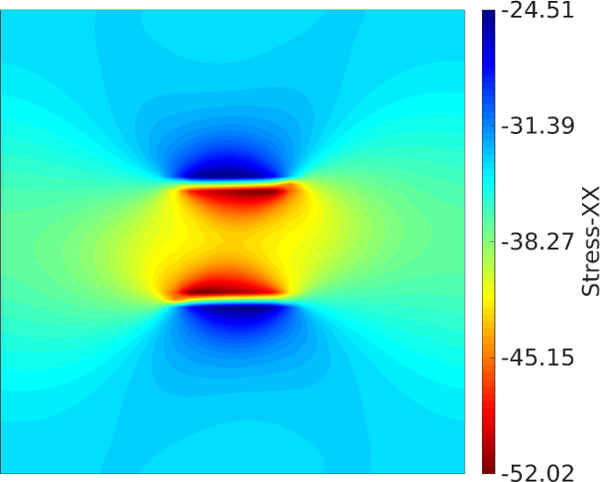}} \hspace*{1mm}
	{\includegraphics[height=4.0cm]{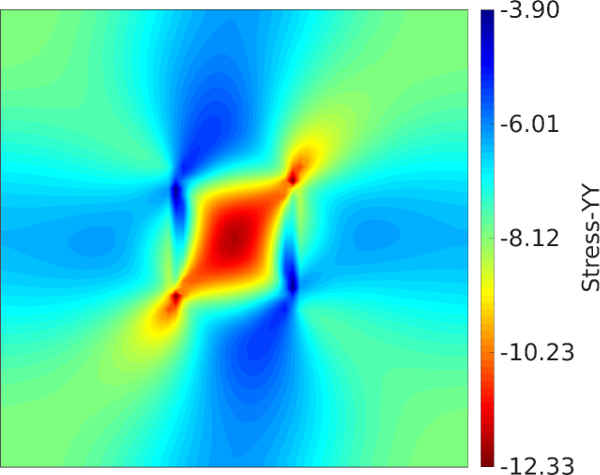}} \hspace*{1mm}
	{\includegraphics[height=4.0cm]{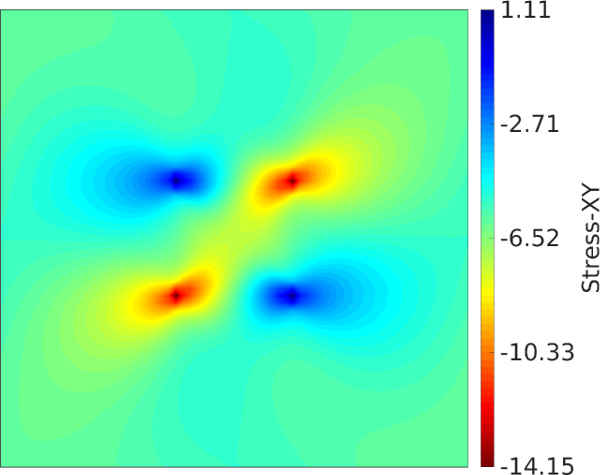}} \\[4mm] 
	{\includegraphics[height=4.0cm]{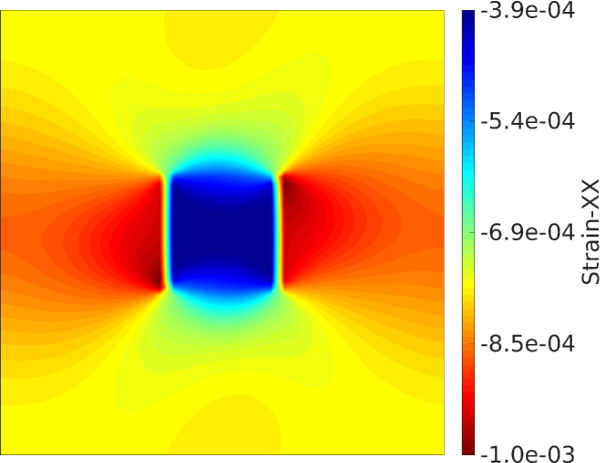}} \hspace*{1mm}
	{\includegraphics[height=4.0cm]{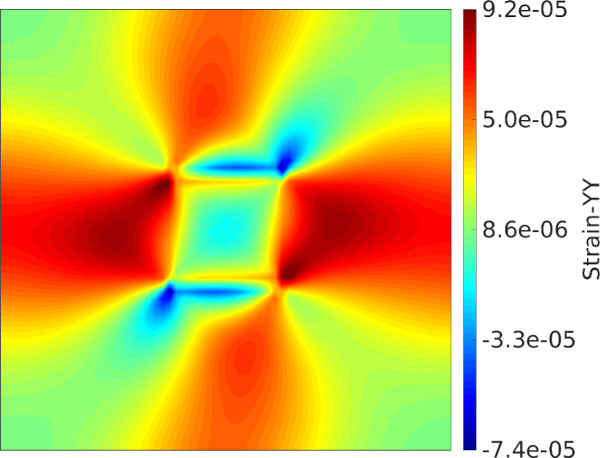}} \hspace*{1mm}
        {\includegraphics[height=4.0cm]{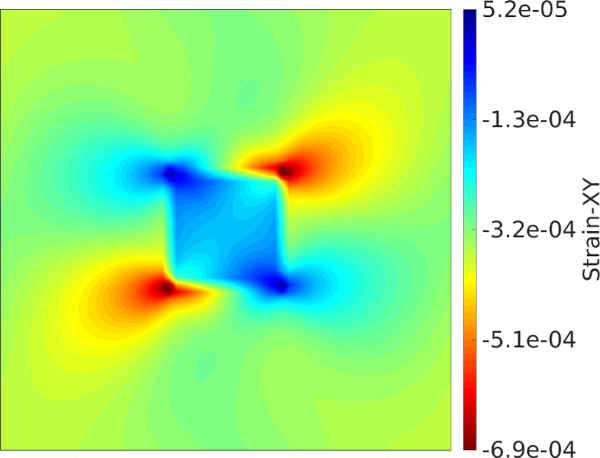}}
\end{minipage}        
\caption{Matrix-inclusion problem for beam bending. On the microdomain for the macro Gauss point at [$21.1 \, \text{mm}$, $21.1 \, \text{mm}$] (first row) stress $\sigma_{xx}$, $\sigma_{yy}$ and $ \tau_{xy}$ in $[\text{N/mm}^2]$ and (bottom row) strain $\varepsilon_{xx}$, $\varepsilon_{yy}$ and $\gamma_{xy}$ are displayed.}
	\label{fig:MatIncl_StressLL_StrainLL}
\end{Figure}  

Refined information is obtained by the contour plots of relevant stress components and strain components for the two points considered in the lower left of the beam, Fig.~\ref{fig:MatIncl_StressLL_StrainLL}, and for the upper right of the beam, Fig.~\ref{fig:MatIncl_StressRU_StrainRU}.
\begin{Figure}[htbp]
\begin{minipage}{16.5cm}  
	\centering
	 {\includegraphics[height=4.0cm]{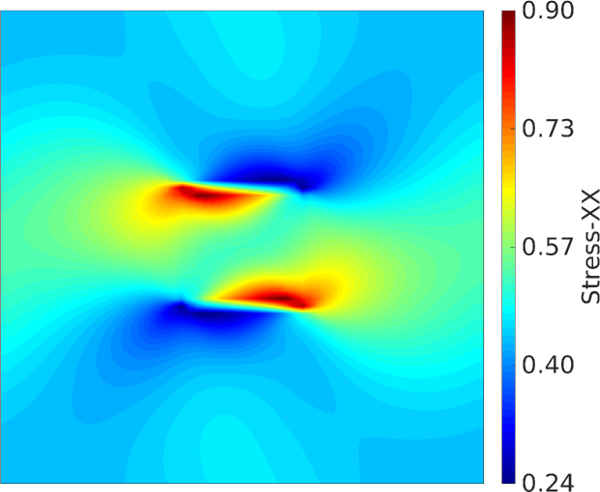}} \hspace*{1mm}
	 {\includegraphics[height=4.0cm]{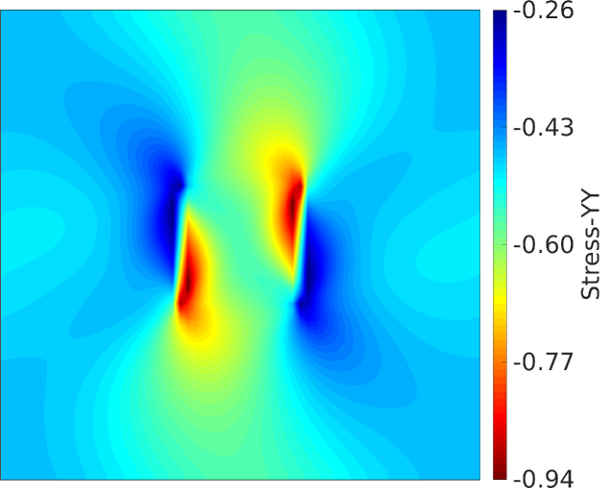}} \hspace*{1mm}
	 {\includegraphics[height=4.0cm]{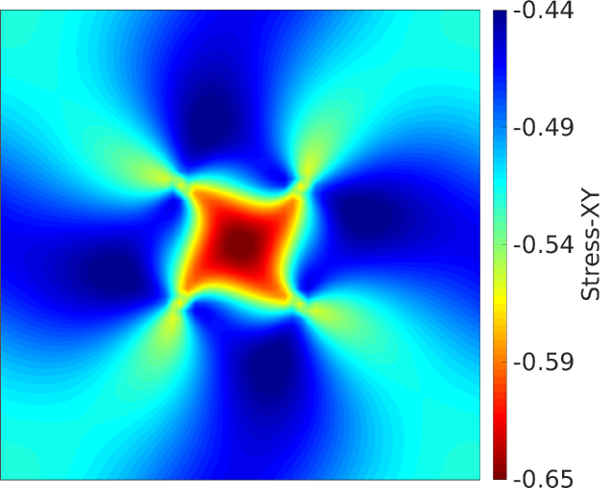}} \\[4mm] 
	 {\includegraphics[height=4.0cm]{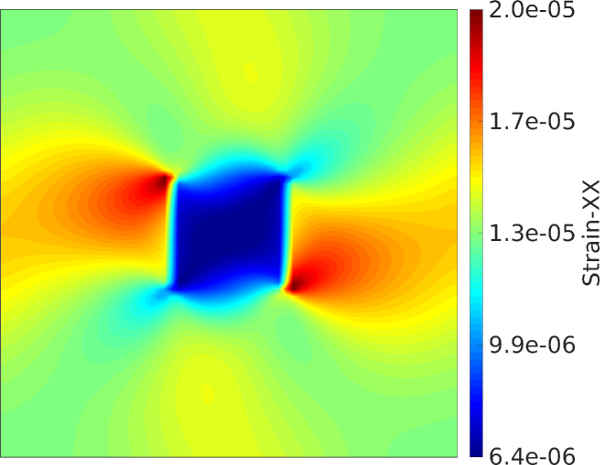}} \hspace*{1mm}
	 {\includegraphics[height=4.0cm]{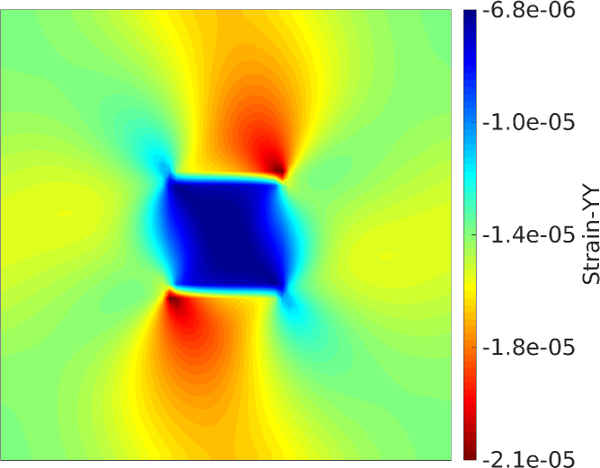}} \hspace*{1mm}
	 {\includegraphics[height=4.0cm]{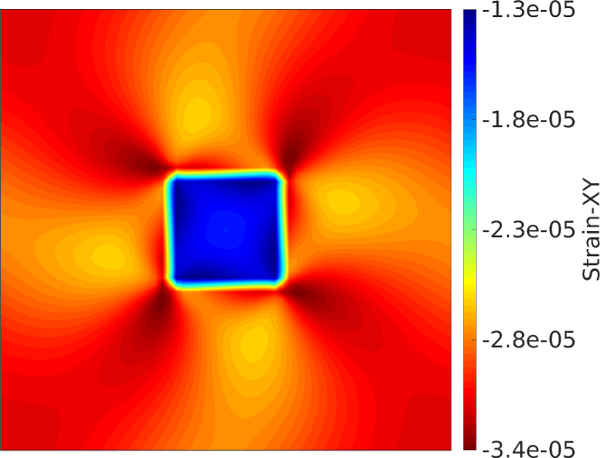}}	 
\end{minipage}  	 
\caption{Matrix-inclusion problem for beam bending. On the microdomain for the macro Gauss point at [$4978.9 \, \text{mm}$, $978.9 \, \text{mm}$] (first row) stress $\sigma_{xx}$, $\sigma_{yy}$ and $ \tau_{xy}$ in $[\text{N/mm}^2]$ and (bottom row) strain $\varepsilon_{xx}$, $\varepsilon_{yy}$ and $\gamma_{xy}$ are displayed.}
	\label{fig:MatIncl_StressRU_StrainRU}
\end{Figure} 

{\bf Remark 5} \quad The $H^1$-norm depends on the chosen physical unit of the primary quantity; the smaller the (length) unit, the larger the weight of the primary quantity compared to the derivative. In the limiting case the $H^1$-norm equals the $L^2$-norm, and the macro error in the $H^1$-norm will show the same order as in the $L^2$-norm -- in contrast to the estimates.  
In order to avoid this falsifying effect, we use in the first example of Sec.~\ref{subsec:matrix-inclusion-problem} the length unit of kilometer instead of millimeter for evaluating the $H^1$-norm.

\subsection{Assessment of superconvergence in the FE-HMM}

Since the present example is uniformly periodic, the homogenized elasticity tensor is constant everywhere including the superconvergent element center. 

\subsubsection{Superconvergence for the macro FEM} 

The error in the energy-norm is computed at mesh nodes of the coarsest discretization, here using a $25 \times 5$ mesh. Then, a hierarchical, uniform mesh refinement is carried out making 4 elements out of one of the coarsest mesh. The reference solution is obtained for $800 \times 160$ elements. The microdomain constantly exhibits a discretization of $32 \times 32$ elements, with $\epsilon = 5\, \text{mm}$.

\begin{Figure}[htbp]
	\centering
	\includegraphics[width=0.7\linewidth]{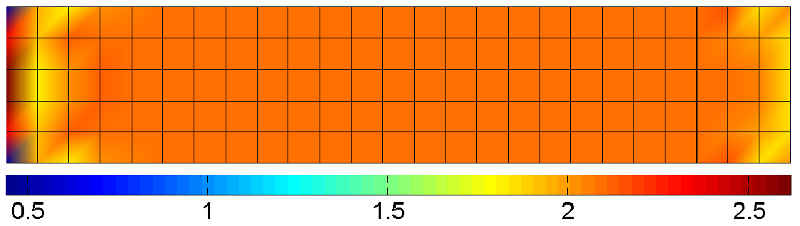}
	\caption{{Superconvergence test.} Convergence order of the local error in the energy-norm for the macrodomain in FE-HMM.}
	\label{fig:SPR_FEHMM}
\end{Figure}

Figure \ref{fig:SPR_FEHMM} displays the local convergence order of the error in the energy-norm as a contour plot in the macrodomain. Obviously, superconvergence shows up in the macro error of FE-HMM in terms of order two (instead of one according to the error estimates of standard FE-HMM) almost everywhere. A reduced order is observed at the sites of load application and reaction forces; the variations beyond and below order 2 are most pronounced at the clamped end.  
\begin{Figure}[htbp]
	\centering
	\includegraphics[width=0.4\linewidth]{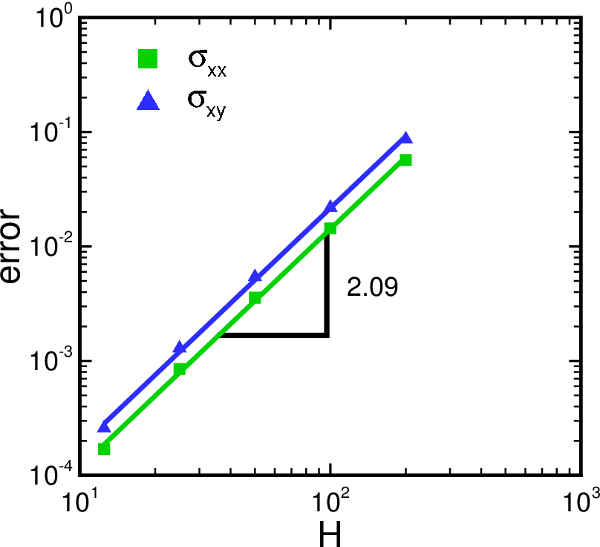}
	\caption{{Superconvergence test.} Convergence of the error in stress components $\sigma_{xx}$ and $\sigma_{xy}$ for the nodes at [$2400 \, \text{mm}, \ 400 \, \text{mm}$].}
	\label{fig:Konvergenzordnung_Makroknoten_SPR}
\end{Figure}
 
Figure~\ref{fig:Konvergenzordnung_Makroknoten_SPR} displays the convergence of the error of the stress components $\sigma_{xx}$ and $\sigma_{xy}$ as a function of the macroelement size. Component $\sigma_{yy}$ was discarded for its virtually vanishing magnitude. The convergence order of 2.09 demonstrates superconvergence of stresses at that particular node.  

\begin{table}[htbp]
	\centering
	\renewcommand{\arraystretch}{1.2}
	\begin{tabular}{l c c c c c c}
		\hline
		macroelements & $25 \times 5$ & $50 \times 10$ & $100 \times 20$ & $200 \times 40$ &  $400 \times 80$ &  \\
		microelements & \multicolumn{5}{c}{$32 \times 32$ } & \phantom{\,} {\bf order}\\
		\hline
		$|| \mathbf{u}^0 - \mathbf{u}^H ||_{A}$    & 4.6020 & 2.3821 & 1.3401 & 0.7034 & 0.2751 & $ \, $0.9888 \\ 
		\hline 
	\end{tabular}
	\caption{{Superconvergence test.} Convergence of the global macro error in the energy-norm, $\epsilon$=5\,mm.}
	\label{tab:Konvergenzordnung_global_FEHMM}
\end{table}

The integral error in the energy-norm is shown in Tab.~\ref{tab:Konvergenzordnung_global_FEHMM}. Superconvergence is not obtained. Instead we see the standard convergence order of one.
The reason for the lost superconvergence is that for the integral error all nodes of the macrodomain are considered. The nodes at the re-entrant corners of the clamped edge exhibit only orders in the range of 0.5, see Fig.~\ref{fig:SPR_FEHMM}, which spoils the overall convergence order.
  
\subsection{Problem with an analytical solution of the homogenized tensor} 
\label{subsec:Optimal-Refinement-Strategy-Numerical-Test}

\begin{Figure}[htbp]
   \begin{minipage}{15cm}  
     \centering
      \includegraphics[width=4cm, angle=0]{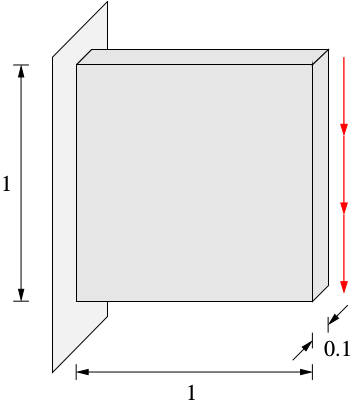} 
   \end{minipage}
\caption{Square plate for the homogenization of a microstructure endowed with an analytical solution.
\label{fig:Optimal-Refinement-Strategy-Numerical-Test}}
\end{Figure}

The BVP of the present example is a square plate of sidelength $l=1\,m$ and thickness $t=0.1\,m$, which is clamped at one edge and loaded by a constant line load of $q_0=1$ $[\text{N/mm}^2]$ at the opposite edge, see Fig.~\ref{fig:Optimal-Refinement-Strategy-Numerical-Test}.

The microheterogeneous elasticity tensor is uniformly periodic 
\begin{displaymath}
\mathbb{A}^{\epsilon}(\bm x) 
=
\left[ \begin{array}{ccc}
     {A}^{\epsilon}_{11}   &   35                        &  0   \\
        35                        & {A}^{\epsilon}_{22}  &  0   \\
        0                         &   0                         &  50  \\
\end{array}\right] \, \dfrac{\text{N}}{\text{m}^2} \qquad \mbox{with} \quad
\begin{tabular}{l}
   ${A}^{\epsilon}_{11} = [500/(5 + 3.5\cdot\mbox{sin}(2\pi x_1/\epsilon))]$  \\
   ${A}^{\epsilon}_{22} = [500/(5 + 3.5\cdot\mbox{cos}(2\pi x_1/\epsilon))]$  
\end{tabular}
\label{mat:A-epsilon-analytisch}
\end{displaymath}

with components in $[N/m^2]$. It was used in \cite{AssyrNonnenmacher2009} for heat conduction in a slightly different format. 

In contrast to the matrix-inclusion problem, the microheterogeneous tensor exhibits an analytical solution, such that 
FE-HMM can here be assessed by an exact solution without drawing to an overkill solution. 

According to \cite{Cioranescu-Donato-BOOK-1999} and \cite{AssyrNonnenmacher2011} the exact homogenized tensor $\mathbb{A}^{0}$ and the exact volumetric mean  
$\langle\mathbb{A}^{\epsilon}\rangle$ is given by

\begin{equation}
\mathbb{A}^{0} 
\approx 
\left[ \begin{array}{ccc}
        100   &   35    &  0   \\
        35    &  140    &  0   \\
        0     &    0    &  50  \\
\end{array}\right] \, \dfrac{\text{N}}{\text{m}^2} \, , 
\qquad
\langle\mathbb{A}^{\epsilon}\rangle
=
\left[ \begin{array}{ccc}
        140   &   35    &  0   \\
        35    &  140    &  0   \\
        0     &    0    &  50  \\
\end{array}\right] \, \dfrac{\text{N}}{\text{m}^2}
\label{mat:A0-analytisch-Aquer}
\end{equation}  

The comparison of ${A}^{0}$ with $\langle\mathbb{A}^{\epsilon}\rangle$ shows that the volumetric mean cannot account for the induced anisotropy due to the laminate character of the microstructure; component $\langle\mathbb{A}^{\epsilon}_{11}\rangle$ considerably deviates from the exact solution.

\subsubsection{Macro- and microconvergence} 

The simulation results 
for the macroconvergence are obtained based on a reference solution using 2560 $\times$ 2560 macroelements and $32 \times 32$ microelements. The results are invariant with respect to $\epsilon$. 
 
Figure \ref{fig:AnalyticalTensor-Macro-Micro-Convergence} (left) displays the total error as a function of the macrodiscretization for a constant microdiscretization. The observed convergence order in the $L^2$-norm of 1.62 instead of 2, and in the $H^1$-norm as well as the energy-norm of almost 0.84 instead of 1 are in reasonable agreement with the a priori estimates. The deviation from the theoretical order can be traced back to notch effects at the clamped end.

The convergence of the micro error 
in the diagram of Fig.~\ref{fig:AnalyticalTensor-Macro-Micro-Convergence} (centre) for various $h/\epsilon$ keeping $H$ fixed confirms the theoretical predictions in quantitative agreement. The reference microsolution is obtained using $50 \times 50$ macroelements and $2048 \times 2048$ microelements.

\begin{Figure}[htbp]
   \begin{minipage}{16.5cm}  
      \centering
     \includegraphics[width=5.0cm, angle=0]{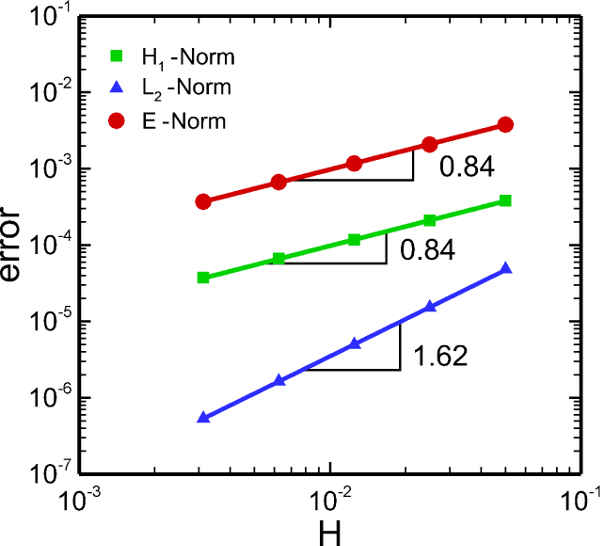} \hspace*{2mm}
     \includegraphics[width=5.0cm, angle=0]{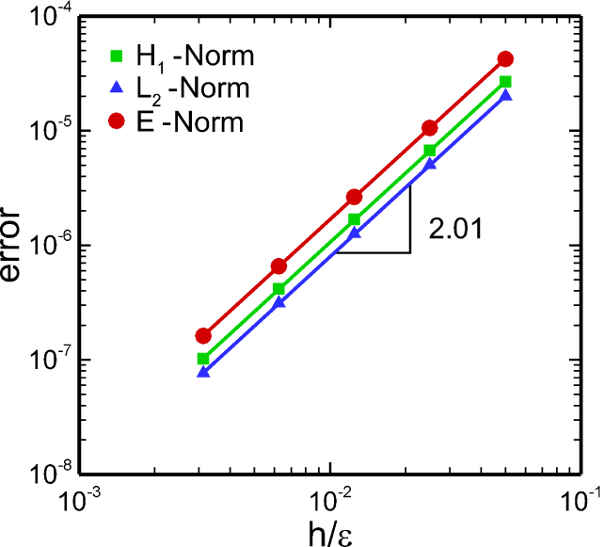} \hspace*{2mm} 
     \includegraphics[width=5.0cm, angle=0]{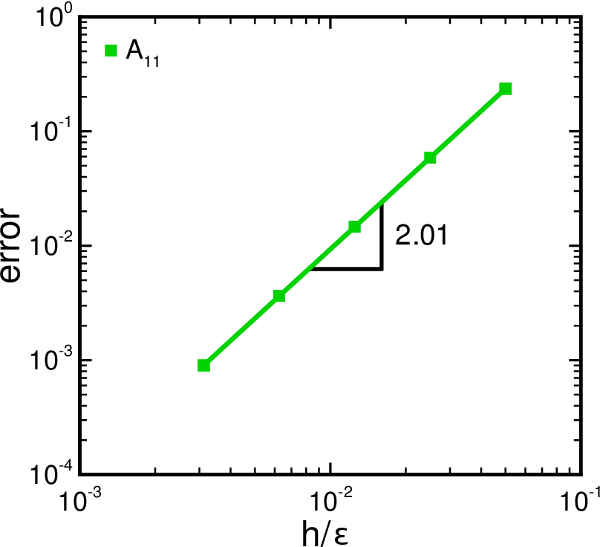} 
   \end{minipage}
\caption{Microstructure with an analytical solution of homogenization. (left) macroconvergence, microconvergence of (centre) the solution, and (right) of component $\mathbb{A}_{11}^{0,h}$ of the homogenized tensor.
\label{fig:AnalyticalTensor-Macro-Micro-Convergence}}
\end{Figure}

\subsubsection{Homogenized elasticity tensor} 

The obtained results by FE-HMM can capture the exact values of the above elasticity tensors in good agreement. A microdiscretization using $160 \times 160$ elements yields   
\begin{equation*}
\mathbb{A}^{0,h} = \begin{bmatrix} 100.0037 & 35 & 0\\ 35 & 140.0280 & 0 \\ 0 & 0 & 50 \end{bmatrix} \dfrac{\text{N}}{\text{m}^2}, \qquad \langle{\mathbb{A}}^{\epsilon}\rangle^{h} = \begin{bmatrix} 140.0280 & 35 & 0\\ 35 & 140.0280 & 0 \\ 0 & 0 & 50 \end{bmatrix} \dfrac{\text{N}}{\text{m}^2} \, .
\end{equation*}

The convergence of the components of the approximate homogenized elasticity tensor $\mathbb{A}_{}^{0,h}$ is 
displayed in the right diagram of Fig.~\ref{fig:AnalyticalTensor-Macro-Micro-Convergence}. The slope indicates order 2 in agreement with the theory, $|\mathbb{A}_{}^{0,h} - \mathbb{A}^{0}| \leq C (h/\epsilon)^2$. The reference solution $\mathbb{A}_{}^{0}$ is obtained on a microdomain with a discretization of 2560 $\times$ 2560 elements. Note that component $\mathbb{A}_{22}^{0,h}$ does not change in the leading 9 digits for meshes finer than 40 $\times$ 40 elements, the reason, why the diagram in Fig.~\ref{fig:AnalyticalTensor-Macro-Micro-Convergence} is restricted to component $\mathbb{A}_{11}^{0,h}$.
 
\subsubsection{Comparison of FE-HMM with FE$^2$} 
Table \ref{tab:Analytical-Tensor-FEHMM-vs-FEM_A0h-vs-FEM_Amean} shows, that FE-HMM and FE$^2$ exhibit virtually the same results. 

\begin{table}[htbp]
\center
\renewcommand{\arraystretch}{1.2}
\begin{tabular}{lrrrrr}
\hline \\[-5mm]         
& microelements                             & \multicolumn{4}{c}{$20 \times 20$ }  \\
& macroelements                             &  $20 \times 20$  & $40 \times 40$  & $80 \times 80$  &  $160 \times 160$   \\[1mm]
\hline \\[-5mm]                       
FE-HMM   &  $|| \bm u ||_{\infty}$ in $10^{-3}$       &  77.7902   &  78.3328   &  78.5746  &  78.6865      \\
         &  $|| \bm u ||_A$ in $10^{-3}$               &  81.6746   &  81.7999   &  81.8388  &  81.8510      \\[1mm]
\hline \\[-5mm]
FE$^2$ $\mathrel{\widehat{=}}$  &  $|| \bm u ||_{\infty}$ in $10^{-3}$       &  77.7902 &  78.3328  &  78.5746  &  78.6865  \\
FE-HMM($\mathbb{A}^{0,h}$)      &  deviation                     &   2.82 $\cdot 10^{-13}$   &  1.50  $\cdot 10^{-13}$ &  2.20 $\cdot 10^{-13}$  &  −1.64 $\cdot 10^{-11}$       \\ 
                            &  $|| \bm u ||_A$ in $10^{-3}$  &  81.6746  &  81.7999   &  81.8388   &  81.8510      \\
                            &  deviation                     &  1.38 $\cdot 10^{-13}$  &   1.16 $\cdot 10^{-13}$ & 1.26 $\cdot 10^{-12}$  &  −8.78 $\cdot 10^{-12}$       \\[1mm]
\hline
\end{tabular}
\caption{{Microstructure with an analytical solution of homogenization.} Convergence of the standard FE-HMM solution in comparison to the FE-HMM($\mathbb{A}^{0,h}$)/FE$^2$ solution and the percentaged deviation from standard FE-HMM.  
\label{tab:Analytical-Tensor-FEHMM-vs-FEM_A0h-vs-FEM_Amean}}
\end{table}
 
\subsubsection{Optimal refinement strategies} 

The a priori error estimates ($\delta/\epsilon \in \mathbb{N}$ along with periodic coupling) enable an optimal uniform refinement strategy for uniform 
macro- ($H$) and micro- ($h$) meshes, see \cite{Assyr2009}.
For $p=q=1$ i.e. bilinear shape functions for the macro- as well as on the micro-FE, we obtain 
\begin{displaymath}
   || \bm u^0 - \bm u^H ||_{L^2(\Omega)} \leq C\left( H^{2} + \left(\dfrac{h}{\epsilon}\right)^{2} \right)\, , \quad 
   || \bm u^0 - \bm u^H ||_{H^1(\Omega)} \leq C\left( H + \left(\dfrac{h}{\epsilon}\right)^{2} \right) \, . 
\end{displaymath}
We denote by $L$ the number of micro elements in each direction of space of the microdomain. Similarly, $M$ is the number of macro elements per direction.
Consequently, the micro element size is $h=\epsilon/L$ and it holds $h_L := h/\epsilon=1/L$ for the micromesh, and $H_M=1/M$ for the macromesh. 
 
Denoting by $N_{mac}$ the number of macro DOF and by $N_{mic}$ the number of micro DOF, the above rates of convergence show that
the following choices are the best uniform refinement strategies that is full order is achieved for minimal computational costs. 
\begin{center}
\begin{tabular}[4]{rlcl}
  $L^2$-norm: &  $N_{mic}=N_{mac}$  & $\Longleftrightarrow$  & $h_L = H_M$  \\
  $H^1$-norm: &  $N_{mic}=\sqrt{N_{mac}}$  & $\Longleftrightarrow$  & $h_L = \sqrt{H_M}$  \\
\end{tabular}
\end{center}
 
\begin{Figure}[htbp]
   \begin{minipage}{16.2cm}  
     \centering
     \includegraphics[height=5.0cm, angle=0]{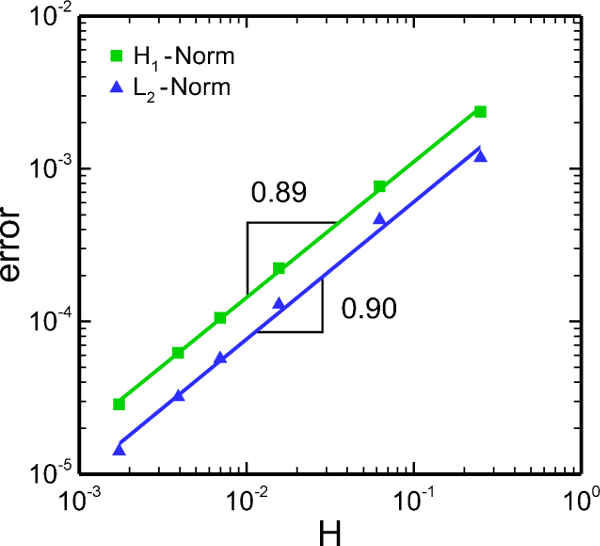} \hspace*{1mm}
     \includegraphics[height=5.0cm, angle=0]{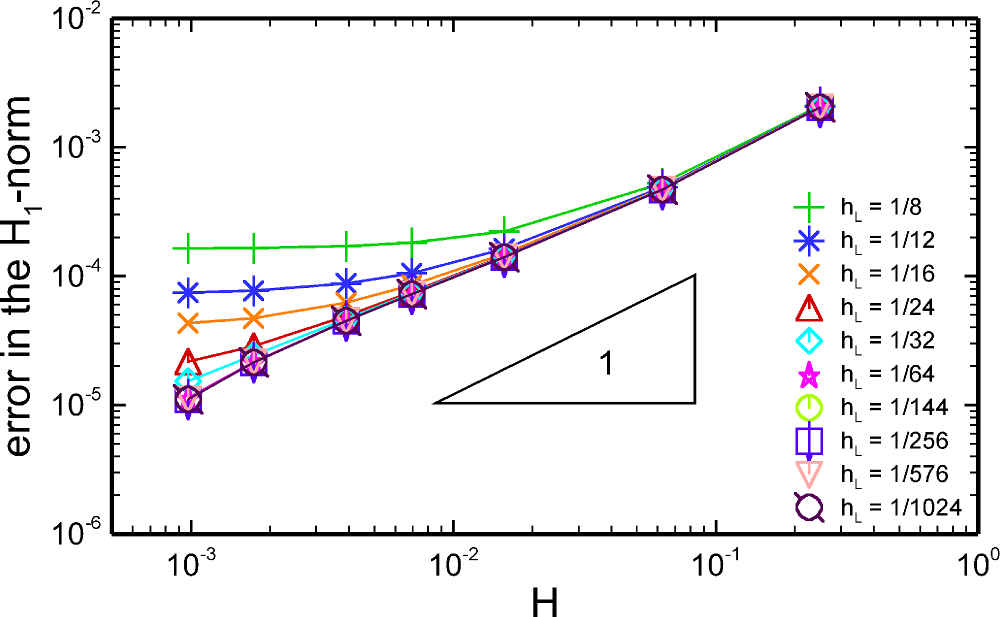} 
   \end{minipage}
\caption{Optimal and suboptimal refinement strategies for the $H^1$-error.
\label{fig:Optimal-vs-Supoptimal-Refinement-Strategy-H1}}
\end{Figure}

\begin{Figure}[htbp]
   \begin{minipage}{16.2cm}  
     \centering
     \includegraphics[height=5.0cm, angle=0]{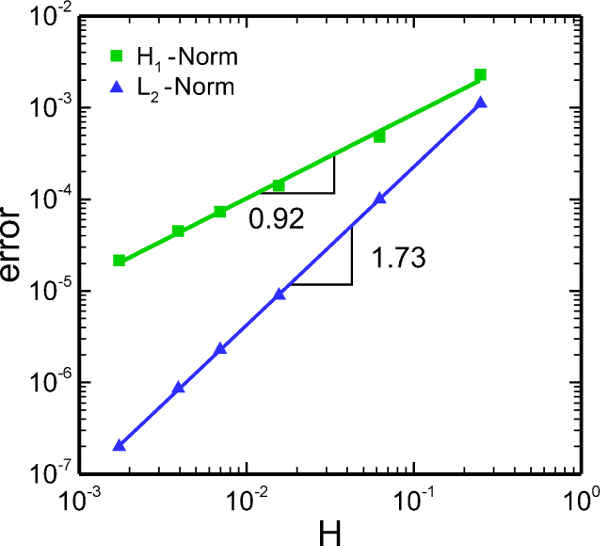} \hspace*{1mm}  
     \includegraphics[height=5.0cm, angle=0]{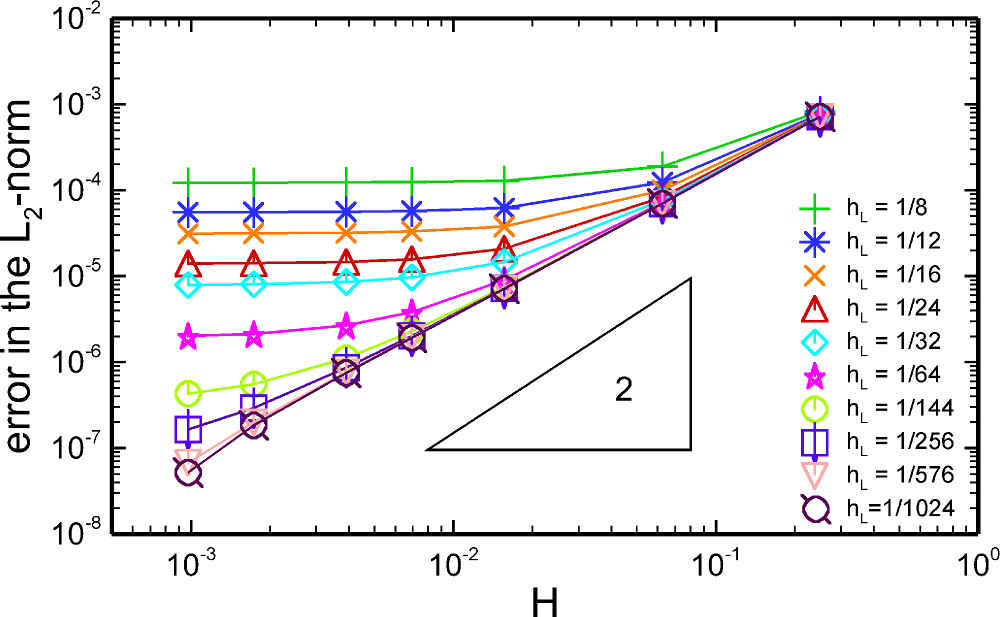} 
   \end{minipage}
\caption{Optimal and suboptimal refinement strategies for the $L^2$-error.
\label{fig:Optimal-vs-Supoptimal-Refinement-Strategy-L2}}
\end{Figure}

We compute the solutions of the problem for uniform macro mesh refinement $H_M=1/8, 1/12, 1/16, \ldots, 1/1024$. 
In the diagrams of Fig.~\ref{fig:Optimal-vs-Supoptimal-Refinement-Strategy-H1} (right) and Fig.~\ref{fig:Optimal-vs-Supoptimal-Refinement-Strategy-L2} (right)
the micromesh is kept fixed for each of the distinct solid lines and is successively refined from one line to the other. 
Optimal refinements clearly follow the ratio $h_L=H_M$  for the error in the $L^2$-norm, Fig.~\ref{fig:Optimal-vs-Supoptimal-Refinement-Strategy-L2}, 
and $h_L=\sqrt{H_L}$ for the error in the $H^1$-norm, Fig.~\ref{fig:Optimal-vs-Supoptimal-Refinement-Strategy-H1}.
The results demonstrate the sharpness of the a priori bounds. 
Similar results for piecewise bilinear FE are reported in \cite{Assyr2009}. 

\subsection{Non-uniformly periodic tensor} 
\label{subsec:Non-uniformly-periodic-tensor}

\begin{Figure}[htbp]
   \begin{minipage}{16cm}  
     \centering
      \includegraphics[width=6cm, angle=0]{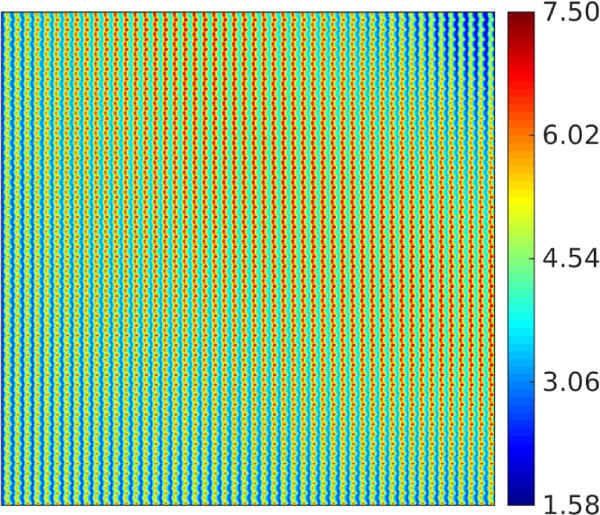} 
   \end{minipage}
\caption{Young's modulus distribution for a non-uniformly periodic elasticity tensor for $\epsilon = 0.005\,$m on a domain $[0, 1]$\,m $\times$ $[0, 1]$\,m.
\label{fig:Non-uniformly-periodic-tensor}}
\end{Figure}

The present example serves the purpose to analyze FE-HMM for a non-uniformly periodic microstructure, which is given by a Young's modulus of 
\begin{displaymath}
    E^{\epsilon} =         \dfrac{1.5 + \mbox{sin}(2\pi x_1/\epsilon)}{1.5 + \mbox{sin}(2\pi x_2/\epsilon)} 
                          + \dfrac{1.5 + \mbox{sin}(2\pi x_2/\epsilon)}{1.5 + \mbox{sin}(2\pi x_1/\epsilon)} 
                          + \mbox{sin}(4x_1 x_2) + 1 \, \quad [\mbox{N/m}^2]\, . 
\end{displaymath}
as displayed for $\epsilon=0.005$\,m in Fig.~\ref{fig:Non-uniformly-periodic-tensor} on the domain $[0, 1]$\,m $\times$ $[0, 1]$\,m, 
see \cite{AssyrNonnenmacher2009}. The coordinate system is centered in the domain, axes align with the boundaries.
 
The macro BVP is almost the same as in the previous example of Sec.~\ref{subsec:Optimal-Refinement-Strategy-Numerical-Test}. 
Only the line load at the free edge is adopted to $q_0=0.01$\,N/m, since the  Young's modulus is significantly smaller than in the previous example.

\subsubsection{Macro- and microconvergence} 

\begin{Figure}[htbp]
   \begin{minipage}{16.5cm}  
      \centering
     \includegraphics[width=5.0cm, angle=0]{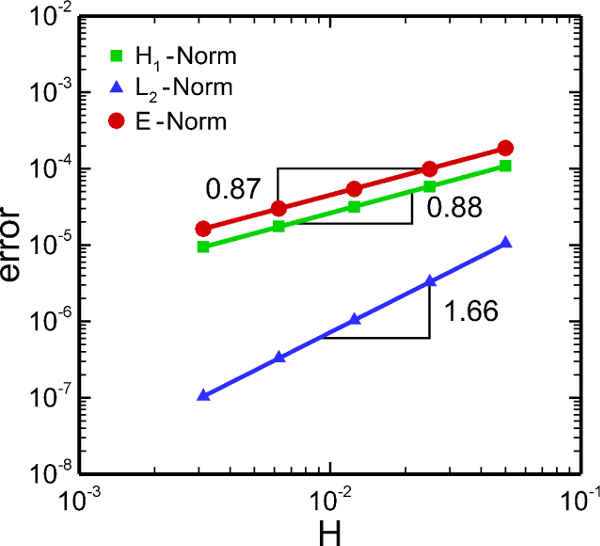} \hspace*{2mm}
     \includegraphics[width=5.0cm, angle=0]{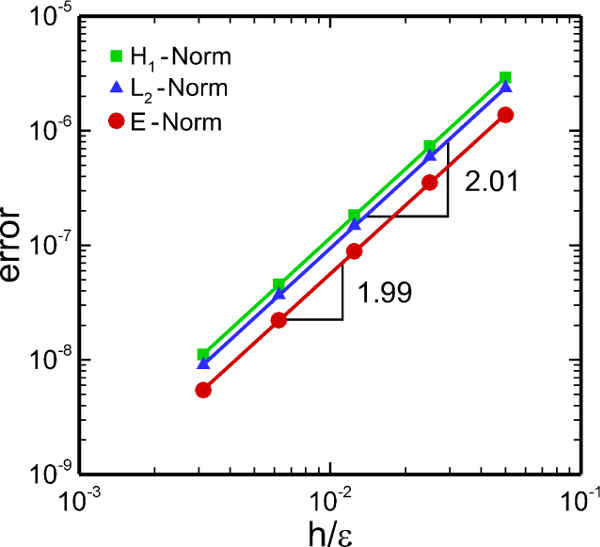} \hspace*{2mm} 
     \includegraphics[width=5.0cm, angle=0]{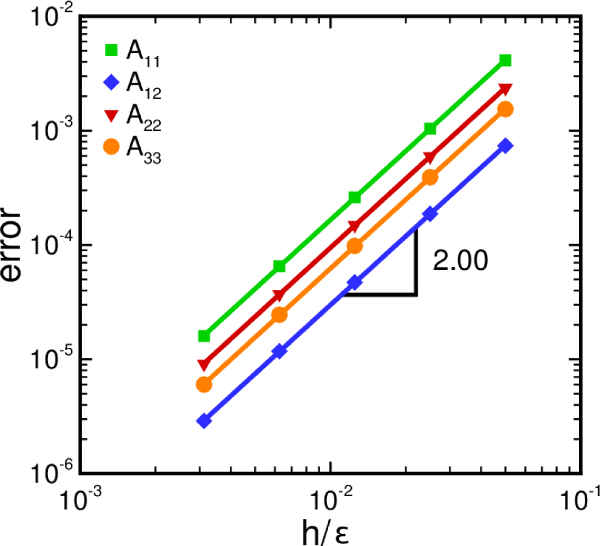}  
   \end{minipage}
\caption{Non-uniformly periodic tensor; (left) macroconvergence, microconvergence of (centre) the solution, and (right) the components of the homogenized tensor $\mathbb{A}^{0,h}$.
\label{fig:Non-uniformly periodic tensor-Macro-Micro-Convergence}}
\end{Figure}

The convergence of the macro error for a fixed microdiscretization is displayed 
in the left diagram of Fig.~\ref{fig:Non-uniformly periodic tensor-Macro-Micro-Convergence}. The corresponding reference solution is obtained by $2560 \times 2560$ macroelements and $40 \times 40$ microelements. 
The convergence order is below the theoretical order in the $L^2$-norm (1.66 instead of 2), and in the $H^1$- as well as in the energy-norm (0.88 instead of 1), which is due to the reduced regularity of the macro BVP, compare the results with those in Sec.~\ref{subsec:Optimal-Refinement-Strategy-Numerical-Test}.

The convergence of the micro error for a fixed macrodiscretization is displayed 
in the centre diagram of Fig.~\ref{fig:Non-uniformly periodic tensor-Macro-Micro-Convergence}. Full order of convergence is obtained for the micro error in the $L^2$-norm and in the $H^1$-norm, respectively. The reference solution for microerror calculation is obtained by $40 \times 40$ macroelements and $2048 \times 2048$ microelements.

\subsubsection{Homogenized elasticity tensor} 

Next, we assess the convergence estimate for the components of the approximate homogenized elasticity tensor $\mathbb{A}_{}^{0,h}$.  
The convergence as a function of the microdiscretization is shown 
in the right diagram of Fig.~\ref{fig:Non-uniformly periodic tensor-Macro-Micro-Convergence}. 
The slope indicates order 2 for all components of the elasticity tensor in the plane strain case in perfect agreement with the theory, \cite{Assyr2009}.

\subsubsection{Comparison of the FE-HMM variants, and with the FE$^2$ method} 

The deviation of the FE-HMM solution based on micro-macro stiffness transfer in comparison with the second version using the homogenized elasticity
tensor is throughout less than $4 \cdot 10^{-11}$ and thus in the range of machine precision, see Tab.~\ref{tab:Non-uniformly periodic tensor-FEHMM-vs-FE2-vs}. Since the latter FE-HMM version coincides
with FE$^2$, the coincidence of FE-HMM and FE$^2$ is underpinned again by numerical means in the present non-uniformly periodic setting.  

\begin{table}[htbp]
\center
\renewcommand{\arraystretch}{1.2}
\begin{tabular}{lcrrrr}
\hline \\[-5mm]         
& microelements                                  &  \multicolumn{4}{c}{$20 \times 20$ }  \\
& macroelements                                  &  $20 \times 20$  & $40 \times 40$  & $80 \times 80$  &  $160 \times 160$   \\[1mm]
\hline \\[-5mm]                       
FE-HMM   &  $|| \bm u ||_{\infty}$ in $10^{-3}$ &   20.0246        & 20.1980         &  20.2810        &  20.3214     \\
         &  $|| \bm u ||_A$ in $10^{-3}$        &   4.15603        &  4.16176        &   4.16353       &   4.16408    \\[1mm]
\hline \\[-5mm] 
FE$^2$ $\mathrel{\widehat{=}}$  &  $|| \bm u ||_{\infty}$ in $10^{-3}$  &   20.0246      &  20.1980        &  20.2810    &  20.3214   \\
FE-HMM($\mathbb{A}^{0,h}$)      &  deviation                            &     3.79  $\cdot 10^{-11}$ &   4.47  $\cdot 10^{-14}$ &  1.41 $\cdot 10^{-14}$  &  -1.29 $\cdot 10^{-13}$    \\ 
                                &  $|| \bm u ||_A$ in $10^{-3}$         &    4.15603                 &   4.16176  &  4.16353   &   4.16408  \\
                                &  deviation                            &     2.73 $\cdot 10^{-14}$  &   2.81 $\cdot 10^{-14}$ &  -8.33 $\cdot 10^{-16}$  &  -2.70 $\cdot 10^{-13}$      \\[1mm]
\hline
\end{tabular}
\caption{{Non-uniformly periodic tensor.} Convergence of the standard FE-HMM solution in comparison to the FE$^2$ solution along with the perecentaged deviation.  
\label{tab:Non-uniformly periodic tensor-FEHMM-vs-FE2-vs}}
\end{table}


%
\subsection{Escher's Bird and Fish tessellation} 
\label{subsec:Eschers-Bird-and-Fish} 
%
\begin{Figure}[htbp]
   \begin{minipage}{16.0cm}  
   \centering  
         \includegraphics[height=5.2cm, angle=0]{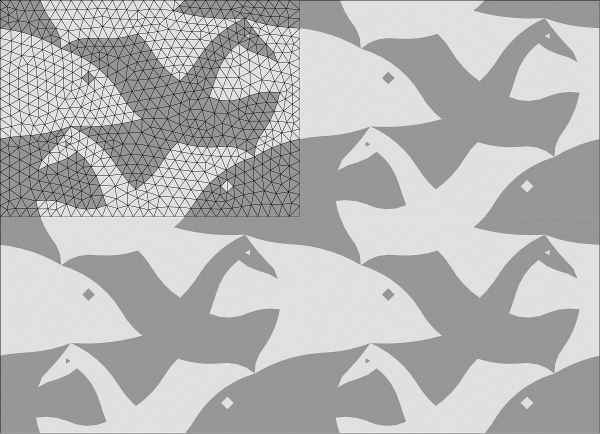}        
   \end{minipage}
\caption{Escher's periodic Bird$/$Fish tessellation: RVE and its discretization using triangular elements with $q=1$.
\label{fig:Escher-Bird-and-Fish-tessellation-discretization}}
\end{Figure}

While Lord Kelvin was concerned with the problem, how to pack equal-sized objects together to fill space with minimal surface area, the graphic artist 
M.C. Escher (1898--1972) was devoted to the question, how to fill the plane by tessellations\footnote{Tessellations are arrangements of closed shapes that completely cover the plane without overlapping and without leaving gaps.} with most intriguing beauty. 

While from 1887 until 1993 the tetrakaidecahedron was continuously but erroneously believed\footnote{The so-called Weaire-Phelan structure was shown in 1993 to solve the Kelvin problem with even smaller surface area.} to be the right answer of Kelvin to his minimization problem, the answer to Escher's question is simple in that it is clearly non-unique at any time, as beauty is in the eye of the beholder\footnote{\emph{''At moments of great enthusiasm it seems to me that no one in the world has ever made something this beautiful and important.''} -- M.C. Escher}.

It is a well-known mathematical result that of all the regular polygons, only the triangle, square, and hexagon can be used for a tessellation.
Escher combined them applying reflections, glide reflections, translations, and rotations to augment his playground for creating a larger 
variety of tessellations.

From Escher's oeuvre in the field which he called ''regular divisions of the plane'' we devote our analysis to the ''Bird/Fish (1941 B)'' tessellation, an artistic artifact that is uniformly periodic and consists of two phases, Fig.~\ref{fig:Escher-Bird-and-Fish-tessellation-discretization}. 
Here we choose $\delta = \varepsilon=28.5$\,mm. The macroscopic BVP is the same as in Sec.~\ref{subsec:matrix-inclusion-problem}, and the Poisson's ratio and Young's moduli exhibit the same values, 
thus showing a contrast of the bright face to the dark face of $E_b/E_d=2.5$.

\subsubsection{Macro and micro convergence} 
\begin{Figure}[htbp]
   \begin{minipage}{16.5cm}  
      \centering                              
     \includegraphics[width=6.0cm, angle=0]{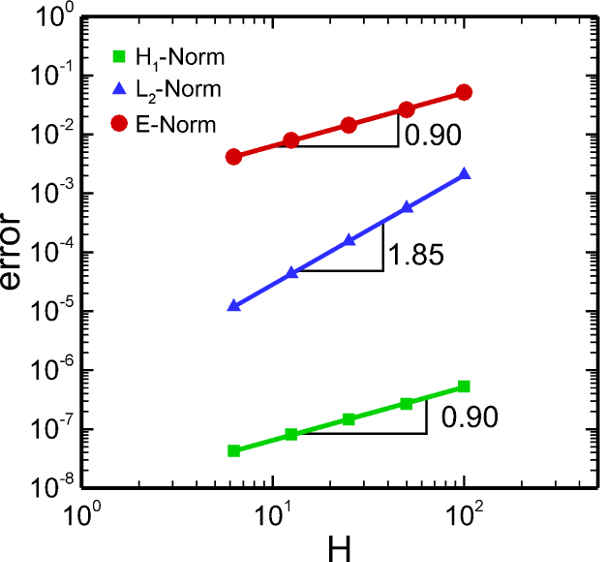} \hspace*{2mm}
     \includegraphics[width=6.0cm, angle=0]{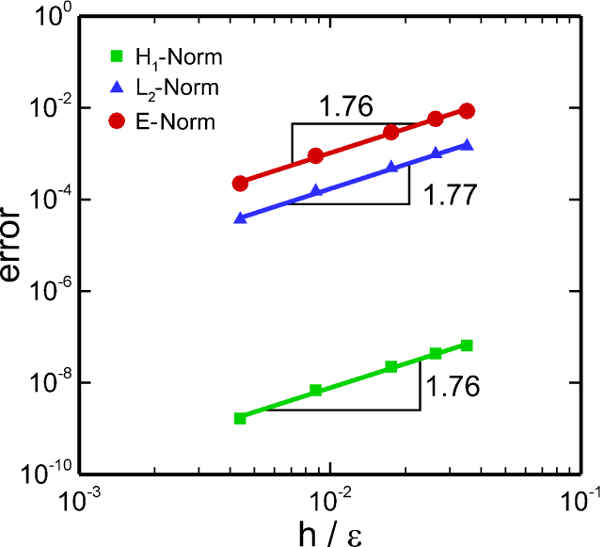}  
   \end{minipage}
\caption{Escher's periodic Bird/Fish tessellation. (Left:) Macroconvergence, (right:) microconvergence.
\label{fig:EscherBirdFish2d-Macro-Micro-Convergence}}
\end{Figure}
Figure~\ref{fig:EscherBirdFish2d-Macro-Micro-Convergence} shows in its left diagram the convergence of the normed error for macroscopic mesh refinement for a reference solution of 4\,096\,000 DOFs (3200$\times$640 macroelements) and with 1462 DOFs on each microdomain. 
Similarly, Fig.~\ref{fig:EscherBirdFish2d-Macro-Micro-Convergence} shows in its right diagram the convergence of the normed errors for microscopic mesh refinement for a reference solution of 16\,000 DOFs on the macrodomain (200$\times$40 macroelements) and 312\,506 DOFs on each microdomain.

The convergence rates for $p=q=1$ reasonably agree with the theoretical values. The deviations are in the reduced regularity of the problem sets on the two scales.
 
\begin{Figure}[htbp]
   \begin{minipage}{16.0cm}  
   \centering  
        \includegraphics[height=5cm, angle=0]{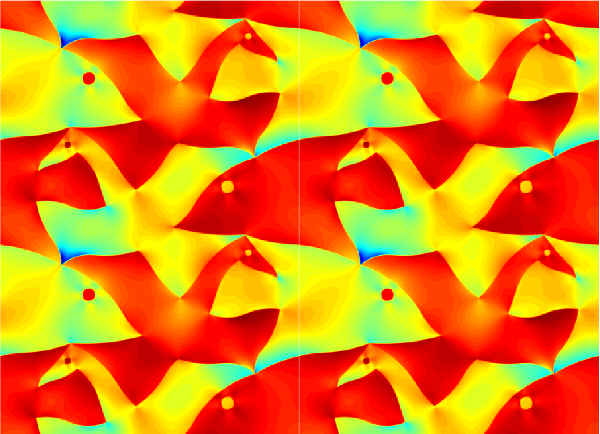} \hspace*{-4.2mm}
        \includegraphics[height=5cm, angle=0]{birdfish_FEM_eps_xy_UR.jpg}  
        \\[4mm]
        \includegraphics[height=5cm, angle=0]{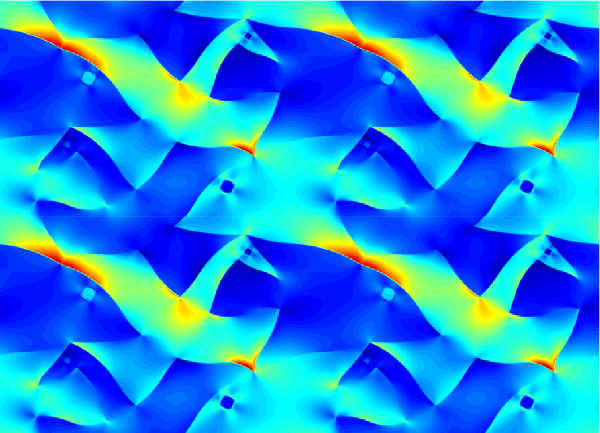} \hspace*{-4.3mm}
        \includegraphics[height=5cm, angle=0]{RO_Sigma_VM_cut_cut_2x2.jpg}  
   \end{minipage}
\caption{Escher's periodic Bird$/$Fish tessellation. Contour plot of (top:) shear strain $\gamma_{xy}$ and of (bottom) von-Mises stress at the macro quadrature point  ($x=4978.9$\,mm, $y=978.9$\,mm).\label{fig:Escher-Bird-and-Fish-ShearStrain}}
\end{Figure}

At the end, we abstain from homogenizing (the elasticities of) bird and fish and discover the microscale simulation results displayed in Fig.~\ref{fig:Escher-Bird-and-Fish-ShearStrain}.
 
%
\section{Summary}

The novel contributions and main results shall be summarized.
\begin{enumerate}
 \item Based on the seminal contributions of E, Engquist and Abdulle
       we have presented a finite element an FE-HMM formulation for linear elasticity in a geometrical linear frame. For that case of a 
       vector-valued field problem we have generalized the FE-HMM micro-to-macro stiffness transfer by deriving corresponding transformation matrices which follow
       from macro element unit displacement states and result in energy minimizers.
       By virtue of the present setting in solid mechanics, a clear mechanical interpretation of FE-HMM as an energy-based homogenization 
       scheme is facilitated. Moreover, this setting opens the door to a conceptual and numerical comparison with other two-scale finite element 
       homogenization concepts. 
 \item A conceptual comparison of the FE-HMM with the FE$^2$ method was presented; they both have in common the equivalence of micro- with 
       macroenergy densities; it was shown that for periodic boundary conditions the macro-homogeneity condition of Hill-Mandel as the 
       theoretical cornerstone of FE$^2$ can be derived from asymptotic homogenization as the theoretical foundation of FE-HMM, which 
       is the deepest common root of the two methods.    
       FE$^2$ in its micro-to-macro stiffness transfer draws on the homogenized fourth order constitutive tensor while keeping the macroscopic 
       bilinear form unaltered. FE-HMM however crucially relies on a modified quadrature formula in the macroscopic bilinear form, which 
       results in a micro-macro stiffness matrix transfer, a setting which enables the a priori error analysis. For that reason the stiffness 
       matrix transfer is a key ingredient of FE-HMM. The homogenized elasticity tensor --coinciding with the one of FE$^2$-- must be 
       seen as a byproduct. 

       The existence of a priori estimates is the particular strength and key novel contribution of FE-HMM to the field of numerical homogenization in general, and 
       in particular to two-scale finite element methods like FE$^2$. In the mathematical cornerstones of FE-HMM Therein are the asymptotic homogenization framework 
       and the setting of numerical quadrature. Here, FE-HMM and its protagonists, W. E, B. Engquist, A. Abdulle and coworkers have 
       brought mathematical rigor into the fully discrete analysis of two-scale finite element methods for homogenization. For FE$^2$ methods, in spite of 
       their considerable progress in a multitude of applications, a thorough mathematical analysis is missing. 
 \item In representative benchmark problems FE-HMM has been assessed for periodic boundary conditions imposed on the microdomains. 
       Here we checked in various norms the a priori error estimates as given by \cite{E-Ming-Zhang-2005} for the general case of elliptic PDEs 
       and for linear elasticity by \cite{Assyr2006}. For sufficient regularity of the BVP, the estimates have been verified for both the macro error and 
       the micro error and for the coefficients of the homogenized elasticity tensor. 
       \\
       Remarkably, superconvergence holds for the micro error in FE-HMM in that stress (and equally the error in the energy-norm and in the $H^1$-norm) converges in the same order as displacements (and equally the error in the $L^2$-norm). 
       In contrast to the notion of superconvergence in  standard, single-scale finite element methods, which holds only at particular points for particular 
       element shapes, the verified estimates for the micro error are generally valid. We have demonstrated by the application of the Superconvergent Patch Recovery (SPR) 
       technique, that for the macro error in FE-HMM standard superconvergence holds at the so-called Barlow points. 
       \\
       Moreover, the existing predictions for FE-HMM concerning the optimal uniform micro-macro refinement strategies on how to obtain full order for minimal computational costs
       are corroborated by the numerical results in the present work. 
 \item Furthermore, the simulations have shown quantitative agreement of the two FE-HMM variants with FE$^2$.
       Of course, for the coincidence of the two methods in practice, FE-HMM findings equally apply for the FE$^2$ method with respect to (i) the micro and macro convergence orders 
       including superconvergence and with respect to (ii) the derived optimal uniform micro-macro mesh refinement strategies.
\end{enumerate}

In conclusion, the present work shall contribute to bridge the still existing gap between the mathematics of multiscale modeling with FE-HMM as an excellent example, and, on the other hand side, the engineering community, in which the popular FE$^2$ method has been growing into different branches of solid mechanics by virtue of its versatility. Here, FE-HMM as a non-standard discretization method, which is endowed with a fully-fledged mathematical structure at least for the linear case, advances numerical homogenization more than by an $\epsilon$ towards reliability as a key ingredient of scientific computations.

 
\bigskip

{\bf Acknowledgements.} The first author acknowledges support by the Deutsche Forschungsgemeinschaft (DFG) 
within the Heisenberg program (grant no. EI 453$/$2-1).



\begin{appendix}

\addcontentsline{toc}{section}{Appendix}
\renewcommand{\thesubsection}{\Alph{section}.\arabic{subsection}}
\renewcommand{\theequation}{\Alph{section}.\arabic{equation}}
\renewcommand{\thefigure}{\Alph{section}.\arabic{figure}}
\renewcommand{\thetable}{\Alph{section}.\arabic{table}}
\newcommand {\ssectapp}{
                        \setcounter{equation}{0}
                        \setcounter{figure}{0}
                        \setcounter{table}{0}
		                \subsection
                        }

\setcounter{equation}{0}

\vfill
\newpage

 
\section{\color{black} Mathematical basis of FE-HMM}
\label{sect:Basis-of-FEHMM}

\subsection{Asymptotic homogenization}
\label{subsec:Homogenization}
 
This part presents the derivation and main results of asymptotic homogenization for linearized elasticity. It is 
a main reference point to underpin the mathematical basis of FE-HMM. For mathematical homogenization by asymptotic expansion
we refer to \cite{Bensoussan-Lions-Papanicolau-BOOK-1976}, \cite{Sanchez-Palencia-BOOK-1980}, \cite{Allaire1992}, \cite{Cioranescu-Donato-BOOK-1999}. 

In the two-scale mathematical homogenization theory, various fields are assumed to depend on two coordinates, $\bm x$ as the coarse- or macroscale position vector, and $\bm y$ as the fine- or microscale position vector. These two vectors are related by $\bm y=\bm x/\epsilon$ with $0<\epsilon \ll 1$. 

The unit cell domain $Y$ is typically chosen a rectangular domain, in which the microstructure is assumed to be locally periodic ($Y$-periodic). Even if the microstructure is not periodic, the response functions $\bm r^{\varepsilon}$ are in either case assumed to be locally periodic. 
 
In addition to the symmetries introduced in Sec.~\ref{subsec:StrongFormLinearElasticity} 
we assume for $\mathbb{A}^{\epsilon}_{ijlm}$
\begin{equation}
 \mathbb{A}^{\epsilon}_{ijlm} (\bm x) = \mathbb{A}_{ijlm} \left(\frac{\bm x}{\epsilon}\right) = \mathbb{A}_{ijlm} (\bm y) 
 \,\, \mbox{is $Y$-periodic for any $i, j, l, m=1, \ldots, d$}.
\end{equation}
The derivative of a response function $\bm r^{\epsilon}(\bm x)$ with respect to the macroscale variable $x_i$ can be calculated by the chain rule according to 
\begin{equation}
 \dfrac{d}{dx_i} = \dfrac{\partial}{\partial x_i} + \dfrac{1}{\epsilon} \dfrac{\partial}{\partial y_i} \, .
 \label{eq:total-derivative}
\end{equation}
For the case of linear elasticity we consider the yet unknown response functions $\bm u^{\epsilon}$, $\varepsilon_{ij}^{\epsilon}$, $\sigma_{ij}^{\epsilon}$ and their expansions into series
\begin{eqnarray}
     u_i^{\epsilon} (\bm x, \bm y) &=& u_i^{0} (\bm x) + \epsilon \, u_i^{1} (\bm x, \bm y) + \epsilon^2  \, u_i^{2} (\bm x, \bm y)  + \mathcal{O}(\epsilon^3) 
    \label{eq:AsympExp-u} \\
    \varepsilon_{ij}^{\epsilon} (\bm x, \bm y) &=& \dfrac{1}{\epsilon} \, \varepsilon_{ij}^{-1} (\bm x, \bm y) + \varepsilon_{ij}^{0} (\bm x, \bm y) + \epsilon \, \varepsilon_{ij}^{1} (\bm x, \bm y) 
                                  + \mathcal{O}(\epsilon^2) 
    \label{eq:AsympExp-strain} \\
    \sigma_{ij}^{\epsilon} (\bm x, \bm y) &=& \dfrac{1}{\epsilon} \, \sigma_{ij}^{-1} (\bm x, \bm y) + \sigma_{ij}^{0} (\bm x, \bm y) + \epsilon \, \sigma_{ij}^{1} (\bm x, \bm y) 
                                  + \mathcal{O}(\epsilon^2) \,.
    \label{eq:AsympExp-sigma}
\end{eqnarray}
The strain expansion \eqref{eq:AsympExp-strain} is obtained by inserting the expansion for $u_i^{\epsilon}(\bm x)$ into \eqref{eq:kinematical-relation} with an account of \eqref{eq:total-derivative} and with the following nomenclature for the strain components of various ''orders''
\begin{equation}
  \varepsilon_{ij}^{-1} = u_{(ij)y}^{0}\, , \qquad \varepsilon_{ij}^{s} = u_{(ij)x}^{s} + u_{(ij)y}^{s+1} \quad \text{for} \,\, s=0,1,\ldots\,  
  \label{eq:strain-relations-various-orders}
\end{equation}
where the abbreviation $u^{k}_{(ij)y}= 1/2 (\partial u_i^k / \partial y_j +  \partial u_j^k / \partial y_i)$ (and analogously for $u^{k}_{(ij)x}$) is employed.
{\color{black} For each scale an elasticity relation is assumed to hold}  
\begin{equation}
   \sigma_{ij}^{s} (\bm x, \bm y) = \mathbb{A}^{\epsilon}_{ijkl} (\bm x) \, \varepsilon_{kl}^{s} (\bm x, \bm y) \quad \text{for} \,\, s=-1, 0, 1, \, \ldots \, .
   \label{eq:stress-strain-relations-various-orders}
\end{equation}
  
Inserting \eqref{eq:AsympExp-sigma} into the balance of linear momentum \eqref{eq:StrongFormMicro-2}$_1$ yields 
\begin{equation}
\label{eq:balance-linear-momentum-asymp-exp}
                \dfrac{1}{\epsilon^{2}}  \dfrac{\partial \sigma_{ij}^{-1}}{\partial y_j}  +  
                \dfrac{1}{\epsilon} \left( \dfrac{\partial \sigma_{ij}^{-1}}{\partial x_j} + \dfrac{\partial \sigma_{ij}^0}{\partial y_j}\right) +
                \left( \dfrac{\partial \sigma_{ij}^0}{\partial x_j} + \dfrac{\partial \sigma_{ij}^1}{\partial y_j} + f_i^{} \right) +
                \mathcal{O}(\epsilon) = 0 \, ,
\end{equation}
from which, according to the order, different balance equations follow (by multiplying \eqref{eq:balance-linear-momentum-asymp-exp} by $\epsilon^{n}$ and
considering the limit $\epsilon  \rightarrow 0^+$, each for $n=2,1,0$) 
\begin{eqnarray}
  \mathcal{O}\left(\epsilon^{-2}\right):  & & \dfrac{\partial \sigma_{ij}^{-1}}{\partial y_j} = 0 
  \label{AsympExp-BalanceEq-1}\\
  \mathcal{O}\left(\epsilon^{-1}\right):  & & \dfrac{\partial \sigma_{ij}^{-1}}{\partial x_j} + \dfrac{\partial \sigma_{ij}^0}{\partial y_j} = 0 
  \label{AsympExp-BalanceEq-2}\\
  \mathcal{O}\left(1 \right):  & & \dfrac{\partial \sigma_{ij}^0}{\partial x_j} + \dfrac{\partial \sigma_{ij}^1}{\partial y_j} + f_i^{} = 0 \, .
  \label{AsympExp-BalanceEq-3}
\end{eqnarray}

For the Dirichlet and Neumann boundary conditions \eqref{eq:StrongFormMicro-2}$_{2,3}$ we obtain by \eqref{eq:AsympExp-u} and \eqref{eq:AsympExp-sigma}
\begin{eqnarray}
   u_i^{0} + \epsilon \, u_i^{1} + \epsilon^2  \, u_i^2 + \mathcal{O}(\epsilon^3) &=& \bar u_i  \qquad \mbox{on} \quad \partial \mathcal{B}_{D} 
      \label{eq:AsympExp-Dirichlet}            \\
  \left[ \dfrac{1}{\epsilon} \, \sigma_{ij}^{-1} + \sigma_{ij}^{0} + \epsilon \, \sigma_{ij}^{1} 
                                  + \mathcal{O}(\epsilon^2) \right] n_j &=& \bar t_i  \qquad \mbox{on} \quad \partial \mathcal{B}_{N} \, .
      \label{eq:AsympExp-Neumann}    
\end{eqnarray}

The boundary conditions for various orders of $\epsilon$ are identified from \eqref{eq:AsympExp-Dirichlet} and \eqref{eq:AsympExp-Neumann} according to 
\begin{equation}
\label{eq:AsympExp-Dirichlet-BCs}
\renewcommand{\arraystretch}{1.6}
\left.\begin{array}{rcl}
 \mathcal{O}\left(1 \right):  & & u_i^0 = \bar u_i   \\
 \mathcal{O}\left(\epsilon^{n}\right):  & & u_i^n = 0 \, , n \in \mathbb{N}    \\
\end{array}\right\} \;  \quad \mbox{on} \quad \partial \mathcal{B}_{D} \, ,
\end{equation}

\begin{equation}
\label{eq:AsympExp-Neumann-BCs}
\renewcommand{\arraystretch}{1.6}
\left.\begin{array}{rcl}
 \mathcal{O}\left(1 \right):  & & \sigma_{ij}^0 n_j = \bar t_i   \\
 \mathcal{O}\left(\epsilon^{n}\right):  & & \sigma_{ij}^n = 0 \, , n = -1,1,2, \ldots    \\
\end{array}\right\} \;  \quad \mbox{on} \quad \partial \mathcal{B}_{N} \, .
\end{equation}

It can be shown, that $\partial u^{0}_{i}/\partial y_j=0$, which implies $\bm u^0 = \bm u(\bm x)$ as already anticipated in \eqref{eq:AsympExp-u}, and moreover, that $\sigma_{ij}^{-1}= 0$. Thus, the first term in the expansion of $\bm u^{\epsilon}$, i.e. $\bm u^{0}$, is the coarse-scale/macroscale displacement field, and according to \eqref{eq:AsympExp-Dirichlet-BCs}$_1$ and \eqref{eq:AsympExp-Neumann-BCs}$_1$ the Dirichlet as well as Neumann boundary conditions \eqref{eq:AsympExp-Dirichlet}, \eqref{eq:AsympExp-Neumann} refer to the macroproblem. \\
For $\sigma_{ij}^{-1}=0$ it follows from the $\mathcal{O}\left(\epsilon^{-1}\right)$ equilibrium, \eqref{AsympExp-BalanceEq-2}, that
\begin{equation}
   \dfrac{\partial \sigma_{ij}^0}{\partial y_j} = 0 \, ,
\end{equation}
which yields along with the kinematical relations \eqref{eq:strain-relations-various-orders}$_2$ and the {\color{black} elastic constitutive law} \eqref{eq:stress-strain-relations-various-orders} for $s=0$ 
\begin{equation}
    \dfrac{\partial}{\partial y_j} \Big[ \mathbb{A}^{\epsilon}_{ijkh} \underbrace{\left( u^{0}_{(kh)x} + u^{1}_{(kh)y} \right)}_{\varepsilon^0_{kh}} \Big] = 0 \, .
    \label{eq:AsympExp-Macro-Equlibrium_0}
\end{equation}
A separation of variables $\bm x$ and $\bm y$ in \eqref{eq:AsympExp-Macro-Equlibrium_0} is achieved by the ansatz
\begin{equation}
    u^{1}_{k}(\bm x, \bm y) = \chi_k^{lm} (\bm y) \, u^{0}_{(lm)x} \, ,
    \label{eq:AsympExp-DisplInfluenceFunction}
\end{equation}
where $\chi_k^{lm}(\bm y)$ is the so-called first-order displacement influence function, which is locally periodic ($Y$-periodic) and symmetric $\chi_k^{lm}(\bm y) = \chi_k^{ml}(\bm y)$.
Inserting \eqref{eq:AsympExp-DisplInfluenceFunction} into \eqref{eq:AsympExp-Macro-Equlibrium_0} along with the identity $u^{0}_{(kh)x} = \mathbb{I}_{khlm} u^{0}_{(lm)x}$ where 
$\mathbb{I}_{khlm}= 1/2 \, (\delta_{kl} \delta_{hm} + \delta_{km} \delta_{hl})$ results in 
\begin{equation}
   \dfrac{\partial}{\partial y_j} \Big[ \mathbb{A}^{\epsilon}_{ijkh} \left( \mathbb{I}_{khlm} + \dfrac{\partial {\chi}_k^{lm}(\bm y)}{\partial y_h} \right) \Big] u^{0}_{(kh)x} = 0 \, ,
    \label{eq:AsympExp-Macro-Equlibrium_1}  
\end{equation}
which must hold for arbitrary macro strain $u^{0}_{(kh)x}$ leading to the cell problem, which reads: 

Find for any $l, m \in \{1, \ldots, d \}$ the functions ${\chi}_k^{lm}(\bm y)$ as the solutions of the system 
\begin{equation}
\label{eq:cell-problem}
\renewcommand{\arraystretch}{1.6}
\left.\begin{array}{rcl}
 && - \dfrac{\partial}{\partial y_j} \left(\mathbb{A}^{\epsilon}_{ijkh} \dfrac{\partial {\chi}_k^{lm}(\bm y)}{\partial y_h} \right) 
    = \dfrac{\partial \mathbb{A}^{\epsilon}_{ijlm}}{\partial y_j} \qquad \mbox{in} \,\, Y, \,\, i=1, \ldots, d \, , \\ 
 && \quad {\chi}_k^{lm} \, \mbox{is $Y$-periodic} \, .    
\end{array}\right\} \; 
\end{equation}
 
Since the field ${\chi}_k^{lm}$ enters the cell problem by $\nabla_y \bm \chi^{lm}$, it is determined up to a constant vector.
For definiteness, however, uniqueness can be achieved by imposing the normalization condition 
\begin{equation}
 \int_{Y} {\chi}_k^{lm}(\bm y) \, d \bm y=0 \, .
\label{eq:Normalization-for-chi}
\end{equation}
It can be shown, that if \eqref{eq:Normalization-for-chi} holds, 
the microstrain can be calculated from macrostrain by means of an elastic strain influence function $E_{kl}^{mn} (\bm y)$
\begin{equation}
   \varepsilon_{kl}^0 = E_{kl}^{mn}(\bm y) \, \bm u^0_{(mn)x} (\bm x) := \varepsilon_{kl} (\bm x, \bm y)
\end{equation}
with 
\begin{equation}
  \label{eq:Elastic-strain-influence-function}
  E_{kl}^{mn} (\bm y) = \mathbb{I}_{klmn} + \dfrac{\partial \chi_k^{mn}(\bm y)}{\partial y_l}   
\end{equation}
thereby realizing a separation of variables similar to \eqref{eq:AsympExp-DisplInfluenceFunction}.

\bigskip

With an eye towards the representation of the microproblem in the discrete case of FE-HMM we prefer to rewrite the cell problem \eqref{eq:cell-problem} as
\begin{equation}
    \dfrac{\partial}{\partial y_j} \Big(\mathbb{A}^{\epsilon}_{ijkh} \, \varepsilon_{kh}({\bm \chi}^{lm} + \bm I^{lm}) \Big)   
    = 0 \label{eq:rewrite-cell-problem-3} \, .   
\end{equation}
In \eqref{eq:rewrite-cell-problem-3} $\bm I^{lm}$ is a vector-valued function given by 
\begin{equation}
  I_k^{lm}(y) = y_m \delta_{kl}, \quad \delta_{kl} \, \mbox{is the Kronecker symbol}.
\end{equation}
The weak form of the cell problem in the format of \eqref{eq:rewrite-cell-problem-3} reads 
\begin{equation}
  \int_Y \mathbb{A}^{\epsilon}_{ijkh}(\bm y) \, \varepsilon_{kh}({\bm \chi}^{lm} + \bm I^{lm}) : \varepsilon_{ij}(\bm w) \, d \bm y =  0\, , \qquad \bm w \in \mathcal{W}_{per}(Y)\, ,
  \label{eq:Weak-form-cell-problem}
\end{equation}
where $\mathcal{W}_{per}$ is given by 
\begin{equation}
 \mathcal{W}_{per}(Y) = \left\{ \bm w; \bm w \in H^1_{per}(Y)^d\, ; \,\, \int_{Y} w_i \, d \bm y = 0 \, , \,\, i=1, \ldots, d\right\}.
 \label{eq:asymptotic_homo_PBCs_def}
\end{equation}
Averaging the microscale elasticity law, i.e. \eqref{eq:stress-strain-relations-various-orders} for $s=0$, $\sigma^0_{ij}= \mathbb{A}^{\epsilon}_{ijkh} \left( u^{0}_{(kh)x} + u^{1}_{(kh)y} \right)$ with an account of 
\eqref{eq:AsympExp-DisplInfluenceFunction} leads to 
\begin{equation}
\label{eq:homogenized-elasticity-relation}
 \langle{\sigma}^0_{ij}\rangle =  \mathbb{A}^{0}_{ijkh} \, \varepsilon^0_{kh} 
\end{equation}
where $\langle{\sigma}^0_{ij}\rangle= 1/|Y| \int_{Y} \sigma^0_{ij}(\bm x, \bm y) \, d \bm y$ and where the homogenized tensor $\mathbb{A}^{0}(\bm x)=(\mathbb{A}^{0}_{ijkl})$ is identified as
\begin{equation}
  \mathbb{A}^{0}_{ijlm}  =  \dfrac{1}{|Y|} \int_Y  \mathbb{A}^{\epsilon}_{ijlm}(\bm y) + \mathbb{A}^{\epsilon}_{ijkh}(\bm y) \dfrac{\partial {\chi}_k^{lm}(\bm y)}{\partial y_h} \, d \bm y \, ,
  \label{eq:HomogenizedElasticityTensor}
\end{equation}
which equally verifies the above symmetries of elasticity.

Integration the second term in \eqref{AsympExp-BalanceEq-3} over the unit cell domain yields
\begin{equation}
\label{eq:VolumeAverageDSigma1Dy}
   \dfrac{1}{|Y|} \int_{Y} \dfrac{\partial \sigma_{ij}^1}{\partial y_j} d \bm y = \dfrac{1}{|Y|} \int_{\partial \mathcal{B}} \sigma_{ij}^1 \, n_j \, d A = 0 \, .
\end{equation}
With \eqref{eq:VolumeAverageDSigma1Dy}, the boundary conditions \eqref{eq:AsympExp-Dirichlet-BCs}$_1$ and \eqref{eq:AsympExp-Neumann-BCs}$_1$, integration of the $\mathcal{O}\left(1\right)$ balance equation \eqref{AsympExp-BalanceEq-3} over the unit cell results in the the macroscopic BVP 
\begin{equation}
\label{eq:Homogenized-Strong-Form-derived}
\renewcommand{\arraystretch}{1.6}
\left.\begin{array}{rcl}
 - \dfrac{\partial}{\partial x_j} \left(\langle{\sigma}^0_{ij}\rangle \right) &=& \langle f_i \rangle  \qquad \mbox{in} \quad \mathcal{B}  \\
  u_i^{0} &=& \langle \bar u_i\rangle_{\Gamma} \qquad \mbox{on} \quad \partial \mathcal{B}_{D}    \\
  \langle {\sigma}^0_{ij} \rangle \, n_j &=& \langle \bar{t}_i \rangle_{\Gamma}   \qquad \mbox{on} \quad \partial \mathcal{B}_{N} \\
\end{array}\right\} \; .
\end{equation}
In \eqref{eq:Homogenized-Strong-Form-derived}$_{1}$  $\langle f_i \rangle$ is defined as the volume average of $f_i$, $\langle{\bm f(\bm x)}\rangle= 1/|Y| \int_{Y} \bm f(\bm x, \bm y) \, d \bm y$. In \eqref{eq:Homogenized-Strong-Form-derived}$_{2,3}$ $\langle \bullet \rangle_{\Gamma}$ defines a surface average according to 
\begin{equation}
\label{eq:SurfaceAveraging}
\langle \bar u_i\rangle_{\Gamma}:= \dfrac{1}{|\Gamma_D|} \int_{\Gamma_D} \bar{u}_i \, dA \, , \qquad 
\langle \bar{t}_i\rangle_{\Gamma}:= \dfrac{1}{|\Gamma_N|} \int_{\Gamma_N} \bar{t}_i \, dA \, , 
\end{equation}
where $\Gamma$ is that particular surface area ($\Gamma \in \partial \mathcal{B}$) to which the averaging is applied. 
Hence, $\langle \bullet \rangle_{\Gamma}$ is the surface counterpart to the volume average $\langle {\sigma}^0_{ij}\rangle$.
 
\bigskip

{\bf Remarks} \quad 
Three remarks are in order. First, that $\bm I^{lm}$ is to be understood as a displacement vector which induces by its definition a unitary strain $\bm \varepsilon(\bm I^{lm})$ with components 
$\varepsilon_{ij}^{lm}=1/2 (\partial I_i^{lm}/\partial y_j + \partial I_j^{lm}/\partial y_i)$. The corresponding strain vectors in Voigt notation read
\begin{center}
\begin{tabular}{c c c}
 $\bm \varepsilon^{11} = (1,0,0,0,0,0)^T$ & $\bm \varepsilon^{22} = (0,1,0,0,0,0)^T$  & $\bm \varepsilon^{33} = (0,0,1,0,0,0)^T$ \\
 $\bm \varepsilon^{12} = (0,0,0,1,0,0)^T$ & $\bm \varepsilon^{13} = (0,0,0,0,1,0)^T$  & $\bm \varepsilon^{23} = (0,0,0,0,0,1)^T$ \\
\end{tabular}
\end{center}
where $\bm \varepsilon^{lm}=(\varepsilon^{lm}_{11}, \varepsilon^{lm}_{22}, \varepsilon^{lm}_{33}, 2\varepsilon^{lm}_{12}, 2\varepsilon^{lm}_{13}, 2\varepsilon^{lm}_{23})^T$,  \cite{Lukkassen1995}.
In analogy to the homogenization of a second order PDE for a scalar-valued field problem, \cite{Assyr2006}, {\color{black}$(\bm \varepsilon^{lm})_{l,m=1}^d$ is the canonical basis of the symmetric, six-dimensional strain space.}

Second, introducing with $\bm U^{lm}:= \bm \chi^{lm} + \bm I^{lm}$ the superposition of the $Y$-periodic fluctuations with the (prescribed) displacement field that induces a unitary, homogeneous strain on the cell, it follows
for \eqref{eq:HomogenizedElasticityTensor}
\begin{eqnarray}
    \mathbb{A}^{0}_{ijlm} &=& \dfrac{1}{|Y|} \int_Y  \mathbb{A}^{\epsilon}_{ijlm} + \mathbb{A}^{\epsilon}_{ijkh}  \dfrac{\partial {\chi}_k^{lm}}{\partial y_h} \, d \bm y \\
                 &=& \dfrac{1}{|Y|} \int_Y  \mathbb{A}^{\epsilon}_{ijlm} + \mathbb{A}^{\epsilon}_{ijkh}  \, \varepsilon_{kh} (\bm \chi^{lm}) \, d \bm y \\ 
                 &=& \dfrac{1}{|Y|} \int_Y  \mathbb{A}^{\epsilon}_{ijlm} + \mathbb{A}^{\epsilon}_{ijkh}  \, \varepsilon_{kh} (\bm U^{lm} - \bm I^{lm}) \, d \bm y \\ 
                 &=& \dfrac{1}{|Y|} \int_Y  \mathbb{A}^{\epsilon}_{ijkh} \, \varepsilon_{kh} (\bm U^{lm}) \, d \bm y = \dfrac{1}{|Y|} \int_Y \sigma^0_{ij} (\bm U^{lm}) \, d \bm y \,. 
\end{eqnarray}
Hence, component $\mathbb{A}^{0}_{ijlm}$ of the homogenized elasticity tensor equals the averaged micro stress for the solution $\bm U^{lm}$ of the cell problem. 
\\
\\
Third, a note on notation. As indicated by its three indices $\chi_k^{lm}$ is a third order tensor; in \eqref{eq:AsympExp-Macro-Equlibrium_1} for example, the partial derivative of 
$\chi_k^{lm}$ with respect to $y_h$ yields a fourth order tensor, which is in dimensional consistency with the fourth order unity tensor in \eqref{eq:AsympExp-Macro-Equlibrium_1}. 
In \eqref{eq:rewrite-cell-problem-3} however, ${\bm \chi}^{lm}$ as well as $\bm I^{lm}$ are to be understood as vector-valued functions, where the superscripts $lm$ are frozen, and only one index is active. Then the calculation of strain from a displacement field as e.g. in \eqref{eq:rewrite-cell-problem-3} is a reasonable operation. Hence, the change in notation shall reflect this characteristics.

\subsection{General convergence} 
From homogenization theory the following convergence results and estimates are known.
\begin{enumerate}
 \item The microdisplacements $\bm u^{\epsilon}$, solution of \eqref{eq:StrongFormMicro-2}, converge for $\epsilon \rightarrow 0$ --usually in a weak sense--
       to $\bm u^0$, the solution of the homogenized problem \eqref{eq:VariationalFormHomogenizedProblem}, see \cite{Cioranescu-Donato-BOOK-1999}, Sec.10. 
       For the corresponding strong error estimate it holds
       \begin{equation}
         \label{eq:estimate_convergence_u_eps-against-u_0}
             || \bm  u^{\epsilon} - \bm u^0 ||_{L^2(\mathcal{B})} \leq C\epsilon \, .
       \end{equation}

 \item  Due to the oscillations of the fine scale solution, strong error estimates 
       in the $H^1$-norm can usually not be obtained since the gradients of the oscillations are in general not $\mathcal{O}(\epsilon)$ quantities.
       The homogenized solution requires a correction through information of the fine scale. This can be done defining the corrector given by 
       \begin{equation}
             \label{eq:definition_u_1-by-corrector}
             \bm  u^1(\bm x, \bm x/\epsilon) = \sum_{j=1}^d \omega^l(\bm x, \bm x/\epsilon) \dfrac{\partial \bm u^0(x)}{\partial x_j}\, ,
       \end{equation}   
       where $\omega^l$ is a $d \times d$ matrix given by $\omega^l = (\omega^l_{km})_{1\leq k,m \leq d} = \left( \chi^{lm}_{k} \right)_{1\leq k,m \leq d}$, 
       and where $\bm \chi^{lm}(\bm x, \bm x/\epsilon)$ is the solution of the cell problem. 
       With this correction an error estimate in the $H^1$-norm can be obtained
       \begin{equation}
             \label{eq:estimate_convergence_u_eps-against-u_01-in-H1-norm}
             || \bm  u^{\epsilon} - (\bm u^0 + \epsilon \, \bm u^1(\bm x, \bm x/\epsilon))||_{H^1(\mathcal{B})} \leq C \sqrt{\epsilon} \, ,
       \end{equation}
       where a boundary layer term yields a $\sqrt{\epsilon}$ instead of a $\varepsilon$ convergence rate, for a thorough analysis see \cite{Cioranescu-Donato-BOOK-1999}, Sec. 4.1.  
\end{enumerate}

\vfill
\newpage

\section{Numerical Data for Simulation Examples}
\label{sect:NumDataSimExamples}


\begin{table}[htbp]
\center
\renewcommand{\arraystretch}{1.2}
\begin{tabular}{llrrrrr}
\hline \\[-5mm]         
               & macroelements      &  $50 \times 10$  & $100 \times 20$  & $200 \times 40$  &  $400 \times 80$ &  $800 \times 160$  \\
$\epsilon$  & microelements      & \multicolumn{5}{c}{$32 \times 32$ }  \\[1mm] 
\hline \\[-5mm]                       
$5$      &  $|| \bm u^H ||_{\infty}$                      &  11.7997   &  11.8477   &  11.8609   &  11.8646   &   11.8658     \\
         &  $|| \bm u^H ||_A$                             &  1080.26   &  1082.41   &  1082.98   &  1083.13   &   1083.17     \\
         &  $|| \bm u^0 - \bm u^H ||_{L^2}$               &  703.5393  & 189.9063   &  52.2047   &  14.5294   &    3.9509     \\
         &  $|| \bm u^0 - \bm u^H ||_{H^1}$ in $10^{-3}$  & 179.1201   &  89.1681   &  48.5004   &  26.8540   &   13.9814     \\ 
         &  $|| \bm u^0 - \bm u^H ||_{A}$                 &  34.6706   &  17.3223   &   9.4171   &   5.2122   &    2.7175     \\[1mm]         
$0.5$    &  $|| \bm u ||_{\infty}$  &   11.7997 &  11.8477  &  11.8609  &  11.8646  &  11.8658   \\
         &  $|| \bm u ||_A$         &   1080.26 &  1082.41  &  1082.98  &  1083.13  &  1083.17   \\
\hline \\[-5mm]
$0.05$   &  $|| \bm u ||_{\infty}$  &   11.7997 &  11.8477  &  11.8609  &  11.8646  &  11.8658   \\ 
         &  $|| \bm u ||_A$         &   1080.26 &  1082.41  &  1082.98  &  1083.13  &  1083.17   \\[1mm]
\hline
\end{tabular}
\caption{{\bf Matrix-inclusion problem for beam bending.} Macroconvergence of FE-HMM
keeping the microdiscretization fixed, variation of $\epsilon$.  
\label{tab:Matrix-inclusion-2d-FEHMM-Macroconvergence}}
\end{table}

\vspace*{0cm}

\begin{table}[htbp]
\center
\renewcommand{\arraystretch}{1.2}
\begin{tabular}{llrrrrr}
\hline \\[-5mm]         
               & macroelements      &  \multicolumn{5}{c}{$200 \times 40$}       \\
$\epsilon$  & microelements      &  $16 \times 16$    & $32 \times 32$  & $64 \times 64$  &  $128 \times 128$ &  $256 \times 256$   \\[1mm] 
\hline \\[-5mm]                       
$5$      &  $|| \bm u^H ||_{\infty}$        &  11.8545   &  11.8609  &   11.8630  & 11.8636  & 11.8638 \\
         &  $|| \bm u^H ||_A$               &  1082.68   &  1082.98  &   1083.07  & 1083.10  & 1083.11 \\
         &  $|| \bm u^0 - \bm u^H ||_{L^2}$ &   102.4091 &  32.6310  &   10.0525  &  3.0164  & 0.8806  \\
         &  $|| \bm u^0 - \bm u^H ||_{H^1}$ in $10^{-4}$  & 48.0344  &   15.2866  &  4.6909  & 1.3997  &  0.4056  \\[1mm]   
\hline \\[-5mm]
$0.5$    &  $|| \bm u ||_{\infty}$  &   11.8545  &  11.8609   &  11.8630  &   11.8636  &  11.8638  \\
         &  $|| \bm u ||_A$         &   1082.68  &  1082.98   &  1083.07  &   1083.10  &  1083.11  \\
\hline \\[-5mm]
$0.05$   &  $|| \bm u ||_{\infty}$  &   11.8545  &   11.8609  &   11.8630 &   11.8636  &  11.8638  \\
         &  $|| \bm u ||_A$         &   1082.68  &  1082.98   &  1083.07  &   1083.10  &  1083.11  \\[1mm]
\hline
\end{tabular}
\caption{{\bf Matrix-inclusion problem for beam bending.} Microconvergence of FE-HMM
keeping the macrodiscretization fixed, variation of $\epsilon$. 
\label{tab:Matrix-inclusion-2d-FEHMM-Microconvergence}}
\end{table}


\begin{table}[htbp]
\center
\renewcommand{\arraystretch}{1.2}
\begin{tabular}{llrrrrr}
\hline \\[-5mm]         
${A}_{ij}^{0,h}$  &  \multicolumn{5}{c}{microelements}       \\
$ij$     &  $16 \times 16$    & $32 \times 32$  & $64 \times 64$  &  $128 \times 128$ &  $256 \times 256$   \\[1mm] 
\hline \\[-5mm]                       
$11$     &   46721.57  &   46700.78   &  46694.06  &   46691.95  &  46691.31  \\
$12$     &   11662.05  &   11666.10   &  11667.41  &   11667.79  &  11667.90  \\
$22$     &   46721.57  &   46700.78   &  46694.06  &   46691.95  &  46691.31  \\
$33$     &   17443.96  &   17436.62   &  17434.23  &   17433.48  &  17433.55  \\[1mm]         
\hline
\end{tabular}
\caption{{\bf Matrix-inclusion problem for beam bending.} Convergence of $\mathbb{A}^{0,h}$ components. 
\label{tab:Matrix-inclusion-2d-FEHMM-Convergence-of-elasticityTensor}}
\end{table}


\begin{table}[htbp]
\center
\renewcommand{\arraystretch}{1.2}
\begin{tabular}{llrrrrr}
\hline \\[-5mm]         
               & macroelements      &  $20 \times 20$  & $40 \times 40$  & $80 \times 80$  &  $160 \times 160$ &  $320 \times 320$  \\
$\epsilon$  & microelements      & \multicolumn{5}{c}{$40 \times 40$ }  \\[1mm] 
\hline \\[-5mm]                       
$0.005$  &  $|| \bm u^H ||_{\infty}$ in $10^{-3}$  &  77.8920  &  78.4353  &  78.6774  &  78.7895  &  78.8427    \\
         &  $|| \bm u^H ||_A$ in $10^{-3} $         &  81.7244  &  81.8499  &  81.8888  &  81.9011  &  81.9050    \\
         &  $|| \bm u^0 - \bm u^H ||_{L^2}$ in $10^{-6}$  &  47.9037    & 15.2437   &  4.9547  &  1.6312  &  0.5347  \\
         &  $|| \bm u^0 - \bm u^H ||_{H^1}$ in $10^{-5}$  &  38.2124    & 20.8187   & 11.7254  &  6.6562  &  3.7066  \\ 
         &  $|| \bm u^0 - \bm u^H ||_{A}$ in $10^{-4}$    &  38.0001    & 20.7804   & 11.7107  &  6.6520  &  3.7103  \\[1mm]         
\hline \\[-5mm]
$0.0005$ &  $|| \bm u ||_{\infty}$ in $10^{-3}$  &  77.8920  &  78.4353 &  78.6774  &  78.7895   &   78.8427  \\
         &  $|| \bm u ||_A$ in $10^{-3}$         &  81.7244  & 81.8499  &  81.8888  &  81.9011   &   81.9050  \\
\hline \\[-5mm]
$0.00005$ &  $|| \bm u ||_{\infty}$ in $10^{-3}$ &  77.8920  & 78.4353  &  78.6774  &  78.7895   &   78.8427   \\ 
          &  $|| \bm u ||_A$ in $10^{-3}$        &   81.7244 &  81.8499 &  81.8888  &  81.9011   &   81.9050 
  \\[1mm]
\hline
\end{tabular}
\caption{{\bf Microstructure with an analytical solution of homogenization.} Macroconvergence of FE-HMM keeping the microdiscretization fixed, variation of $\epsilon$.  
\label{tab:AnalyticalTensor-2d-FEHMM-Macroconvergence}}
\end{table}

\begin{table}[htbp]
\center
\renewcommand{\arraystretch}{1.2}
\begin{tabular}{llrrrrr}
\hline \\[-5mm]         
               & macroelements      &  \multicolumn{5}{c}{$40 \times 40$}       \\
$\epsilon$  & microelements      &  $16 \times 16$    & $32 \times 32$  & $64 \times 64$  &  $128 \times 128$ &  $256 \times 256$   \\[1mm] 
\hline \\[-5mm]                       
$0.005$      &  $|| \bm u^H ||_{\infty}$ in $10^{-3}$        &  78.3328   &  78.4353  &   78.4611  & 78.4675  & 78.4691 \\
             &  $|| \bm u^H ||_A$ in $10^{-3}$               &  81.7999   &  81.8499  &   81.8624  & 81.8656  & 81.8664 \\
             &  $|| \bm u^0 - \bm u^H ||_{L^2}$ in $10^{-7}$ &  199.8699  &  50.1244  &   12.5275  & 3.1182   & 0.7653  \\
             &  $|| \bm u^0 - \bm u^H ||_{H^1}$ in $10^{-7}$ &  267.2828  &  67.0309  &   16.7529  & 4.1700   & 1.0234  \\ 
             &  $|| \bm u^0 - \bm u^H ||_{A}$ in $10^{-7}$   &  422.0078  & 105.7584  &   26.4273  & 6.5777   & 1.6142  \\[1mm]         
\hline \\[-5mm]
$0.0005$    &  $|| \bm u ||_{\infty}$ in $10^{-3}$ &   78.3328  &  78.4353   &  78.4611  &   78.4675  &  78.4691  \\
            &  $|| \bm u ||_A$        in $10^{-3}$ &   81.7999  &  81.8499   &  81.8624  &   81.8656  &  81.8664  \\
\hline \\[-5mm]
$0.00005$   &  $|| \bm u ||_{\infty}$ in $10^{-3}$ &   11.8545  &  78.4353   &  78.4611 &   78.4675  &  78.4691  \\
            &  $|| \bm u ||_A$        in $10^{-3}$ &   81.7999  &  81.8499   &  81.8624  &   81.8656  &  81.8664  \\[1mm]
\hline
\end{tabular}
\caption{{\bf Microstructure with an analytical solution of homogenization.} Microconvergence of FE-HMM keeping the macrodiscretization fixed, variation of $\epsilon$. 
\label{tab:AnalyticalTensor-2d-FEHMM-Microconvergence}}
\end{table}

\begin{table}[htbp]
\center
\renewcommand{\arraystretch}{1.2}
\begin{tabular}{llrrrrr}
\hline \\[-5mm]         
                 &  \multicolumn{5}{c}{microelements}       \\
${A}_{ij}^{0,h}$   &  $20 \times 20$    & $40 \times 40$  & $80 \times 80$  &  $160 \times 160$ &  $320 \times 320$   \\[1mm] 
\hline \\[-5mm]                       
$11$     &    100.235     &  100.059           & 100.015   &  100.004        &  100.001        \\
$22$     &   140.028      &  140.028           & 140.028   &  140.028        &  140.028 
\\[1mm] 
\hline
\end{tabular}
\caption{{\bf Microstructure with an analytical solution of homogenization.} Convergence of the coefficients of $\mathbb{A}^{0,h}$. 
\label{tab:AnalyticalTensor-2d-FEHMM-Convergence-of-elasticityTensor}}
\end{table}

\begin{table}[htbp]
\center
\renewcommand{\arraystretch}{1.4}
\begin{tabular}{llccccc}
\hline
       &  macroelements                          			   &  $16\times 16$   & $64 \times 64$  & $144\times 144$  &  $256\times 256$ &  $576\times 576$ \\
$H^1$  &  microelements                          			   &  $4 \times 4$    &  $8 \times 8$   & $12\times 12$    &  $16\times 16$   &  $24\times 24$   \\
\hline
       & $|| \bm u^0 - \bm u^H ||_{L^2(\Omega)}$ in $10^{-5}$  &  46.1045         &  12.8707        &  5.6976          &   3.1917         &   1.4097    \\
       & $|| \bm u^0 - \bm u^H ||_{H^1(\Omega)}$ in $10^{-5}$  &  76.6287         &  22.2750        & 10.5488          &   6.2132         &   2.8595  \\
\hline
$L^2$  & microelements                           			   &  $16\times 16$   & $64 \times 64$  & $144\times 144$  &  $256\times 256$ &  $576\times 576$ \\
       & $|| \bm u^0 - \bm u^H ||_{L^2(\Omega)}$ in $10^{-7}$  &  994.5764        &  89.4812        &    22.7339       &   8.5862         &  1.9925  \\
       & $|| \bm u^0 - \bm u^H ||_{H^1(\Omega)}$ in $10^{-5}$  &   47.9858        &  14.1055        &     7.2641       &   4.4981         &  2.1436  \\
\hline 
\end{tabular}
\newline 
\caption{{\bf Optimal refinement strategies.} Comparison of the $L^2$ and $H^1$ errors for the FE-HMM.
\label{tab:optimal-convergence-L2-H1-errors}}
\end{table}

\begin{table}[htbp]
\center
\renewcommand{\arraystretch}{1.2}
\begin{tabular}{llrrrrr}
\hline \\[-5mm]         
               & macroelements      &  $20 \times 20$  & $40 \times 40$  & $80 \times 80$  &  $160 \times 160$ &  $320 \times 320$  \\
$\epsilon$  & microelements      & \multicolumn{5}{c}{$40 \times 40$ }  \\[1mm] 
\hline \\[-5mm]                       
$0.005$  &  $|| \bm u^H ||_{\infty}$ in $10^{-3}$         &  20.0360  & 20.2095 &  20.2925 &  20.3330  & 20.3530  \\
         &  $|| \bm u^H ||_A$  in  $10^{-3} $             &   4.1572  &  4.1630   &  4.1647   &  4.1653   &  4.1654       \\
         &  $|| \bm u^0 - \bm u^H ||_{L^2}$ in $10^{-7}$  & 104.7913  &  32.6314 & 10.3133  &  3.2907  & 1.0449   \\
         &  $|| \bm u^0 - \bm u^H ||_{H^1}$ in $10^{-6}$  & 108.7447  &  58.0419 & 31.7047  & 17.4515  & 9.4527    \\ 
         &  $|| \bm u^0 - \bm u^H ||_{A}$ in $10^{-5}$    &  18.5865  &   9.9356 &  5.4413  &  3.0054  & 1.6343    \\[1mm]         
\hline \\[-5mm]
$0.0005$ &  $|| \bm u ||_{\infty}$ in $10^{-3}$  &  20.0360 & 20.2095  & 20.2925 &  20.3330  & 20.3530   \\
         &  $|| \bm u ||_A$ in $10^{-3}$         &  4.1572  & 4.1630   & 4.1647  &  4.1653   &  4.1654   \\
\hline \\[-5mm]
$0.00005$ &  $|| \bm u ||_{\infty}$ in  $10^{-3}$ &   20.0360 & 20.2095  & 20.2925 &  20.3330 &  20.3530  \\ 
          &  $|| \bm u ||_A$ in $10^{-3}$         &   4.1572  &  4.1630  &  4.1647   &  4.1653    &    4.1654   \\[1mm]
\hline
\end{tabular}
\caption{{\bf Non-uniformly periodic tensor.} Macroconvergence of FE-HMM keeping the microdiscretization fixed, variation of $\epsilon$.  
\label{tab:Non-uniformly periodic tensor-2d-FEHMM-Macroconvergence}}
\end{table}

\begin{table}[htbp]
\center
\renewcommand{\arraystretch}{1.2}
\begin{tabular}{llrrrrr}
\hline \\[-5mm]         
               & macroelements      &  \multicolumn{5}{c}{$40 \times 40$ }  \\  
$\epsilon$  & microelements      & $20 \times 20$  & $40 \times 40$  & $80 \times 80$  &  $160 \times 160$ &  $320 \times 320$  \\[1mm] 
\hline \\[-5mm]                       
$0.005$   &  $|| \bm u^H ||_{\infty}$ in $10^{-3}$  &  20.1980         & 20.2095   &   20.2123  &   20.2131  &  20.2133 \\
          &  $|| \bm u^H ||_A$  in $10^{-3} $        &   4.16176        &  4.16295  &    4.16325 &   4.16333  & 4.16335  \\
         &  $|| \bm u^0 - \bm u^H ||_{L^2}$ in $10^{-8}$  &  235.0888   & 59.0407   &   14.7612  &   3.6745   & 0.9018   \\
         &  $|| \bm u^0 - \bm u^H ||_{H^1}$ in $10^{-8}$  &  290.3249   & 72.9140   &   18.2298  &   4.5379   & 1.1137   \\ 
         &  $|| \bm u^0 - \bm u^H ||_{A}$ in $10^{-8}$    &  137.4293   & 35.1447   &    8.8671  &   2.2173   & 0.5454   \\[1mm]         
\hline \\[-5mm]
$0.0005$ &  $|| \bm u ||_{\infty}$ in $10^{-3}$           &  20.1980    &  20.2095  & 20.2123   & 20.2131    &   20.2133  \\
         &  $|| \bm u ||_A$ in  $10^{-3}$                  &  4.16176    &  4.16295  & 4.16325   &  4.16333   &   4.16335  \\
\hline \\[-5mm]
$0.00005$ &  $|| \bm u ||_{\infty}$ in $10^{-3}$ &   20.1980 &  20.2095  & 20.2123   &  20.2131  &  20.2133    \\ 
          &  $|| \bm u ||_A$ in $10^{-3}$        &   4.16176 &   4.16295 &  4.16325  &   4.16333 &   4.16335   \\[1mm]
\hline
\end{tabular}
\caption{{\bf Non-uniformly periodic tensor.} Microconvergence of FE-HMM keeping the macrodiscretization fixed, variation of $\epsilon$. 
\label{tab:Non-uniformly periodic tensor-Microconvergence}}
\end{table}
 
\begin{table}[htbp]
\center
\renewcommand{\arraystretch}{1.2}
\begin{tabular}{llrrrrr}
\hline \\[-5mm]         
${A}_{ij}^{0,h}$  &  \multicolumn{5}{c}{microelements}       \\
 $ij$                    &  $20 \times 20$ & $40 \times 40$  & $80 \times 80$  &  $160 \times 160$ &  $320 \times 320$   \\[1mm] 
\hline \\[-5mm]                       
11     &   3.74326      &  3.74018  &  3.73940  &   3.73921  & 3.73916 \\
12     &   0.94034      &  0.93979  &  0.93965  &   0.93961  & 0.93961 \\
22     &   3.99944      &  3.99767  &  3.99723  &  3.99712   & 3.99709 \\
33     &   1.40440      &  1.40325  &  1.40295  &  1.40288   & 1.40286 \\[1mm]         
\hline
\end{tabular}
\caption{
{\bf Non-uniformly periodic tensor.} Convergence of the coefficients of $\mathbb{A}^{0,h}$. 
\label{tab:Non-uniformly periodic tensor-2d-FEHMM-Convergence-of-elasticityTensor}
}
\end{table}

\begin{table}[htbp]
\center
\renewcommand{\arraystretch}{1.2}
\begin{tabular}{lcccccc}
\hline \\[-5mm]         
Macro DOFs        &  $1000$  & $4000$  & $16000$  &  $64000$ &  $256000$  \\
Micro DOFs / $h$  & \multicolumn{5}{c}{$1462$ / 1}  \\[1mm] 
\hline \\[-5mm]                       

$|| \bm u^0 - \bm u^H ||_{L_2}\cdot 10^{-2}$   &   0.2042  &  0.0554   &  0.0153    &   0.0043   &    0.0012     \\
$|| \bm u^0 - \bm u^H ||_{H_1}\cdot 10^{-6}$   &  0.5304   &   0.2665   &  0.1459    &   0.0811   &    0.0423     \\ 
$|| \bm u^0 - \bm u^H ||_{A}\cdot 10^{-2}$     &  5.2007    &  2.6193    &  1.4314    &   0.7948   &    0.4154     \\[1mm]         

\hline
\end{tabular}
\caption{{\bf Escher's periodic Bird and Fish tesselation.} Macro-convergence of FE-HMM keeping the microdiscretization fixed. \label{tab:Escher-birdfish-2d-FEHMM-Macroconvergence}}
\end{table}

\begin{table}[htbp]
\center
\renewcommand{\arraystretch}{1.2}
\begin{tabular}{lccccc}
\hline \\[-5mm]         
Macro DOFs      &  \multicolumn{4}{c}{$16000$ }  \\
Micro DOFs / $h$     &  1462/1  & 2486/0.75 & 5502/0.5 & 21670/0.25 &  82856/0.125 \\[1mm] 
\hline \\[-5mm]                       

$|| \bm u^0 - \bm u^H ||_{L_2}\cdot 10^{-2}$   &  0.1438   &  0.0969   &  0.0491   &  0.0150    &   0.0036        \\
$|| \bm u^0 - \bm u^H ||_{H_1}\cdot 10^{-7}$   &  0.6395   &  0.4349   &  0.2202   &  0.0677    &   0.0166        \\ 
$|| \bm u^0 - \bm u^H ||_{A}\cdot 10^{-2}$     &  0.8506   &  0.5755   &  0.2925   &  0.0900    &   0.0221        \\[1mm]         

\hline
\end{tabular}
\caption{{\bf Escher's periodic Bird and Fish tesselation.} Micro-convergence of FE-HMM keeping the macrodiscretization fixed. \label{tab:Escher-birdfish-2d-FEHMM-Microconvergence}}
\end{table}


\end{appendix}

\bibliographystyle{plainnat}

\begin{thebibliography}{99}

\expandafter\ifx\csname url\endcsname\relax\def\url#1{\texttt{#1}}\fi
 

\bibitem[Abdulle (2005)]{Assyr2005}
A. Abdulle, 
\newblock On a-priori error analysis of fully discrete Heterogeneous Multiscale FEM, 
\newblock {SIAM Multiscale Model. Simul.} {4,2} (2005) 447--459.
                  
\bibitem[Abdulle (2006)]{Assyr2006}
A. Abdulle, 
\newblock Analysis of the heterogeneous multiscale FEM for problems in elasticity, 
\newblock {Math. Models Methods Appl. Sci.} {16(4)} (2006) 615--635.
                               
\bibitem[Abdulle (2009)]{Assyr2009}
A. Abdulle,
\newblock The Finite Element Heterogeneous Multiscale Method: a computational strategy for multiscale PDEs,
\newblock {Math. Sci. Appl., Vol. 31} {31} (2009) 133--181.
 
\bibitem[Abdulle and Attinger (2006)]{AbdulleAttinger2006}
A. Abdulle, S. Attinger,
\newblock Numerical methods for transport problems in microdevices,
In: I. Lirkov, S. Margenov, and J. Wa\'{sa}sniewski (Eds.): Lecture Notes in Computer Science 3743, pp. 67--75, 2006.
Springer-Verlag Berlin Heidelberg.
                  
\bibitem[Abdulle and Bai (2013)]{AssyrBai2013}
A. Abdulle, Y. Bai, 
\newblock Adaptive reduced basis finite element heterogeneous multiscale method,
\newblock {Comput. Methods Appl. Mech. Engrg.} {257} (2013) 203--220.
 
\bibitem[Abdulle and Nonnenmacher (2009)]{AssyrNonnenmacher2009}
A. Abdulle, A. Nonnenmacher,
\newblock A short and versatile finite element multiscale code for homogenization problems,
\newblock {Comput. Methods Appl. Mech. Engrg.} {198} (2009) 2839--2859.
                                
\bibitem[Abdulle and Nonnenmacher (2011)]{AssyrNonnenmacher2011}
A. Abdulle, A. Nonnenmacher,
\newblock Adaptive finite element heterogeneous multiscale method for homogenization problems,
\newblock {Comput. Methods Appl. Mech. Engrg.} {200} (2011) 2710--2726.

\bibitem[Abdulle and Schwab (2005)]{AssyrSchwab2005}
A. Abdulle, C. Schwab,
\newblock Heterogeneous multiscale FEM for diffusion problems on rough surfaces,
\newblock {Multiscale Model. Simul.} {3(1)} (2005) 195--220.
 
\bibitem[Abdulle et al. (2012)]{Assyr-etal2012}
A. Abdulle, W. E., B. Engquist, E. Vanden-Eijnden, 
The heterogeneous multiscale method, 
\newblock {Acta Numer.} {466} (2012) 1--87.  

\bibitem[Allaire (1992)]{Allaire1992}
G. Allaire,
Homogenization and two-scale convergence. 
\newblock {SIAM J. Math. Anal} {23} (1992) 1482--1518.  
                               
\bibitem[Amelang et al. (2015)]{AmelangVenturiniKochmann2015}
J.S. Amelang, G.N. Venturini, D.M. Kochmann, 
\newblock Summation rules for a fully nonlocal energy-based quasicontinuum method,
\newblock {J. Mech. Phys. Solids} {82} (2015) 378--413.
                                        
\bibitem[Barlow (1976)]{Barlow1976}
J. Barlow,
Optimal stress locations in finite element models, 
\newblock {Int. J. Numer. Methods Eng.} {10} (1976) 243--251.
                                        
\bibitem[Bensoussan et al. (1976)]{Bensoussan-Lions-Papanicolau-BOOK-1976}
A. Bensoussan, J.L. Lions, G. Papanicolau,
Asymptotic Analysis for Periodic Structures,
\newblock  North-Holland, Amsterdam (1976). 
                                        
\bibitem[Bertsekas (1996)]{Bertsekas-BOOK-1996}
D.P. Bertsekas,
Constrained Optimization and Lagrange Multiplier Methods,
\newblock Athena Scientific, Belmont (1996). 
                  
\bibitem[Braess (2012)]{Braess-BOOK-1997}
D. Braess,
Finite Elements,
\newblock  Cambridge University Press, Cambridge, UK  (1997). 
                  
\bibitem[Ciarlet (1978)]{Ciarlet-BOOK-1978}
P.G. Ciarlet,
The Finite Elements Method for Elliptic Problems,
\newblock  North Holland (1978). 
                               
\bibitem[Cioranescu and Donato (1999)]{Cioranescu-Donato-BOOK-1999}
D. Cioranescu, P. Donato,
An Introduction to Homogenization,
\newblock Oxford University Press, New York (1999).

\bibitem[Duan et al. (2005)]{Duan-Wang-Huang-Karihaloo2005}
H.L. Duan, J. Wang, Z.P. Huang, and B.L. Karihaloo, 
\newblock Size-dependent effective elastic constants of solids containing nano-inhomogeneities with interface stress.
\newblock {J. Mech. Phys. Solids} {53(7)} (2005) 1574--1596. 
    
\bibitem[E and Engquist (2003)]{E-Engquist-2003}
W. E, B. Engquist,
The heterogeneous multi-scale methods,
\newblock {Commun. Math. Sci.} {1} (2003) 87--132.
 
\bibitem[E, Engquist, and Huang (2003)]{E-Engquist-Huang-2003}
W. E, B. Engquist, Z. Huang,
Heterogeneous multiscale method: A general methodology for multiscale modeling,
\newblock {Phys. Rev. B: Condens. Matter} {67} (2003) 092101. 
                  
\bibitem[E et al. (2007)]{E-Engquist-Li-Ren-VandenEijnden2007}
W. E, B. Engquist, X. Li, W. Ren, E. Vanden-Eijnden,
\newblock Heterogeneous Multiscale Methods: A Review,
\newblock {Commun. Comput. Phys} {2} (2007) 367--450.
                  
\bibitem[E et al. (2005)]{E-Ming-Zhang-2005}
W. E, P. Ming, P. Zhang,
\newblock Analysis of the heterogeneous multi-scale method for elliptic homogenization problems,
\newblock {J. Amer. Math. Soc} {18} (2005) 121--156.
 
\bibitem[Eidel and Stukowski (2009)]{EidelStukowski2009}
B. Eidel, A. Stukowski,
\newblock  A variational formulation of the quasicontinuum method based on energy sampling in clusters,
\newblock {J. Mech. Phys. Solids} {57} (2009) 87--108.
                  
\bibitem[Eidel (2009)]{Eidel2009}
B. Eidel,
\newblock  Coupling atomistic accuracy with continuum effectivity for predictive simulations in materials research 
- the Quasicontinuum method,
\newblock {Int. J. Mater. Res.} {100} (2009) 1503--1512.
                  
\bibitem[Eidel et al. (2010)]{Eidel-etal-2010}
B. Eidel, A. Hartmaier, P. Gumbsch, 
Atomistic Simulation Methods and their Application on Fracture, in: R. Pippan, P. Gumbsch (Eds.), Multiscale Modelling of Plasticity and Fracture by Means of Dislocation Mechanics, CISM International Centre for Mechanical Sciences, (2010) 1--57. 
                  
\bibitem[Eidel and Fischer (2016)]{EidelFischer2016}
B. Eidel, A. Fischer,
\newblock  The heterogeneous multiscale finite element method FE-HMM for the homogenization of linear elastic solids,
\newblock {PAMM} {16} (2016) 521--522.

\bibitem[Feyel and Chaboche (2000)]{FeyelChaboche2000}
F. Feyel, J.L. Chaboche, 
\newblock FE2 multiscale approach for modelling the elastoviscoplastic behaviour of long fibre SiC/Ti composite materials,
\newblock {Comput. Methods Appl. Mech. Engrg.} {183} (2000) 309--330.
                               
\bibitem[Fish et al. (1999)]{Fish-etal-damage1999}
J. Fish, Q. Yu, K. Shek, 
\newblock Computational damage mechanics for composite materials based on mathematical homogenization,
\newblock {Int. J. Numer. Meth. Eng.} {45} (1999) 1657--1679.
                  
\bibitem[Geers et al. (2010a)]{GeersKouznetsovaBrekelmans2010a}
M.G.D. Geers, V.G. Kouznetsova, W.A.M. Brekelmans, 
\newblock Multi-scale computational homogenization: trends and challenges,
\newblock {J. Comput. Appl. Math.} {234} (2010) 2175--2182.
                  
\bibitem[Geers et al. (2010b)]{GeersKouznetsovaBrekelmans2010b}
M.G.D. Geers, V.G. Kouznetsova, W.A.M. Brekelmans, 
\newblock Computational homogenization, in: R. Pippan, P. Gumbsch (Eds.), Multiscale Modelling of Plasticity and Fracture by Means of Dislocation Mechanics, CISM International Centre for Mechanical Sciences, (2010) 327--394. 
                  
\bibitem[Guedes and Kikuchi (1990)]{GuedesKikuchi1990}
J.M. Guedes, N. Kikuchi,
\newblock Preprocessing and postprocessing for materials based on the homogenization method with adaptive
finite element methods,
\newblock {Comput. Methods Appl. Mech. Engrg.} {83} (1990) 143--198.
                               
\bibitem[Hill (1963)]{Hill1963}
R. Hill,
\newblock Elastic properties of reinforced solids: some theoretical principles,
\newblock {J. Mech. Phys. Solids} {11} (1963) 357--372.

\bibitem[Hill (1972)]{Hill1972}
R. Hill,
\newblock On constitutive macro-variables for heterogeneous solids at finite strain,
\newblock {Proc. R. Soc. London, Ser. A} {326} (1972) 131--147.
                  
\bibitem[Hou and Wu (1997)]{HouWu1997}
T. Hou, X. Wu,
\newblock A multiscale finite element method for elliptic problems in composite materials and porous media,
\newblock {J. Comput. Phys.} {134} (1997) 169--189.
                  
\bibitem[Hou et al. (1999)]{Hou-Wu-Cai1999}   
T. Hou, X. Wu, Z. Cai,
\newblock Convergence of a multiscale finite element method for elliptic problems with rapidly oscillating coefficients,
\newblock {Math. Comp.} {68} (1999) 913--943.
                    
\bibitem[Hughes (2000)]{Hughes-BOOK-2000}
T.J.R. Hughes,
The finite element method: linear static and dynamic finite element analysis,
\newblock  Dover Publications, Mineola, New York (2000). 
                  
\bibitem[Javili et al. (2013)]{JaviliChatzigeorgiouSteinmann2013}
A. Javili, Chatzigeorgiou, P. Steinmann,
\newblock Computational homogenization in magneto-mechanics,
\newblock {Int. J. Numer. Meth. Eng.} {50(25--26)} (2013) 4197--4216.

\bibitem[Javili et al. (2017)]{JaviliSteinmannMosler2017}
A. Javili, P. Steinmann, J. Mosler, 
\newblock Micro-to-macro transition accounting for general imperfect interfaces, 
\newblock {Comput. Methods Appl. Mech. Engrg.} {317} (2017) 274--317.

\bibitem[Keip et al. (2014)]{KeipSteinmannSchroeder-2014}
M.-A. Keip, P. Steinmann, J. Schr\"{o}der,  
\newblock Two-scale computational homogenization of electro-elasticity at finite strains, 
\newblock {Comput. Methods Appl. Mech. Engrg.} {278} (2014) 62--79. 
                  
\bibitem[Knap and Ortiz (2001)]{KnapOrtiz2001}
J. Knap,  M. Ortiz,
\newblock An analysis of the quasicontinuum method,
\newblock {J. Mech. Phys. Solids} {49} (2001) 1899--1923.
                  
\bibitem[Kouznetsova et al. (2001)]{Kouznetsova-etal2001}
V. Kouznetsova, W.A.M. Brekelmans, F.P.T. Baaijens, 
\newblock An approach to micro-macro modeling of heterogeneous materials,
\newblock {Comput. Mech.} {27} (2002) 37--48.
                  
\bibitem[Kouznetsova et al. (2002)]{Kouznetsova-etal2002}
V. Kouznetsova, M.G.D. Geers, W.A.M. Brekelmans,
\newblock Multi-scale constitutive modelling of heterogeneous materials with a gradient-enhanced computational homogenization scheme,
\newblock {Int. J. Numer. Meth. Eng.} {54} (2002) 1235--1260.
 
\bibitem[Larsson et al. (2011)]{LarssonRunesson-etal-2011}
F. Larsson, K. Runesson, S. Saroukhani, R. Vafadari, 
\newblock Computational homogenization based on a weak format of micro-periodicity for RVE-problems, 
\newblock Comput. Methods Appl. Mech. Engrg. {1-4} (2011) 11--26. 
                  
\bibitem[Li and E (2005)]{LiE2005}
X. Li,  W. E, 
\newblock Multiscale modeling of the dynamics of solids at finite temperature,
\newblock {J. Mech. Phys. Solids} {53} (2005) 1650--1685.
                   
\bibitem[Lukkassen et al. (1995)]{Lukkassen1995}
D. Lukkassen, L.-E. Persson, P. Wall, 
\newblock Some engineering and mathematical aspects on the homogenization method,
\newblock {Compos. Eng.} {5} (1995) 519--531.
                   
\bibitem[Michel et al. (1999)]{MichelMoulinecSuquet1999}
J.C. Michel, H. Moulinec, P. Suquet, 
\newblock Effective properties of composite materials with periodic microstructure: a computational approach,
\newblock {Comput. Methods Appl. Mech. Engrg.} {172} (1999) 109--143.
                  
\bibitem[Miehe et al. (1999a)]{Miehe-etal-1999a}
C. Miehe,  J. Schr\"{o}der, J. Schotte, 
\newblock Computational homogenization analysis in finite plasticity Simulation of texture development in polycrystalline materials,
\newblock {Comput. Methods Appl. Mech. Engrg.} {171} (1999) 387--418.
                  
\bibitem[Miehe et al. (1999b)]{Miehe-etal-1999b}
C. Miehe,  J. Schotte, J. Schr\"{o}der,
\newblock Computational homogenization analysis in finite plasticity Simulation of texture development in polycrystalline materials,
\newblock {Comput. Mat. Sci.} {16 (1-4)} (1999) 372--382.

\bibitem[Miehe and Koch (2002)]{Miehe-Koch-2002}
C. Miehe,  A. Koch, 
\newblock Computational micro-to-macro transitions of discretized microstructures undergoing small strain,
\newblock {Arch. Appl. Mech.} {71} (2002) 300--317.
                  
\bibitem[Ming and Yue (2006)]{Ming-Yue-2006}
P. Ming, X. Yue, 
\newblock Numerical methods for multiscale elliptic problems,
\newblock {J. Comput. Phys.} {214} (2006) 421--445.
                   
\bibitem[Moulinec and Suquet (1994)]{Moulinec-Suquet-1994}
H. Moulinec, P. Suquet,
\newblock  Fast numerical method for computing the linear and nonlinear properties of composites,
\newblock {CR. Acad. Sci. II} {318} (1994) 1417--1423. 
        
\bibitem[Nonnenmacher (2011)]{Nonnenmacher-Disse-2011}
A. Nonnenmacher, 
\newblock Adaptive Finite Element Methods for Multiscale Partial Differential Equations,
\newblock {PhD Thesis} N$^{\circ}$ 5097, \'{E}cole Polytechnique F\'{e}d\'{e}rale de Lausanne (2011).
                   
\bibitem[Ohlberger (2005)]{Ohlberger2005}
M. Ohlberger,
\newblock A posteriori error estimates for the heterogeneous multiscale finite element method for elliptic homogenization problems, 
\newblock {Multiscale Model. Simul.} {4}(1) (2005) 88--114. 
                  
\bibitem[\"{O}zdemir et al. (2008)]{OezdemirBrekelmansGeers2008}
I. \"{O}zdemir, W.A.M. Brekelmans, M.G.D. Geers,
\newblock FE$^2$ computational homogenization for the thermo-mechanical analysis of heterogeneous solids, 
\newblock {Comput. Methods Appl. Mech. Engrg.} {198} (2008) 602--613. 

\bibitem[Peric et al. (2010)]{Peric-etal-2010}
D. Peri\'{c}, E.A. de Souza Neto, R.A.Feij\'{o}o, M. Partovi, A.J. Carneiro Molina, 
\newblock On micro-to-macro transitions for multi-scale analysis of non-linear heterogeneous materials: unified variational
basis and finite element implementation, 
\newblock {Int. J. Numer. Methods Eng.} {87} (2010) 149--170.                  

\bibitem[Saeb et ali (2016)]{Saeb-Steinmann-Javili-2016}
S. Saeb, P. Steinmann, A. Javili, 
Aspects of computational homogenization at finite deformations: a unifying review from Reuss' to Voigt's bound,
\newblock {Appl. Mech. Rev.} {68} (2016) 1--33. 
 
\bibitem[Sanchez-Palencia (1980)]{Sanchez-Palencia-BOOK-1980}
E. Sanchez-Palencia, 
Non-Homogeneous Media and Vibration Theory, Lecture Notes in Physics, Vol. 127,
\newblock {Springer, Berlin} (1980). 

\bibitem[Schneider et al. (2015)]{Schneider-Ospald-Kabel-2015}
M. Schneider, F. Ospald, M. Kabel, 
Computational homogenization of elasticity on a staggered grid,
\newblock {Int. J. Numer. Meth. Engng.} {105}(9) (2015) 693--720.    

\bibitem[Schr\"{o}der (2009)]{Schröder2009}
J. Schr\"{o}der, 
Derivation of the localization and homogenization conditions for electro-mechanically coupled problems,
\newblock {Comput. Mat. Sci.} {46}(3) (2009) 595--599. 
                  
\bibitem[Schr\"{o}der et al. (2010)]{SchröderBalzaniBrands2010}
J. Schr\"{o}der, D. Balzani, D. Brands,
Approximation of random microstructures by periodic statistically similar representative volume elements based on lineal-path functions,
\newblock {Arch. Appl. Mech.} {81}(7) (2010) 975--997. 
                
\bibitem[Schr\"{o}der (2014)]{Schröder2014}
J. Schr\"{o}der,  
A numerical two-scale homogenization scheme: the FE$^2$-method, 
in: J. Schr\"{o}der, K. Hackl (Eds.), Plasticity and Beyond, CISM International Centre for
Mechanical Sciences, (2014) 1--64. 

\bibitem[Temizer and Wriggers (2008)]{TemizerWriggers2008}
I. Temizer, P. Wriggers, 
On the computation of the macroscopic tangent for multiscale volumetric homogenization problems, 
\newblock {Comput. Methods Appl. Mech. Engrg.} {198}(3) (2008) 495--510. 

\bibitem[Yue and E (2007)]{Yue-E2007}
X.Y. Yue, W. E,
The local microscale problem in the multiscale modelling of strongly heterogeneous media: Effect of boundary conditions and cell size, 
\newblock {J. Comput. Phys.} {222} (2007) 556--572. 

\bibitem[Zienkiewicz and Zhu (1992a)]{SPR}
O.C. Zienkiewicz, J.Z. Zhu,
The superconvergent patch recovery and a posteriori error estimates. part 1: the recovery technique,
\newblock {Int. J. Numer. Methods Eng.} {33} (1992) 1331--1364.

\bibitem[Zienkiewicz and Zhu (1992b)]{SPR2}
O.C. Zienkiewicz, J.Z. Zhu,
The superconvergent patch recovery and a posteriori error estimates. part 2: error estimates and adaptivity,
\newblock {Int. J. Numer. Methods Eng.} {33} (1992) 1365--1382.
 
 
\end{thebibliography}
\end{document}